%% file: article.tex
\documentclass[preprint,3p,times]{elsarticle}

\usepackage[section]{placeins}
\usepackage{multirow}
\usepackage{makecell}
\usepackage{longtable}
\usepackage{booktabs}
\usepackage{supertabular}

\usepackage[bf]{caption2}

\usepackage{amssymb}
\usepackage{amsmath}
\newtheorem{theorem}{Theorem}
\newtheorem{lemma}{Lemma}
\newtheorem{property}{Property}

\newtheorem{corollary}{Corollary}

\newtheorem{remark}{Remark}
\newtheorem{example}{Example}

\journal{Elsevier}

\begin{document}

\begin{frontmatter}

  \title{A modified adaptive improved mapped WENO method}
  
  \author[a]{Ruo Li}
  \ead{rli@math.pku.edu.cn}

  \author[b,c]{Wei Zhong\corref{cor1}}
  \ead{zhongwei2016@pku.edu.cn}

  \cortext[cor1]{Corresponding author}

  \address[a]{CAPT, LMAM and School of Mathematical Sciences, Peking
    University, Beijing 100871, China}

  \address[b]{School of Mathematical Sciences, Peking University,
    Beijing 100871, China}

  \address[c]{Northwest Institute of Nuclear Technology, Xi'an 
  710024, China}

  \input{article_abstract}


\end{frontmatter}

\input{article_introduction}
\input{article_fvm}
\input{article_modify}
\input{article_numerics}
\input{article_conclusion}

\input{article_appendix}

\bibliographystyle{model1b-shortjournal-num-names}
\bibliography{../refs}

\end{document}

%% file: article_abstract.tex
\begin{abstract}

 We propose a new family of mapped WENO schemes by using several 
 adaptive control functions and a smoothing approximation of the 
 signum function. The proposed schemes introduce the adaptivity and 
 admit an extensive permitted range of the parameters in the mapping 
 functions. Consequently, they have the capacity to achieve optimal 
 convergence rates, even near critical points. Particularly, one of 
 these new schemes with fine-tuned parameters illustrates a 
 significant advantage when solving problems with discontinuities. 
 It produces numerical solutions with high resolution without 
 generating spurious oscillations, especially for long output times.


\end{abstract}


%% file: article_introduction.tex
\section{Introduction}
\label{secIntroduction}
In recent decades, the essentially non-oscillatory (ENO) schemes
\cite{ENO1987JCP71, ENO1987V24, ENO1986, ENO1987JCP83} and weighted
ENO (WENO) schemes \cite{ENO-Shu1988, ENO-Shu1989, WENO-LiuXD,
WENO-JS, WENOoverview} have been developed quite successfully to
solve the hyperbolic conservation laws in the form
\begin{equation}
    \begin{array}{ll}
      \dfrac{\partial \mathbf{u}}{\partial t} +
      \displaystyle\sum\limits_{\alpha = 1}^{d} \dfrac{\partial
      \mathbf{f}_{\alpha}(\mathbf{u})}{\partial x_{\alpha}} = 0,
      & x_{\alpha} \in \mathbb{R}, t > 0,
\end{array}
\label{governingEquation}
\end{equation}
where $\mathbf{u} = (u_{1}, u_{2}, \cdots, u_{m}) \in \mathbb{R}^{m}$
are the conserved variables and $\mathbf{f}_{\alpha}: \mathbb{R}^{m} 
\rightarrow \mathbb{R}^{m}$, $\alpha = 1,2,\cdots,d$ are the 
Cartesian components of flux.  Within the general framework 
(referred as WENO-JS below) of smoothness indicators and non-linear 
weights proposed by Jiang and Shu \cite{WENO-JS}, many successful 
works have improved the accuracy and efficiency of WENO schemes 
\cite{WENO-M, WENO-IM, WENO-PM, WENO-Z, WENO-Z01, WENO-Z02, 
WENO-Accuracy, WENO-GlobalAccuracy}.

Henrick et al. \cite{WENO-M} found that the fifth-order WENO-JS 
scheme is only third-order accurate at critical points of order 
$n_{\mathrm{cp}} = 1$ in the smooth regions, where $n_{\mathrm{cp}}$ 
denotes the order of the critical point; e.g., $n_{\mathrm{cp}} = 1$ 
corresponds to $f'=0, f'' \neq 0$ and $n_{\mathrm{cp}} = 2$ 
corresponds to $f'=0, f'' = 0, f''' \neq 0$, etc. They derived 
necessary and sufficient conditions \cite{WENO-M} on the weights for
optimality of the order. Then, by introducing a mapping function to 
the original weights of the WENO-JS scheme, they developed the 
WENO-M method \cite{WENO-M} to achieve the optimal order of accuracy 
in smooth regions even with critical points. Feng et al. 
pointed out that \cite{WENO-PM}, when the mapping function of the 
WENO-M scheme is used for solving the problems with discontinuities, 
it may amplify the effect from the non-smooth stencils, thereby 
causing a potential loss of accuracy near discontinuities. In order 
to address this issue, they devised the WENO-PM$k$ scheme 
\cite{WENO-PM} by proposing a piecewise polynomial mapping function 
with two additional requirements, that is, $g'(0)=0$ and $g'(1)=0$ 
($g(x)$ denotes the mapping function), to the original criteria in 
\cite{WENO-M}. Also, these two additional requirements were 
considered to be very important to decrease the effect from the 
non-smooth stencils \cite{WENO-RM260}, and they were used in the 
construction of the WENO-RM($mn0$) scheme \cite{WENO-RM260}. 
Similarly, requirements $g'(0) = 1$ and $g'(1) = 1$ were employed 
when the WENO-PPM$n$($n=4,5,6$) \cite{WENO-PPM5} schemes were 
constructed. Although the mapping function in the WENO-PM$k$ scheme 
decreases the effect from the non-smooth stencils, it is not smooth 
enough \cite{WENO-RM260} because it is only piecewise continuous. 
Furthermore, the WENO-PM$k$ scheme may generate the oscillations 
near discontinuities \cite{WENO-RM260}. Recently, Feng et al. 
\cite{WENO-IM} proposed a new family of mapping functions, that has 
two parameters $k$ and $A$, to improve the WENO-M method. They 
called the corresponding improved mapped WENO scheme WENO-IM($k$, $A$
). The WENO-IM($k$, $A$) scheme with proper parameters can achieve 
optimal order of accuracy near critical points for any 
$(2r - 1)$th-order WENO schemes. Moreover, the recommended version 
of the WENO-IM($k,A$) scheme, that is, WENO-IM($2, 0.1$), provides 
better numerical solutions \cite{WENO-IM} with less dissipation and 
higher resolution for the fifth-order WENO method than the WENO-JS, 
WENO-M and WENO-Z \cite{WENO-Z} schemes. However, it is required 
that the parameter $k$ has to be an even integer \cite{WENO-IM} in 
the WENO-IM($k, A$) method, which may prevent one from finding the 
best version from the family of the WENO-IM($k$, $A$) schemes. 
Besides, it is easy to verify that the mapping function of the 
WENO-IM($k,A$) method could not satisfy the requirements, namely 
$g'(0) = 0$ or $g'(0) = 1$, to decrease the effect from the 
non-smooth stencils, which may lead to the numerical solution with 
non-physical oscillations near discontinuities \cite{WENO-RM260}. It 
was demonstrated by numerical experiments \cite{WENO-RM260} that the 
seventh- and ninth- order WENO-IM($2,0.1$) schemes generate the 
oscillations obviously when solving the linear advection problem 
with discontinuities proposed by Jiang and Shu \cite{WENO-JS} for a 
long output time. Actually, we found that the fifth-order WENO-IM($2
, 0.1$) also generates the oscillations when solving the linear 
advection problem with discontinuities by taking a long output time 
and a \textit{bigger} grid number, and we show the numerical results 
in subsection \ref{subsec:keyEx} of this paper. Therefore, the goal 
of this study is to design a new mapped WENO scheme that can obtain 
high resolution without generating spurious oscillations when 
solving problems with discontinuities for long output times.

In order to achieve this goal, by introducing several adaptive 
control functions and using a smoothing approximation of the signum 
function, we propose a group of new mapping functions, which 
satisfies: first, it has a very similar form to the mapping function 
of the WENO-IM($k,A$) method, but the permitted range of the 
parameter $k$ is extended from positive even integers to any 
positive integers; second, it is smooth enough and keeps $g'(0) =
1$ and $g'(1) = 1$ so that it can decrease the effect from the 
non-smooth stencils; third, it introduces the adaptive property, so 
that it may prevent the corresponding mapped WENO schemes from 
generating non-physical oscillations near the discontinuities in 
long time simulations. We prove that the optimal order of accuracy 
at or near the critical points in smooth regions can be recovered by 
using the new mapping functions. Extensive numerical experiments 
show that, our proposed schemes perform satisfactorily for those 
benchmark problems in the references. Observations from the given 
numerical results of the one-dimensional linear advection problems 
with discontinuities at long output times (discussed carefully and 
presented detailly in subsection \ref{subsec:keyEx} of this paper) 
manifest that two of the proposed schemes are able to obtain high 
resolution without producing spurious oscillations, and this is a 
major improvement over other fifth-order WENO schemes: WENO-JS, 
WENO-M and WENO-IM(2, 0.1). Also, one of these two schemes provides 
improved behavior on calculating one-dimensional Euler system cases. 
Furthermore, we will see that the advantage of this scheme (which is 
our recommondation) seems more salient in two-dimensional Euler 
system cases.

The rest of this paper is organized as follows. In Section 2, we give
a brief description of the finite volume method and the procedures 
of the WENO-JS \cite{WENO-JS}, WENO-M \cite{WENO-M} and WENO-IM($k,A$
) \cite{WENO-IM}  schemes to clarify our major concern. In Section 
3, the details on how we construct the modified adaptive improved 
mapped WENO method, referred as WENO-MAIM$i$ later on, are 
presented. In Section 4, some numerical experiments are presented to 
compare the performances of different WENO methods. Finally, some 
concluding remarks are made in Section 5.


%% file: article_fvm.tex
\section{Description of finite volume WENO methods}
\label{secDescriptionFVM-WENO}
In this section, we first review the implementation of the finite 
volume method \cite{LeVeque-FVM} and then recall the essentials of 
the classic WENO-JS scheme proposed in \cite{WENO-JS}, along with 
the mapped version WENO-M introduced in \cite{WENO-M} and the 
WENO-IM($k, A$) scheme introduced in \cite{WENO-IM}.

\subsection{Finite volume method}
\label{subsecFVM}
Considering the one-dimensional scalar case, we rewrite Eq.(\ref
{governingEquation}) as follows,
\begin{equation}
\begin{array}{ll}
\dfrac{\partial u}{\partial t} + \dfrac{\partial f(u)}{\partial x} = 
0, & x \in [x_{l}, x_{r}], t > 0.
\end{array}
\label{1DScalarGoverningEquation}
\end{equation}
For simplicity, we assume that the computational domain is 
distributed into smaller uniform cells $I_{j} = [x_{j - 1/2}, x_{j + 
1/2}]$, where $\Delta x = x_{j + 1/2} - x_{j - 1/2}$ is the mesh 
width, $x_{j \pm 1/2} = x_{j} \pm \frac{\Delta x}{2}$ are the 
interfaces of $I_{j}$ and $x_{j} = \frac{1}{2}(x_{j + 1/2} + x_{j - 1
/2})$ are the cell centers. After some simple mathematical 
manipulations, we can approximate Eq.(\ref{1DScalarGoverningEquation}
) by the following finite volume conservative formulation
\begin{equation*}
\dfrac{\mathrm{d} \bar{u}_{j}(t)}{\mathrm{d} t} \approx - \dfrac{1}{
\Delta x}\big( \hat{f}_{j + 1/2} - \hat{f}_{j - 1/2} \big), 
\end{equation*}
where $\bar{u}_{j}(t)$ is the numerical approximation to the cell 
average $\bar{u}_{j} = \dfrac{1}{\Delta x}\int _{x_{j - 1/2}}^{x_{j +
1/2}}u(\xi, t) \mathrm{d} \xi$ and the numerical flux $\hat{f}_{j + 
1/2}$ is a function of $u(x, t)$ at the cell boundary, namely, $u_{j 
+ 1/2}^{\pm}$, defined by
\begin{equation}
\hat{f}_{j + 1/2} = \hat{f}(u_{j + 1/2}^{-}, u_{j + 1/2}^{+}),
\label{numericalFlux}
\end{equation}
where $\hat{f}(u^{-}, u^{+})$ is a monotone numerical flux. In this 
paper, we take the global Lax-Friedrichs flux $\hat{f}(a, b) = \dfrac
{1}{2}\big[ f(a) + f(b) - \alpha (b - a) \big]$, where $\alpha = \max
_{u} \lvert f'(u) \rvert$ is a constant and the maximum is taken 
over the whole range of $u$. In Eq.(\ref{numericalFlux}), $u_{j + 1/2
}^{\pm}$ can be obtained from some kind of WENO reconstruction, 
which is detailed in the following subsections. For the systems of 
conservation laws, a local characteristic decomposition is used in 
the reconstruction, and \cite{WENO-JS,WENOoverview} are referred to 
for more details.

\subsection{The classic WENO-JS reconstruction}
\label{subsecWENO-JS}
The three $3$rd-order approximations of $u(x_{j + 1/2}, t)$ in the 
left-biased substencils $S_{3}^{s} = \big\{ x_{j +s - 2}, x_{j + s - 
1}, x_{j + s} \big\}$ are defined as
\begin{equation}
u^{s,-}_{j + 1/2} = \sum\limits_{i = 0}^{2}c_{si}\bar{u}_{j + s - i}
, \quad s = 0,1,2,
\label{ENO3}
\end{equation}
where $c_{si}$ are Lagrangian interpolation coefficients (see 
\cite{WENOoverview,WENO-JS}) that depend on parameter $s$ but not on 
the values of $\bar{u}_{j}$. Explicitly, we have
\begin{equation*}
\begin{array}{l}
\begin{aligned}
u^{0,-}_{j + 1/2} &= \dfrac{1}{3}\bar{u}_{j - 2} - \dfrac{7}{6}\bar{u
}_{j - 1} + \dfrac{11}{6}\bar{u}_{j}, \\
u^{1,-}_{j + 1/2} &= -\dfrac{1}{6}\bar{u}_{j - 1} + \dfrac{5}{6}\bar{
u}_{j} + \dfrac{1}{3}\bar{u}_{j + 1},  \\
u^{2,-}_{j + 1/2} &= \dfrac{1}{3}\bar{u}_{j} + \dfrac{5}{6}\bar{u}_{j
 + 1} -\dfrac{1}{6}\bar{u}_{j + 2}.
\end{aligned}
\end{array}
\end{equation*}

The fifth-degree polynomial approximation $u_{j + 1/2}^{-} = u(x_{j +
 1/2}, t) + O(\Delta x^{5})$ is built via the convex combination of 
 the interpolated values $u^{s,-}_{j + 1/2}$ in Eq.(\ref{ENO3})
\begin{equation*}
u_{j + 1/2}^{-} = \sum\limits_{s = 0}^{2}\omega_{s}u_{j + 1/2}^{s,-},
\end{equation*}
where $\omega_{s}$ are called nonlinear weights taking the 
smoothness of the solution into consideration. In the classic 
WENO-JS reconstruction, they are computed as
\begin{equation} 
\omega_{s}^{\mathrm{JS}} = \dfrac{\alpha_{s}^{\mathrm{JS}}}{\sum_{l =
 0}^{2} \alpha_{l}^{\mathrm{JS}}}, \alpha_{s}^{\mathrm{JS}} = \dfrac{
 d_{s}}{(\epsilon + \beta_{s})^{2}}.
\label{nonlinearWeightsWENO-JS}
\end{equation}
Here, $\epsilon$ is a small positive number that is introduced to 
prevent the denominator becoming zero, and $d_{s}$ are the optimal 
weights satisfying $\sum\limits_{s = 0}^{2} d_{s} u^{s,-}_{j + 1/2} =
u(x_{j + 1/2}, t) + O(\Delta x^{5})$. For fifth-order WENO schemes, 
they are given by $d_{0} = 0.1, d_{1} = 0.6, d_{2} = 0.3$. The 
parameters $\beta_{s}$ are called smoothness indicators, which are 
defined as follows
\begin{equation*}
\begin{array}{l}
\begin{aligned}
\beta_{0} &= \dfrac{13}{12}\big(\bar{u}_{j - 2} - 2\bar{u}_{j - 1} + 
\bar{u}_{j} \big)^{2} + \dfrac{1}{4}\big( \bar{u}_{j - 2} - 4\bar{u}_
{j - 1} + 3\bar{u}_{j} \big)^{2}, \\
\beta_{1} &= \dfrac{13}{12}\big(\bar{u}_{j - 1} - 2\bar{u}_{j} + \bar
{u}_{j + 1} \big)^{2} + \dfrac{1}{4}\big( \bar{u}_{j - 1} - \bar{u}_{
j + 1} \big)^{2}, \\
\beta_{2} &= \dfrac{13}{12}\big(\bar{u}_{j} - 2\bar{u}_{j + 1} + \bar
{u}_{j + 2} \big)^{2} + \dfrac{1}{4}\big( 3\bar{u}_{j} - 4\bar{u}_{j 
+ 1} + \bar{u}_{j + 2} \big)^{2}.
\end{aligned}
\end{array}
\end{equation*}

 In smooth regions without critical points, the WENO-JS scheme gives 
 the fifth-order accuracy. However, near critical points, its order 
 of accuracy decreases to third order or even less. For more 
 details, we refer to \cite{WENO-M,WENO-JS}.

\subsection{The mapped WENO-M reconstruction}
\label{subsecWENO-M}
A mapping function of the nonlinear weights $\omega$ was constructed 
by Henrick et al. \cite{WENO-M} to correct the deficiency of the 
WENO-JS method mentioned above. The mapping function is written as
\begin{equation}
\big( g^{\mathrm{M}} \big)_{s}(\omega) = \dfrac{ \omega \big( d_{s} +
(d_{s})^2 - 3d_{s}\omega + \omega^{2} \big) }{ (d_{s})^{2} + (1 - 2d_
{s})\omega }, \quad \quad \omega \in [0, 1],
\label{mappingFunctionWENO-M}
\end{equation}
It is easy to verify that the mapping function $\big( g^{\mathrm{M}} 
\big)_{s}(\omega)$ is a monotonically increasing function in $[0, 1]$
with finite slopes that satisfies the following properties.
\begin{lemma} 
If the mapping function is defined by Eq.(\ref
{mappingFunctionWENO-M}), it satisfies: \\

C1. $0 \leq \big( g^{\mathrm{M}} \big)_{s}(\omega) \leq 1, \big( g^{
\mathrm{M}} \big)_{s}(0) = 0, \big( g^{\mathrm{M}} \big)_{s}(1) = 1$;

C2. $\big( g^{\mathrm{M}} \big)_{s}(\omega) \approx 0$ if $\omega 
\approx 0$; $\big( g^{\mathrm{M}} \big)_{s}(\omega) \approx 1$ if $
\omega \approx 1$;

C3. $\big( g^{\mathrm{M}} \big)_{s}(d_{s}) = d_{s}, \big( g^{\mathrm{
M}} \big)_{s}'(d_{s}) = \big(g^{\mathrm{M}}\big)_{s}''(d_{s}) = 0$. 
\label{lemmaWENO-Mproperties}
\end{lemma}

Mapping function Eq.(\ref{mappingFunctionWENO-M}) is employed to 
obtain the mapped weights as
\begin{equation*}
\omega_{s}^{\mathrm{M}} = \dfrac{\alpha _{s}^{\mathrm{M}}}{\sum_{l = 
0}^{2} \alpha _{l}^{\mathrm{M}}}, \alpha_{s}^{\mathrm{M}} = \big( g^{
\mathrm{M}} \big)_{s}(\omega^{\mathrm{JS}}_{s}),
\end{equation*}

In smooth regions, the WENO-M scheme gives the fifth-order accuracy 
even near the first-order critical points where the first derivative 
vanishes, and \cite{WENO-M} can be referred to for more details.

\subsection{The improved mapped WENO-IM($k, A$) reconstructions}
\label{subsecWENO-IM}
An improved mapped WENO-IM($k, A$) reconstruction was proposed by 
Feng et al. \cite{WENO-IM}. By rewriting the mapping function Eq.(
\ref{mappingFunctionWENO-M}), they obtained a new type of mapping 
function of the form 
\begin{equation}
\big( g^{\mathrm{IM}} \big)_{s}(\omega; k, A) = d_{s} + \dfrac{\big( 
\omega - d_{s} \big)^{k + 1}A}{\big( \omega - d_{s} \big)^{k}A + 
\omega(1 - \omega)}, \quad A > 0, k = 2n, n \in \mathbb{N}^{+},
\label{mappingFunctionWENO-IM}
\end{equation}
and the corresponding improved mapped weights are given by
\begin{equation*}
\omega_{s}^{\mathrm{IM}} = \dfrac{\alpha _{s}^{\mathrm{IM}}}{\sum_{l 
= 0}^{2} \alpha _{l}^{\mathrm{IM}}}, \alpha_{s}^{\mathrm{IM}} = \big(
g^{\mathrm{IM}} \big)_{s}(\omega^{\mathrm{JS}}_{s}; k, A).
\end{equation*}
It is trivial to show that the mapping function Eq.(\ref
{mappingFunctionWENO-M}) belongs to the family of the improved 
mapping functions Eq.(\ref{mappingFunctionWENO-IM}) by choosing $A = 
1, k = 2$. Moreover, the improved mapping functions Eq.(\ref
{mappingFunctionWENO-IM}) have the following properties.
\begin{lemma} 
If the mapping functions are defined by Eq.(\ref
{mappingFunctionWENO-IM}), they satisfy: \\

C1. $\big( g^{\mathrm{IM}} \big)'_{s}(\omega; k, A) \geq 0, \omega 
\in [0, 1]$;

C2. $\big( g^{\mathrm{IM}} \big)_{s}(0; k, A) = 0, \big( g^{\mathrm{
IM}} \big)_{s}(d_{s}; k, A) = d_{s}, \big( g^{\mathrm{IM}} \big)_{s}
(1; k, A) = 1$; 

C3. $\big( g^{\mathrm{IM}} \big)'_{s}(d_{s}; k, A) = \cdots = \big( g
^{\mathrm{IM}} \big)^{(k)}_{s}(d_{s}; k, A) = 0, \big( g^{\mathrm{IM}
} \big)^{(k + 1)}_{s}(d_{s}; k, A) \neq 0$.
\label{lemmaWENO-IMproperties}
\end{lemma}

The detailed proof of Lemma \ref{lemmaWENO-IMproperties} can be 
found in the statement on page 456 in \cite{WENO-IM}.

\subsection{Time discretization}
\label{subsecTime-discretization}
Commonly, WENO schemes are employed in a method of lines (MOL) 
approach, where one discretizes space while leaving time continuous.
Then, following this approach, we turn the PDE Eq.(\ref
{1DScalarGoverningEquation}) into a large number of coupled ODEs, 
resulting in the system of equations
\begin{equation}
\dfrac{\mathrm{d}\bar{u}_{j}(t)}{\mathrm{d} t} = \mathcal{L}(u_{j}), 
\label{ODEs}
\end{equation}
where $\mathcal{L}(u_{j})$ is the result of the application of the 
WENO scheme and is defined as
\begin{equation*}
\mathcal{L}(u_{j}) := - \dfrac{1}{\Delta x}\big( \hat{f}_{j + 1/2} - 
\hat{f}_{j - 1/2} \big). 
\end{equation*}
Throughout this paper, we solve the ODEs system Eq.(\ref{ODEs}) 
using the explicit, third-order, TVD, Runge-Kutta method \cite
{ENO-Shu1988,SSPRK1998,SSPRK2001} as follows
\begin{equation*}
\begin{array}{l}
\begin{aligned}
&u^{(1)} = u^{n} + \Delta t \mathcal{L}(u^{n}), \\
&u^{(2)} = \dfrac{3}{4} u^{n} + \dfrac{1}{4} u^{(1)} + \dfrac{1}{4} 
\Delta t \mathcal{L}(u^{(1)}), \\
&u^{n + 1} = \dfrac{1}{3} u^{n} + \dfrac{2}{3} \Delta t \mathcal{L}(u
^{(2)}),
\end{aligned}
\end{array}
\end{equation*}
where $u^{(1)}$ and $u^{(2)}$ are the intermediate stages, $u^{n}$ 
is the value of $u$ at time level $t^{n} = n\Delta t$, and $\Delta t$
is the time step satisfying some proper CFL condition.


%% file: article_modify.tex
\section{The modified adaptive improved mapped WENO reconstructions}
\label{secWENO-MAIM}
\subsection{The adaptive control functions}
\label{subsec_adaptiveControlFunctions}
In order to obtain the property $g'(0) = g'(1) =1$, we introduce two 
adaptive control functions, namely $f_{i,s}^{\mathrm{Ada}0}$ and $f_{
i,s}^{\mathrm{Ada}1}$, in our new mapping functions which will be 
proposed in a later subsection. For $f_{i,s}^{\mathrm{Ada}0}$ and
$f_{i,s}^{\mathrm{Ada}1}$, the following two requirements need to 
be satisfied: (1) if $\omega$ tends to $0$ or $1$, the product
$\omega^{f_{i,s}^{\mathrm{Ada}0}}(1 - \omega)^{f_{i,s}^{\mathrm{Ada}
1}}$ tends to $0$ rapidly; (2) if $\omega$ tends to $d_{s}$, the 
product $\omega^{f_{i,s}^{\mathrm{Ada}0}}(1 - \omega)^{f_{i,s}^{
\mathrm{Ada}1}}$ will be relatively far from $0$. By the 
constraints of these requirements, we design the following five 
types of $f_{i,s}^{\mathrm{Ada}0}$ and $f_{i,s}^{\mathrm{Ada}1}$, 
where the subscript $i = 1, 2, 3, 4, 5$ stands for Type1, Type2, 
Type3, Type4 and Type5, respectively.

Type 1:
\begin{equation}
f_{1,s}^{\mathrm{Ada}0} = \dfrac{d_{s}}{m_{s}\omega + \epsilon_{
\mathrm{A}}},f_{1,s}^{\mathrm{Ada}1} = \dfrac{1 - d_{s}}{m_{s}(1 - 
\omega) + \epsilon_{\mathrm{A}}},
\label{eq:fAda-type1}
\end{equation}
where $m_{s} \in \Big[ \frac{\alpha_{s}}{k + 1}, M \big)$ with $M$ 
being a finite positive constant real number. Here, $\alpha_{s}$ is 
a positive constant that depends on only the parameters $r$ and $s$. 
A detailed description of $\alpha_{s}$ is given in Property \ref
{property:helperFunctions02} and its proof below, and the 
recommended values of $\alpha_{s}$ for different order WENO schemes 
are provided in Table \ref{table:recommended_alpha_s} of Appendix A. 
In Eq.(\ref{eq:fAda-type1}), $\epsilon_{\mathrm{A}}$ is a very small 
positive number to prevent the denominator becoming zero, and the 
same holds for the following Eq.(\ref{eq:fAda-type3})(\ref
{eq:fAda-type4}). We will drop it in the theoretical analysis of 
this paper for simplicity.

Type 2:
\begin{equation}
f_{2,s}^{\mathrm{Ada}0} = f_{2,s}^{\mathrm{Ada}1} = \left\{
\begin{aligned}
\begin{array}{rl}
Q \cdot k, &\omega \leq \mathrm{CFS}_{s}, \\
1.0,   &\mathrm{CFS}_{s} < \omega < 1.0 - \dfrac{1.0-d_{s}}{d_{s}} 
\times \mathrm{CFS}_{s},\\
Q \cdot k, &\omega \geq 1.0 - \dfrac{1.0-d_{s}}{d_{s}} \times \mathrm
{CFS}_{s},
\end{array}
\end{aligned}\right.
\label{eq:fAda-type2}
\end{equation}
where $Q \geq k^{-1}$, $k$ is the same as in Eq.(\ref
{mappingFunctionMAIM}) and the \textbf{C}ontrol \textbf{F}actor of 
\textbf{S}moothness $\mathrm{CFS}_{s} \in (0, d_{s}]$.

Type 3:
\begin{equation}
f_{3,s}^{\mathrm{Ada}0} = f_{3,s}^{\mathrm{Ada}1} = \dfrac{\max
\limits_{0 \leq j \leq 2}\Big(\beta_{j}\Big)}{\min\limits_{0 \leq j 
\leq 2}\Big(\beta_{j}\Big) + \epsilon_{\mathrm{A}}}.
\label{eq:fAda-type3}
\end{equation}

Type 4:
\begin{equation}
f_{4,s}^{\mathrm{Ada}0} = f_{4,s}^{\mathrm{Ada}1} = \dfrac{\max
\limits_{0 \leq j \leq 2}\bigg( \dfrac{\omega_{j}}{d_{j}} \bigg)}{
\min\limits_{0 \leq j \leq 2}\bigg( \dfrac{\omega_{j}}{d_{j}} \bigg) 
+ \epsilon_{\mathrm{A}}}.
\label{eq:fAda-type4}
\end{equation}

Type 5:
\begin{equation}
f_{5,s}^{\mathrm{Ada}0} = f_{5,s}^{\mathrm{Ada}1} = C,
\label{eq:fAda-type5}
\end{equation}
where $C$ is a constant and $C \geq 1$.

\begin{remark}
For simplicity, we mainly consider the case of Type 1 as an example 
in the following theoretical analysis. However, we present the 
mapping function curves and show the convergence order of the 
corresponding WENO schemes, as well as the numerical results, for 
the cases of Type 1 to Type 4. Here, we define Type 5 for 
completeness: it will be used only in the discussion of the mapping 
function curves as it is not truly adaptive.
\label{remark:fAda}
\end{remark}

\subsection{A smoothing approximation to the signum function}
\label{subsec_sgFunction}
To extend the range of the parameter $k$ in the mapping functions
of WENO-IM($k,A$), we design a smoothing approximation to the 
well-known non-smoothing signum function as
\begin{equation}
sg\big( x, \delta \big) = \left\{ 
\begin{array}{ll}
\begin{aligned}
&\dfrac{x}{|x|}, & |x| \geq \delta, \\ 
&\dfrac{x}{\big( \delta ^2 - x^2 \big)^{k + 3} + \ |x|}, & |x| < 
\delta,
\end{aligned}
\end{array} \right.
\label{sgFunction}
\end{equation} 
where the constant $\delta > 0$ and $\delta \rightarrow 0$, $k \in 
\mathbb{N}^{+}$. It is trivial to verify the following properties. 
\begin{property}
The smoothing approximation $sg(x, \delta)$ is $(k + 2)$th-order 
differentiable.
\label{propertySg1}
\end{property}

\begin{property}
The smoothing approximation $sg(x, \delta)$ has the following 
properties: \\

P1. $sg(x, \delta)$ is a bounded smoothing function; 

P2. if $x \neq 0$, then $sg(x, \delta) \neq 0$; 

P3. $sg(x, \delta)$ and $x$ always have the same sign or we always 
have $sg(x, \delta) \cdot x^{2n - 1} \geq 0, \forall n \in \mathbb{N}
^{+}$;

P4. $\frac{\mathrm{d}}{\mathrm{d}x}sg(x, \delta) \geq 0$; 

P5. $sg(0, \delta) = 0$ and $\frac{\mathrm{d}}{\mathrm{d}x}sg(0, 
\delta) \neq 0$; 

P6. if $\lvert x \rvert > \delta$, then $\frac{\mathrm{d}}{\mathrm{d}
x}sg(x, \delta) = 0$.
\label{propertySg2}
\end{property}

As the proofs of the above two properties are simple, we do not 
provide them here.

\subsection{A class of modified adaptive improved mapping functions}
\label{subsec_MAIM-mappingFunction}
Now, employing the adaptive control functions $f_{i,s}^{\mathrm{Ada}0
},f_{i,s}^{\mathrm{Ada}1}$ and the smoothing approximation $sg(x,
\delta)$, we propose a new family of modified adaptive mapping 
functions as
\begin{equation}
\big( g^{\mathrm{MAIM}i} \big)_{s}\big( \omega \big) = d_{s} + \dfrac
{f_{s}^{\mathrm{MAIM}} \cdot (\omega - d_{s})^{k + 1}}{f_{s}^{\mathrm
{MAIM}} \cdot (\omega - d_{s})^{k} + \omega^{f_{i,s}^{\mathrm{Ada}0}}
 (1 - \omega)^{f_{i,s}^{\mathrm{Ada}1}}}, \qquad k \in N^{+}, A > 0,
\label{mappingFunctionMAIM}
\end{equation}
where $f_{s}^{\mathrm{MAIM}}$ is a function of $\omega, A, k, d_{s}, 
\delta$ and it is defined as
\begin{equation}
f_{s}^{\mathrm{MAIM}} = A \bigg( \dfrac{1 + (-1)^{k}}{2} + \frac{1 + 
(-1)^{k + 1}}{2} \cdot sg\big( \omega - d_{s}, \delta \big) \bigg).
\label{fimFunction}
\end{equation}

Notably, in our modified adaptive mapping functions, the parameter 
$k$ is no longer limited to even integers but to all positive 
integers. Additionally, if one sets $k = 2n, n \in \mathbb{N}^{+}, A 
> 0$ and $f_{i,s}^{\mathrm{Ada}0} = f_{i,s}^{\mathrm{Ada}1} = 1$ 
(choosing $i = 5$ and $C = 1$), the family of mapping functions for 
WENO-IM($k, A$) in \cite{WENO-IM} is obtained immediately. 
It means that the family of mapping functions of Feng et al. 
\cite{WENO-IM} belongs to our new family of modified adaptive 
mapping functions (\ref{mappingFunctionMAIM}). Naturally, the 
mapping function proposed by Herick et al. in \cite{WENO-M} also 
belongs to the new family of modified adaptive mapping functions by 
setting $k = 2, A = 1$ and $f_{i,s}^{\mathrm{Ada}0} = f_{i,s}^{
\mathrm{Ada}1} = 1$.

Before giving Theorem \ref{theorem_maimMappingFunction} and its 
proof, we provide the following lemma and properties.
\begin{lemma}
The function $f_{s}^{\mathrm{MAIM}}$ satisfies:\\

C1. it is a bounded smoothing function; 

C2. if $\omega \neq d_{s}$, then $f_{s}^{\mathrm{MAIM}} \neq 0$;

C3. $f_{s}^{\mathrm{MAIM}} \cdot (\omega - d_{s})^{k} \geq 0$; 

C4. $\frac{\mathrm{d}f_{s}^{\mathrm{MAIM}}}{\mathrm{d}\omega} \cdot (
\omega - d_{s})^{k + 1} \geq 0$; 

C5. if $k = 2n, n \in \mathbb{N}^{+}$, then $f_{s}^{\mathrm{MAIM}}
\Big\lvert_{\omega = d_{s}} \neq 0$, and if $k = 2n - 1, n \in 
\mathbb{N}^{+}$, then $f_{s}^{\mathrm{MAIM}}\Big\lvert_{\omega = d_{s
}} = 0$ and $\frac{\mathrm{d}f_{s}^{\mathrm{MAIM}}}{\mathrm{d}\omega}
\Big\lvert_{\omega = d_{s}} \neq 0$; 

C6. $f_{s}^{\mathrm{MAIM}}$ is at least $(k + 2)$th-order 
differentiable with respect to $\omega$.
\label{Lemma_fimFunction}
\end{lemma}
\textbf{Proof.}

(1) First, if $k = 2n, n \in \mathbb{N}^{+}$, we directly have $f_{s}
^{\mathrm{MAIM}} = A$. As $A$ is a constant, $A > 0$, and $k = 2n, n 
\in \mathbb{N}^{+}$; thus, i) $f_{s}^{\mathrm{MAIM}}$ is bounded and 
smooth, then $C\mathit{1}$ is true; ii) for $\forall \omega \in [0, 1
]$, we have $f_{s}^{\mathrm{MAIM}} = A > 0$, then $C\mathit{2}$ is 
true; iii) $C\mathit{3}$ is true as $f_{s}^{\mathrm{MAIM}} \cdot (
\omega - d_{s})^{k} = A \cdot (\omega - d_{s})^{2n} \geq 0$; iv) $C
\mathit{4}$ is true as $\frac{\mathrm{d}f_{s}^{\mathrm{MAIM}}}{
\mathrm{d}\omega} \cdot (\omega - d_{s})^{k + 1} = \frac{\mathrm{d}A}
{\mathrm{d}\omega} \cdot (\omega - d_{s})^{k + 1} = 0$; v) for $
\forall \omega \in [0, 1]$, we have $f_{s}^{\mathrm{MAIM}} = A \neq 
0$, then $C\mathit{5}$ is true; vi) $C\mathit{6}$ is true as $\frac{
\mathrm{d}^{q}f_{s}^{\mathrm{MAIM}}}{\mathrm{d}\omega ^{q}} = \frac{
\mathrm{d}^{q} A}{\mathrm{d}\omega ^{q}} = 0, \forall q \in \mathbb{N
}^{+}$.

(2) Then, we give the proof in the case that $k = 2n - 1, n \in 
\mathbb{N}^{+}$. At present, we can rewrite Eq.(\ref{fimFunction}) as
\begin{equation*}
f_{s}^{\mathrm{MAIM}} = A \cdot sg(\omega - d_{s}, \delta).
\end{equation*}
As $A$ is a constant, $A > 0$, and $k = 2n - 1, n \in \mathbb{N}^{+}$
, letting $x = \omega - d_{s}$ and substituting it into Eq.(\ref
{sgFunction}), we obtain the following results trivially: i) 
according to $P\mathit{1}$ of Property \ref{propertySg2}, we know 
that $f_{s}^{\mathrm{MAIM}}$ is bounded and smooth, then
$C\mathit{1}$ is true; ii) if $\omega \neq d_{s} \Rightarrow x \neq 0
\Rightarrow$ $P\mathit{2}$ of Property \ref{propertySg2} is 
satisfied, then $f_{s}^{\mathrm{MAIM}} \neq 0$ and $C\mathit{2}$ is 
true; iii) if $k = 2n - 1, n \in \mathbb{N}^{+} \Rightarrow$ $P
\mathit{3}$ of Property \ref{propertySg2} is satisfied, then $f_{s}^{
\mathrm{MAIM}} \cdot (\omega - d_{s})^{k} \geq 0$ and $C\mathit{3}$ 
is true; iv) according to $P\mathit{4}$ of Property \ref{propertySg2}
, we have $\frac{\mathrm{d}f^{\mathrm{MAIM}}}{\mathrm{d}\omega} = sg(
x, \delta) \frac{\mathrm{d}A}{\mathrm{d}x} + A \frac{\mathrm{d}sg(x, 
\delta)}{\mathrm{d}x} \geq 0$, and it is easy to find that $(\omega -
 d_{s})^{k + 1} \geq 0, (k = 2n - 1, n \in \mathbb{N}^{+})$, thus $
 \frac{\mathrm{d}f^{\mathrm{MAIM}}}{\mathrm{d}\omega} \cdot (\omega -
  d_{s})^{k + 1} \geq 0$ and $C\mathit{4}$ is true; v) as $\omega = d
 _{s}$, namely, $x = 0$, so $P\mathit{5}$ of Property \ref
 {propertySg2} is satisfied, then we have $f_{s}^{\mathrm{MAIM}}\Big
 \lvert_{\omega = d_{s}} = A \cdot sg(\omega - d_{s}, \delta)\Big
 \lvert_{\omega = d_{s}} = 0$ and $\frac{\mathrm{d}f_{s}^{\mathrm{
 MAIM}}}{\mathrm{d}\omega} \Big\lvert_{\omega = d_{s}} = \Big(A\frac{
 \mathrm{d}}{\mathrm{d}x}sg(x, \delta) + sg(x, \delta)\frac{\mathrm{d
 }A}{\mathrm{d}x} \Big)\Big\lvert_{x = 0} \neq 0$, thus,
 $C\mathit{5}$ is true; vi) as $A$ is a constant, we obtain $\frac{
 \mathrm{d}^{m}A}{\mathrm{d}x^{m}} = 0, \forall m \in \mathbb{N}^{+}$
 , and according to Property \ref{propertySg1}, we have $sg(x, \delta
 ) \in C^{k + 2}\big( \mathbb{R} \big)$, thus, $C\mathit{6}$ is true 
 as $\frac{\mathrm{d}^{q} f_{s}^{\mathrm{MAIM}}}{\mathrm{d}\omega ^{q
 }} = A\frac{\mathrm{d}^{q}}{\mathrm{d}x^{q}}sg(x, \delta), q = 1, 2
 , \cdots, k + 2.$
$\hfill\square$ \\
\begin{property}
Let
\begin{equation*}
\displaystyle q_{s}(\omega) = (\omega - d_{s})\bigg( \frac{d_{s}}{
\omega^{2}}\big( 1 - \ln \omega \big) + \frac{1 - d_{s}}{(1 - \omega)
^{2}}\big( \ln(1 - \omega) - 1 \big) \bigg), \quad \omega \in (0, 1)
\end{equation*}
then, the maximum value of $q_{s}(\omega)$ is finite for all optimal 
weights $d_{s}$ in various ($2r - 1$)th-order WENO schemes with $r 
\in \{2, \cdots, 9 \}$ (the values of all these $d_{s}$ can be found 
in Table 3 of \cite{veryHighOrderWENO}).
\label{property:helperFunctions01}
\end{property}
\textbf{Proof.}

After simple mathematical manipulations, we obtain
\begin{equation*}
q_{s}'(\omega) = \dfrac{d_{s}}{\omega^{2}}\Big( 1 - \ln \omega \Big) 
+ \dfrac{1 - d_{s}}{(1 - \omega)^{2}}\Big( \ln(1 - \omega) - 1 \Big) 
+ (\omega - d_{s})\Bigg( - \dfrac{2d_{s}}{\omega^{3}}\Big( 1 - \ln 
\omega \Big) - \dfrac{d_{s}}{\omega^{3}} + \dfrac{2(1 - d_{s})}{(1 - 
\omega)^{3}}\Big( \ln (1 - \omega) - 1 \Big) - \dfrac{1 - d_{s}}{(1 -
 \omega)^{3}} \Bigg).
\end{equation*}
Clearly, $q_{s}'(\omega)$ is continuous in $(0,1)$. Therefore, if we 
can verify that there is one and only one value of 
$\omega \in (0, 1)$, namely, $\omega_{\mathrm{crit}}$, satisfying $q_
{s}'(\omega_{\mathrm{crit}}) = 0$, $q_{s}'(\omega) > 0$, $\forall 
\omega \in (0, \omega_{\mathrm{crit}})$ and $q_{s}'(\omega) < 0$, $
\forall \omega \in (\omega_{\mathrm{crit}}, 1)$, we can prove 
Property \ref{property:helperFunctions01}. Unfortunately, this 
direct theoretical verification is challenging. However, for a fixed 
value of $d_{s}$, we can easily obtain the solution of $q_{s}'(\omega
) = 0$ and get the curve of $q_{s}'(\omega)$ by numerical means. 
After extensive calculations using software MATLAB, we found that, 
for every optimal weight $d_{s}$ in various $(2r - 1)$th-order WENO 
schemes with $r \in \{2,\cdots,9\}$, the requirements of $q_{s}'(
\omega) = 0$, $q_{s}'(\omega) > 0$ for $\omega \in (0, \omega_{
\mathrm{crit}})$ and $q_{s}'(\omega) < 0$ for $\omega \in (\omega_{
\mathrm{crit}}, 1)$ are all satisfied. We show the maximum value of
$q_{s}(\omega)$ for various $(2r - 1)$th-order WENO schemes with $r 
\in \{2,\cdots,9\}$ in Table \ref{table:recommended_alpha_s} in 
Appendix A. 
$\hfill\square$

\begin{property}
Let $T_{s}(\omega) = \omega^{\frac{d_{s}}{m_{s}\omega}} \cdot \Big( 1
- \omega \Big)^{\frac{1 - d_{s}}{m_{s}(1 - \omega)}}, Q_{s}(\omega) =
k + 1 - \dfrac{q_{s}(\omega)}{m_{s}}$ and $P_{s}(\omega) = T_{s}(
\omega)Q_{s}(\omega)$, then, $T_{s}(\omega), Q_{s}(\omega), P_{s}(
\omega)$ possess the following properties: \\

P1. $T_{s}(\omega) \geq 0, \forall \omega \in [0, 1]$ and $\lim
\limits_{\omega \to 0^{+}} T_{s}(\omega)= \lim\limits_{\omega \to 1^{
-}} T_{s}(\omega)= 0$; 

P2. if $m_{s} \in \Big[ \frac{\alpha_{s}}{k + 1}, M \Big)$, where $M$
is a finite positive constant real number and $\alpha_{s} > \max\Big
( 0, \max\limits_{0 < \omega < 1} q_{s}(\omega) \Big)$, then $Q_{s}(
\omega) > 0, \forall \omega \in (0, 1)$;

P3. $\lim\limits_{\omega \rightarrow 0^{+}} P_{s}(\omega) = \lim
\limits_{\omega \rightarrow 1^{-}} P_{s}(\omega) = 0$.
\label{property:helperFunctions02}
\end{property}
\textbf{Proof.}

(1) It is easy to verify $P\mathit{1}$ of Property \ref
{property:helperFunctions02}, as $\omega \in [0, 1]$ and $\frac{d_{s}
}{m_{s}\omega} > 0, \frac{1 - d_{s}}{m_{s}(1 - \omega)} > 0$. 

(2) Clearly, when $q_{s}(\omega) \leq 0$, we have $Q_{s}(\omega) > 0$
as $m_{s} > 0$. If $q_{s}(\omega) > 0$, as $m_{s} \in \Big[ \frac{
\alpha_{s}}{k + 1}, M \Big)$, according to Property \ref
{property:helperFunctions01}, we have 
\begin{equation*}
Q_{s}(\omega) = k + 1 - \dfrac{q_{s}(\omega)}{m_{s}} \geq k + 1 - 
\dfrac{k + 1}{\alpha_{s}}q_{s}(\omega) > (k + 1)\Bigg( 1 - \dfrac{q_{
s}(\omega)}{\max\limits_{0 < \omega < 1}q_{s}(\omega)} \Bigg) \geq (k
+ 1)\Bigg( 1 - \dfrac{\max\limits_{0 < \omega < 1}q_{s}(\omega)}{\max
\limits_{0 < \omega < 1}q_{s}(\omega)} \Bigg) = 0.
\end{equation*}
Thus, $P\mathit{2}$ of Property \ref{property:helperFunctions02} is 
true.

(3) If $n > 0$, by employing L'Hospital's Rule, it is trivial to 
show that
\begin{equation*}
\lim\limits_{x \to 0^{+}} \Bigg( x^{n}\ln \dfrac{1}{x} \Bigg) = \lim
\limits_{x \to 0^{+}} \dfrac{\ln \dfrac{1}{x}}{\bigg( \dfrac{1}{x} 
\bigg)^{n}} = \lim\limits_{t \to + \infty} \dfrac{\ln t}{t^{n}} = 
\lim\limits_{t \to + \infty} \dfrac{1}{nt^{n}} = 0.
\end{equation*} 
Then, as $\lim\limits_{\omega \to 0^{+}} \Big( \frac{d_{s}}{m_{s}
\omega} - 2 \Big) \gg 0$ and $\lim\limits_{\omega \to 1^{-}} \Big( 
\frac{1 - d_{s}}{m_{s}(1 - \omega)} - 2 \Big) = \lim\limits_{\tilde{
\omega} \to 0^{+}} \Big( \frac{1 - d_{s}}{m_{s}\tilde{\omega}} - 2 
\Big) \gg 0$, we obtain
\begin{equation*}
\lim\limits_{\omega \to 0^{+}} \Bigg( \omega^{\frac{d_{s}}{m_{s}
\omega}}\dfrac{1}{\omega^{2}}\Big(-\ln \omega \Big) \Bigg) = \lim
\limits_{\omega \to 0^{+}} \Bigg( \omega^{\frac{d_{s}}{m_{s}\omega} -
2}\ln \dfrac{1}{\omega} \Bigg) = 0,
\end{equation*}
\begin{equation*}
\lim\limits_{\omega \to 1^{-}} \Bigg( \Big(1 - \omega\Big)^{\frac{1 -
d_{s}}{m_{s}(1 - \omega)}}\dfrac{1}{( 1 -\omega )^{2}}\ln (1 - \omega
) \Bigg) = \lim\limits_{\tilde{\omega} \to 0^{+}} \Bigg( - \tilde{
\omega}^{\frac{1 - d_{s}}{m_{s}\tilde{\omega}} - 2}\ln \dfrac{1}{
\tilde{\omega}} \Bigg) = 0.
\end{equation*}
Now, we can prove $P\mathit{3}$ of Property \ref
{property:helperFunctions02} as
\begin{equation*}
\begin{array}{lll}
\lim\limits_{\omega \to 0^{+}} P_{s}(\omega) &=& \lim\limits_{\omega 
\to 0^{+}} \Bigg\{ \omega^{\frac{d_{s}}{m_{s}\omega}} \cdot \Big( 1 -
\omega \Big)^{\frac{1 - d_{s}}{m_{s}(1 - \omega)}} \Bigg( k + 1 - 
\dfrac{1}{m_{s}} (\omega - d_{s})\bigg( \frac{d_{s}}{\omega^{2}}\big
( 1 - \ln \omega \big) + \frac{1 - d_{s}}{(1 - \omega)^{2}}\big( \ln
(1 - \omega) - 1 \big) \bigg) \Bigg) \Bigg\} \\ 
\quad &=& \lim\limits_{\omega \to 0^{+}} \Bigg\{ (k + 1)\cdot\Big( 1 
- \omega \Big)^{\frac{1 - d_{s}}{m_{s}(1 - \omega)}} \times \Bigg( 
\omega^{\frac{d_{s}}{m_{s}\omega}} \Bigg) \Bigg\} \\
\quad & \quad +& \lim\limits_{\omega \to 0^{+}} \Bigg\{ -\dfrac{d_{s}
}{m_{s}}(\omega - d_{s}) \cdot \Big( 1 - \omega \Big)^{\frac{1 - d_{s
}}{m_{s}(1 - \omega)}} \times \Bigg( \omega^{\frac{d_{s}}{m_{s}\omega
} - 2} \Bigg)\Bigg\} \\
\quad & \quad +& \lim\limits_{\omega \to 0^{+}} \Bigg\{ -\dfrac{d_{s}
}{m_{s}}(\omega - d_{s}) \cdot \Big( 1 - \omega \Big)^{\frac{1 - d_{s
}}{m_{s}(1 - \omega)}} \times \Bigg( \omega^{\frac{d_{s}}{m_{s}\omega
}}\dfrac{1}{\omega^{2}}\Big(-\ln \omega \Big) \Bigg)\Bigg\} \\
\quad & \quad +& \lim\limits_{\omega \to 0^{+}} \Bigg\{ -\dfrac{1 - d
_{s}}{m_{s}}(\omega - d_{s}) \cdot \Big( 1 - \omega \Big)^{\frac{1 - 
d_{s}}{m_{s}(1 - \omega)}} \cdot \dfrac{\ln (1 - \omega) - 1}{(1 - 
\omega)^{2}} \times \Bigg( \omega^{\frac{d_{s}}{m_{s}\omega}} \Bigg)
\Bigg\} \\
\quad &=0,&
\end{array}
\end{equation*}
and
\begin{equation*}
\begin{array}{lll}
\lim\limits_{\omega \to 1^{-}} P_{s}(\omega) &=& \lim\limits_{\omega 
\to 1^{-}} \Bigg\{ \omega^{\frac{d_{s}}{m_{s}\omega}} \cdot \Big( 1 -
\omega \Big)^{\frac{1 - d_{s}}{m_{s}(1 - \omega)}} \Bigg( k + 1 - 
\dfrac{1}{m_{s}} (\omega - d_{s})\bigg( \frac{d_{s}}{\omega^{2}}\big
( 1 - \ln \omega \big) + \frac{1 - d_{s}}{(1 - \omega)^{2}}\big( \ln(
1 - \omega) - 1 \big) \bigg) \Bigg) \Bigg\} \\ 
\quad &=& \lim\limits_{\omega \to 1^{-}} \Bigg\{ (k + 1)\cdot\omega^{
\frac{d_{s}}{m_{s}\omega}} \times \Bigg( \Big( 1 - \omega \Big)^{
\frac{1 - d_{s}}{m_{s}(1 - \omega)}} \Bigg) \Bigg\} \\
\quad & \quad +& \lim\limits_{\omega \to 1^{-}} \Bigg\{ -\dfrac{d_{s}
}{m_{s}}(\omega - d_{s}) \cdot \omega^{\frac{d_{s}}{m_{s}\omega}} 
\cdot \dfrac{1 - \ln\omega}{\omega^{2}} \times \Bigg( \Big( 1 - 
\omega \Big)^{\frac{1 - d_{s}}{m_{s}(1 - \omega)}} \Bigg) \Bigg\} \\ 
\quad & \quad +& \lim\limits_{\omega \to 1^{-}} \Bigg\{ -\dfrac{1 - d
_{s}}{m_{s}}(\omega - d_{s}) \cdot \omega^{\frac{d_{s}}{m_{s}\omega}}
\times \Bigg( \Big(1 - \omega\Big)^{\frac{1 - d_{s}}{m_{s}(1 - \omega
)}}\dfrac{1}{( 1 -\omega )^{2}}\ln (1 - \omega) \Bigg) \Bigg\} \\
\quad & \quad +& \lim\limits_{\omega \to 1^{-}} \Bigg\{ -\dfrac{1 - d
_{s}}{m_{s}}(\omega - d_{s}) \cdot \omega^{\frac{d_{s}}{m_{s}\omega}}
\times \Bigg( - \Big(1 - \omega\Big)^{\frac{1 - d_{s}}{m_{s}(1 - 
\omega)} - 2} \Bigg) \Bigg\} \\
\quad &=0.&
\end{array}
\end{equation*}
$\hfill\square$

\begin{theorem}
If the mapping functions are defined by Eq.(\ref{eq:fAda-type1}) and 
Eq.(\ref{mappingFunctionMAIM}), they satisfy:\\

C1. $\big( g^{\mathrm{MAIM}1} \big)_{s}(0) = 0, \big( g^{\mathrm{MAIM
}1} \big)_{s}(1) = 1, \big( g^{\mathrm{MAIM}1} \big)_{s}(d_{s}) = d_{
s}$;

C2. $\big( g^{\mathrm{MAIM}1} \big)'_{s}(\omega) \geq 0$, if $\omega 
\in (0, 1)$;

C3. If $k = 2n, n \in \mathbb{N}^{+}$, then $\big( g^{\mathrm{MAIM}1}
 \big)'_{s}(d_{s}) = \big( g^{\mathrm{MAIM}1} \big)''_{s}(d_{s}) = 
 \cdots = \big( g^{\mathrm{MAIM}1} \big)^{(k)}_{s}(d_{s}) = 0, \big( 
 g^{\mathrm{MAIM}1} \big)^{(k + 1)}_{s}(d_{s}) \neq 0$, and if $k = 2
 n - 1, n \in \mathbb{N}^{+}$, then $\big( g^{\mathrm{MAIM}1} \big)'_
 {s}(d_{s}) = \big( g^{\mathrm{MAIM}1} \big)''_{s}(d_{s}) = \cdots = 
 \big( g^{\mathrm{MAIM}1} \big)^{(k + 1)}_{s}(d_{s}) = 0, \big( g^{
 \mathrm{MAIM}1} \big)^{(k + 2)}_{s}(d_{s}) \neq 0$;

C4. $\big( g^{\mathrm{MAIM}1} \big)'_{s}(0^{+}) = \big( g^{\mathrm{
MAIM}1} \big)'_{s}(1^{-}) = 1$ for $m_{s} \in \Big[ \frac{\alpha_{s}}
{k + 1}, M \Big)$, where $M$ is a finite positive constant real 
number and $\alpha_{s} > \max\Big( 0, \max\limits_{0 < \omega < 1} q_
{s}(\omega) \Big)$.
\label{theorem_maimMappingFunction}
\end{theorem}
\textbf{Proof.}

As the bounded smoothing functions $f_{s}^{\mathrm{MAIM}}$ satisfy 
conditions $C\mathit{2}$ and $C\mathit{3}$ of Lemma \ref
{Lemma_fimFunction} and $d_{s} \in (0, 1)$, we have
\begin{equation*}
f_{s}^{\mathrm{MAIM}} \cdot (\omega - d_{s})^{k} + \omega^{f_{1,s}^{
\mathrm{Ada}0}} (1 - \omega)^{f_{1,s}^{\mathrm{Ada}1}} > 0, \quad 
\forall \omega \in [0, 1].
\end{equation*}
Therefore, the denominator of Eq.(\ref{mappingFunctionMAIM}) will 
never be zero; in other words, Eq.(\ref{mappingFunctionMAIM}) will 
always make sense. Thus, one can design the mapping functions 
according to Eq.(\ref{eq:fAda-type1}) and Eq.(\ref
{mappingFunctionMAIM}). Next, we prove $C1 \sim C4$ of Theorem \ref
{theorem_maimMappingFunction}.

(1) According to the $C\mathit{2}$ of Lemma \ref{Lemma_fimFunction} 
and $d_{s} \in (0, 1)$, we can obtain $f_{s}^{\mathrm{MAIM}}\Big
\lvert_{\omega = 0} \neq 0 $ and $f_{s}^{\mathrm{MAIM}}\Big\lvert_{
\omega = 1} \neq 0$; then,
\begin{equation*}
\begin{array}{l}
\begin{aligned}
&\big(g^{\mathrm{MAIM}1}\big)_{s}(\omega)\Big\lvert_{\omega = 0} = d_
{s} + \dfrac{f_{s}^{\mathrm{MAIM}}\cdot (\omega - d_{s})^{k + 1}}{f_{
s}^{\mathrm{MAIM}}\cdot (\omega - d_{s})^{k}} \Bigg\lvert_{\omega = 0
} = d_{s} - d_{s} = 0, \\
&\big(g^{\mathrm{MAIM}1}\big)_{s}(\omega)\Big\lvert_{\omega = 1} = d_
{s} + \dfrac{f_{s}^{\mathrm{MAIM}}\cdot (\omega - d_{s})^{k + 1}}{f_{
s}^{\mathrm{MAIM}}\cdot (\omega - d_{s})^{k}} \Bigg\lvert_{\omega = 1
} = d_{s} + (1 - d_{s}) = 1.
\end{aligned}
\end{array}
\end{equation*}
From the $C\mathit{1}$ of Lemma \ref{Lemma_fimFunction}, we know $f_{
s}^{\mathrm{MAIM}} \cdot (\omega - d_{s})^{k + 1} \Big\lvert_{\omega 
= d_{s}} = 0$, so
\begin{equation*}
\big(g^{\mathrm{MAIM}1}\big)_{s}(\omega)\Big\lvert_{\omega = d_{s}} =
 d_{s} + \dfrac{0}{d_{s}(1 - d_{s})}= d_{s}.
\end{equation*}

(2) As the parameter $\epsilon_{\mathrm{A}}$ in Eq.(\ref
{eq:fAda-type1}) is a very small number used only to prevent the 
denominator becoming zero, we drop it in the theoretical analysis. 
Then, taking the derivative of $\big( g^{\mathrm{MAIM}1} \big)_{s}(
\omega)$ with respect to $\omega$, we obtain
\begin{equation}
\big( g^{\mathrm{MAIM}1} \big)'_{s}(\omega) = \dfrac{\mathrm{P}_{num}
(\omega)}{\mathrm{P}_{den}(\omega)} = \dfrac{(f_{s}^{\mathrm{MAIM}})^
{2}\cdot(\omega - d_{s})^{2k} + \dfrac{\mathrm{d} f_{s}^{\mathrm{MAIM
}}}{\mathrm{d}\omega}\cdot(\omega - d_{s})^{k + 1}T_{s}(\omega) + f_{
s}^{\mathrm{MAIM}}\cdot(\omega - d_{s})^{k}P_{s}(\omega)}{\bigg(f_{s}
^{\mathrm{MAIM}}\cdot(\omega - d_{s})^{k} + T_{s}(\omega) \bigg)^{2}}
,
\label{eq:derivativeOfMappingfunctionWENO-MAIM}
\end{equation}
where $P_{s}(\omega), T_{s}(\omega), Q_{s}(\omega)$ are the same as 
in Property \ref{property:helperFunctions02}.

For $m_{s} \in \Big[ \frac{\alpha_{s}}{k + 1}, M \Big)$, according 
to $P\mathit{1}$ and $P\mathit{2}$ of Property \ref
{property:helperFunctions02}, we have $T_{s}(\omega) \geq 0, Q_{s}(
\omega) > 0, \forall \omega \in (0, 1)$. By employing the conditions
$C\mathit{3}$ and $C\mathit{4}$ of Lemma \ref{Lemma_fimFunction} and 
considering $(f_{s}^{\mathrm{MAIM}})^{2} \geq 0, (\omega - d_{s})^{2k
} \geq 0$, we conclude that
\begin{equation*}
\big( g^{\mathrm{MAIM}1} \big)'_{s}(\omega) \geq 0, \quad \forall 
\omega \in (0, 1).
\end{equation*}

(3) As $f_{s}^{\mathrm{MAIM}} \in C^{k + 2}\big( \mathbb{R} \big)$, 
which has been provided by $C\mathit{6}$ of Lemma \ref
{Lemma_fimFunction}, $(g^{\mathrm{MAIM}1})^{(q)}_{s}(\omega)$ always 
makes sense for $q = 1, 2, \cdots, k + 2$.

i) If $k = 2n, n \in \mathbb{N}^{+}$, we know that $f_{s}^{\mathrm{
MAIM}} \neq 0$; then, $(\omega - d_{s})^{k}\mid\mathrm{P}_{num}$, 
and $(\omega - d_{s})^{k + 1}\nmid\mathrm{P}_{num}$. As $T_{s}(\omega
)Q_{s}(\omega) > 0, \forall \omega \in (0, 1)$, we obtain
\begin{equation*}
\big( g^{\mathrm{MAIM}1} \big)'_{s}(d_{s}) = \big( g^{\mathrm{MAIM}1}
\big)''_{s}(d_{s}) = \cdots = \big( g^{\mathrm{MAIM}1} \big)^{(k)}_{s
}(d_{s}) = 0, \big( g^{\mathrm{MAIM}1} \big)^{(k + 1)}_{s}(d_{s}) 
\neq 0.
\end{equation*}

ii) If $k = 2n - 1, n \in \mathbb{N}^{+}$, we have $f_{s}^{\mathrm{
MAIM}} = 0$ and $\dfrac{\mathrm{d}f_{s}^{\mathrm{MAIM}}}{\mathrm{d}
\omega} \neq 0$; then, $(\omega - d_{s})^{k + 1}\mid\mathrm{P}_{num}$
, and $(\omega - d_{s})^{k + 2}\nmid\mathrm{P}_{num}$. Similarly, as
$T_{s}(\omega)Q_{s}(\omega) > 0, \forall \omega \in (0, 1)$, we 
obtain
\begin{equation*}
\big( g^{\mathrm{MAIM}1} \big)'_{s}(d_{s}) = \big( g^{\mathrm{MAIM}1}
 \big)''_{s}(d_{s}) = \cdots = \big( g^{\mathrm{MAIM}1} \big)^{(k + 1
 )}_{s}(d_{s}) = 0, \big( g^{\mathrm{MAIM}1} \big)^{(k + 2)}_{s}(d_{s
 }) \neq 0.
\end{equation*}

(4) According to $P\mathit{1}$ and $P\mathit{3}$ of Property \ref
{property:helperFunctions02} and the fact of $P\mathit{6}$ of 
Property \ref{propertySg2}, from Eq.(\ref
{eq:derivativeOfMappingfunctionWENO-MAIM}), it is very easy to obtain
\begin{equation*}
\begin{array}{ll}
\big(g^{\mathrm{MAIM}1}\big)_{s}^{'}(0^{+}) &= \dfrac{(f_{s}^{\mathrm
{MAIM}})^{2}\cdot(\omega - d_{s})^{2k} + \dfrac{\mathrm{d} f_{s}^{
\mathrm{MAIM}}}{\mathrm{d}\omega}\cdot(\omega - d_{s})^{k + 1}\lim
\limits_{\omega \rightarrow 0^{+}}T_{s}(\omega) + f_{s}^{\mathrm{MAIM
}}\cdot(\omega - d_{s})^{k}\lim\limits_{\omega \rightarrow 0^{+}}P_{s
}(\omega)}{\big(f_{s}^{\mathrm{MAIM}}\cdot(\omega - d_{s})^{k} + \lim
\limits_{\omega \rightarrow 0^{+}}T_{s}(\omega) \big)^{2}} \\
\quad &= \dfrac{(f_{s}^{\mathrm{MAIM}})^{2}\cdot(\omega - d_{s})^{2k}
}{\big(f_{s}^{\mathrm{MAIM}}\cdot(\omega - d_{s})^{k} \big)^{2}} = 1.
\end{array}
\end{equation*}
and
\begin{equation*}
\begin{array}{ll}
\big(g^{\mathrm{MAIM}1}\big)_{s}^{'}(1^{-}) &= \dfrac{(f_{s}^{\mathrm
{MAIM}})^{2}\cdot(\omega - d_{s})^{2k} + \dfrac{\mathrm{d} f_{s}^{
\mathrm{MAIM}}}{\mathrm{d}\omega}\cdot(\omega - d_{s})^{k + 1}\lim
\limits_{\omega \rightarrow 1^{-}}T_{s}(\omega) + f_{s}^{\mathrm{MAIM
}}\cdot(\omega - d_{s})^{k}\lim\limits_{\omega \rightarrow 1^{-}}P_{s
}(\omega)}{\big(f_{s}^{\mathrm{MAIM}}\cdot(\omega - d_{s})^{k} + \lim
\limits_{\omega \rightarrow 1^{-}}T_{s}(\omega) \big)^{2}} \\
\quad &= \dfrac{(f_{s}^{\mathrm{MAIM}})^{2}\cdot(\omega - d_{s})^{2k}
}{\big(f_{s}^{\mathrm{MAIM}}\cdot(\omega - d_{s})^{k} \big)^{2}} = 1.
\end{array}
\end{equation*}
$\hfill\square$

\begin{corollary}
If the adaptive mapping function is defined by Eq.(\ref
{mappingFunctionMAIM}), with the functions $f_{i,s}^{\mathrm{Ada}0}, 
f_{i,s}^{\mathrm{Ada}1}$ calculated by Eq.(\ref{eq:fAda-type2}) or 
Eq.(\ref{eq:fAda-type5}), the mapping function has the properties 
proposed in Theorem \ref{theorem_maimMappingFunction}.
\label{corollary:theorem_maimMappingFunction01}
\end{corollary}

\begin{corollary}
In smooth regions, if the adaptive mapping function is defined by 
Eq.(\ref{mappingFunctionMAIM}), with the functions $f_{i,s}^{\mathrm{
Ada}0}, f_{i,s}^{\mathrm{Ada}1}$ calculated by Eq.(\ref
{eq:fAda-type3}) or Eq.(\ref{eq:fAda-type4}), the mapping function 
has the properties proposed in Theorem \ref
{theorem_maimMappingFunction}.
\label{corollary:theorem_maimMappingFunction02}
\end{corollary}

\begin{remark}
The proof of Corollary \ref{corollary:theorem_maimMappingFunction01} 
is very similar to that of Theorem \ref{theorem_maimMappingFunction}
, especially for $\big(g^{\mathrm{MAIM}5}\big)_{s}(\omega)$ with $f_{
5,s}^{\mathrm{Ada}0},f_{5,s}^{\mathrm{Ada}1}$ defined by Eq.(\ref
{eq:fAda-type5}). For Corollary \ref
{corollary:theorem_maimMappingFunction02}, when in smooth regions, 
we can treat $f_{3,s}^{\mathrm{Ada}0},f_{3,s}^{\mathrm{Ada}1}$ in 
Eq.(\ref{eq:fAda-type3}) and $f_{4,s}^{\mathrm{Ada}0},f_{4,s}^{
\mathrm{Ada}1}$ in Eq.(\ref{eq:fAda-type4}) as constants; then, the 
proof is almost the same as that of $\big(g^{\mathrm{MAIM}5}\big)_{s}
(\omega)$ with $f_{5,s}^{\mathrm{Ada}0},f_{5,s}^{\mathrm{Ada}1}$ 
defined by Eq.(\ref{eq:fAda-type5}).
\label{remark:corrollary:theorem_maimMappingFunction}
\end{remark}

\begin{remark}
Notably, for problems with discontinuities, if the adaptive mapping 
function is defined by Eq.(\ref{mappingFunctionMAIM}), with the 
functions $f_{i,s}^{\mathrm{Ada}0}, f_{i,s}^{\mathrm{Ada}1}$ 
calculated by Eq.(\ref{eq:fAda-type3}) or Eq.(\ref{eq:fAda-type4}), 
the mapping function will only have  property $C1$ proposed in 
Theorem \ref{theorem_maimMappingFunction}, and properties $C2, C3, 
C4$ are not satisfied. We verify these conclusions in subsection \ref
{ParametricStudy}.
\label{remark:note:theorem_maimMappingFunction}
\end{remark}

\subsection{The new WENO schemes}
\label{subsecWENO-MAIM}
Now, we give the new mapped weights as follows
\begin{equation*}
\omega_{s}^{\mathrm{MAIM}i} = \dfrac{\alpha _{s}^{\mathrm{MAIM}i}}{
\sum_{l = 0}^{2} \alpha _{l}^{\mathrm{MAIM}i}}, \alpha_{s}^{\mathrm{
MAIM}i}=\big( g^{\mathrm{MAIM}i} \big)_{s}(\omega^{\mathrm{JS}}_{s}).
\end{equation*}
We denote the new family of modified adaptive improved mapped 
schemes using the weights $\omega_{s}^{\mathrm{MAIM}i}$ as 
WENO-MAIM$1$($k, A,m_{s}$), WENO-MAIM$2$($k, A, Q, \mathrm{CFS}_{s}$
), WENO-MAIM$3$($k, A$), WENO-MAIM$4$($k, A$) and WENO-MAIM$5$\\
($k, A, C$), respectively. For the sake of simplicity, we use 
WENO-MAIM$i$ without causing any confusion.

\subsubsection{The role of $f_{s}^{\mathrm{MAIM}}$}
\label{subsubsecRoleOf_fimFunction}
\begin{figure}[ht]
\centering
 \includegraphics[height=0.33\textwidth]{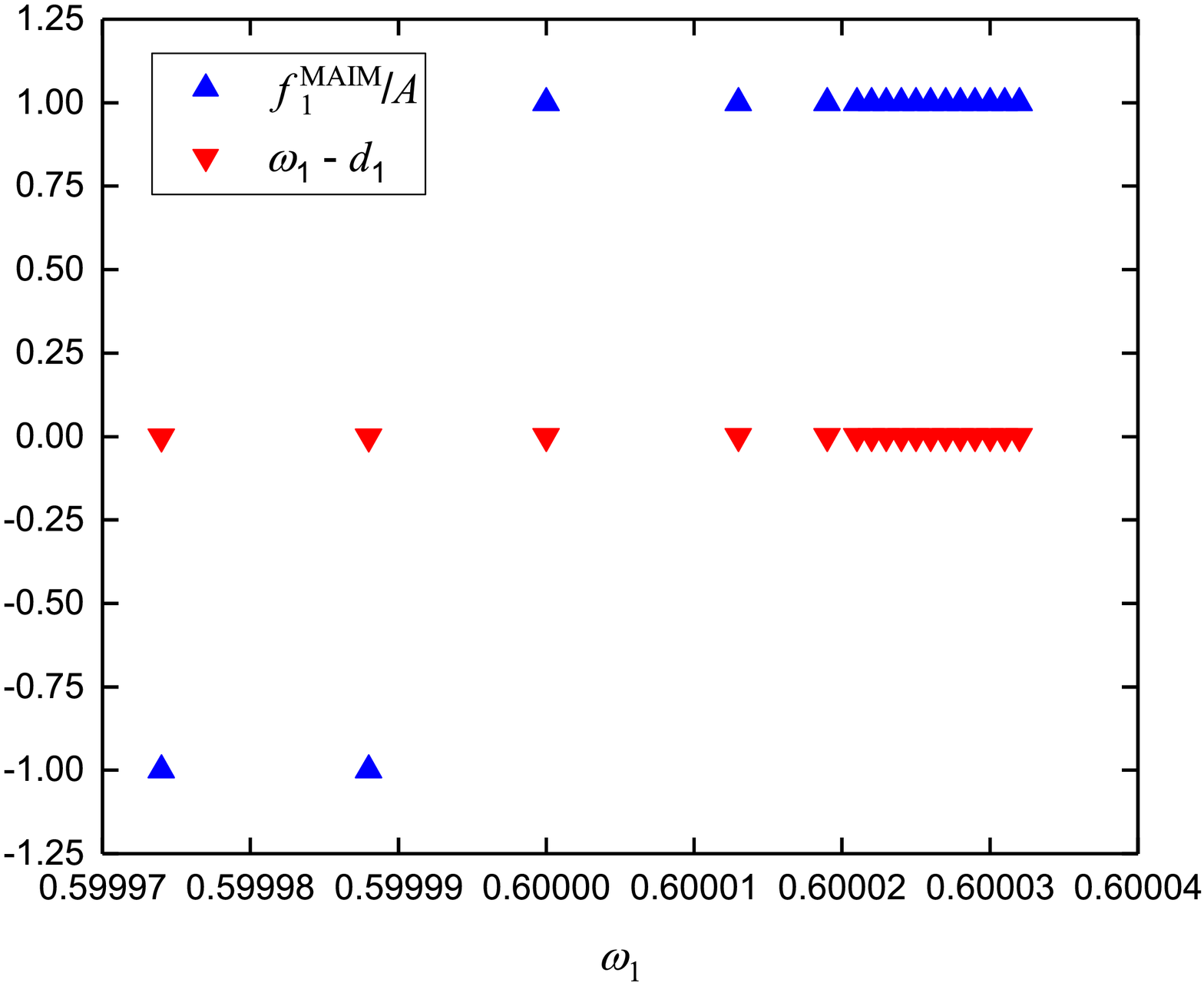}
 \quad
 \includegraphics[height=0.33\textwidth]{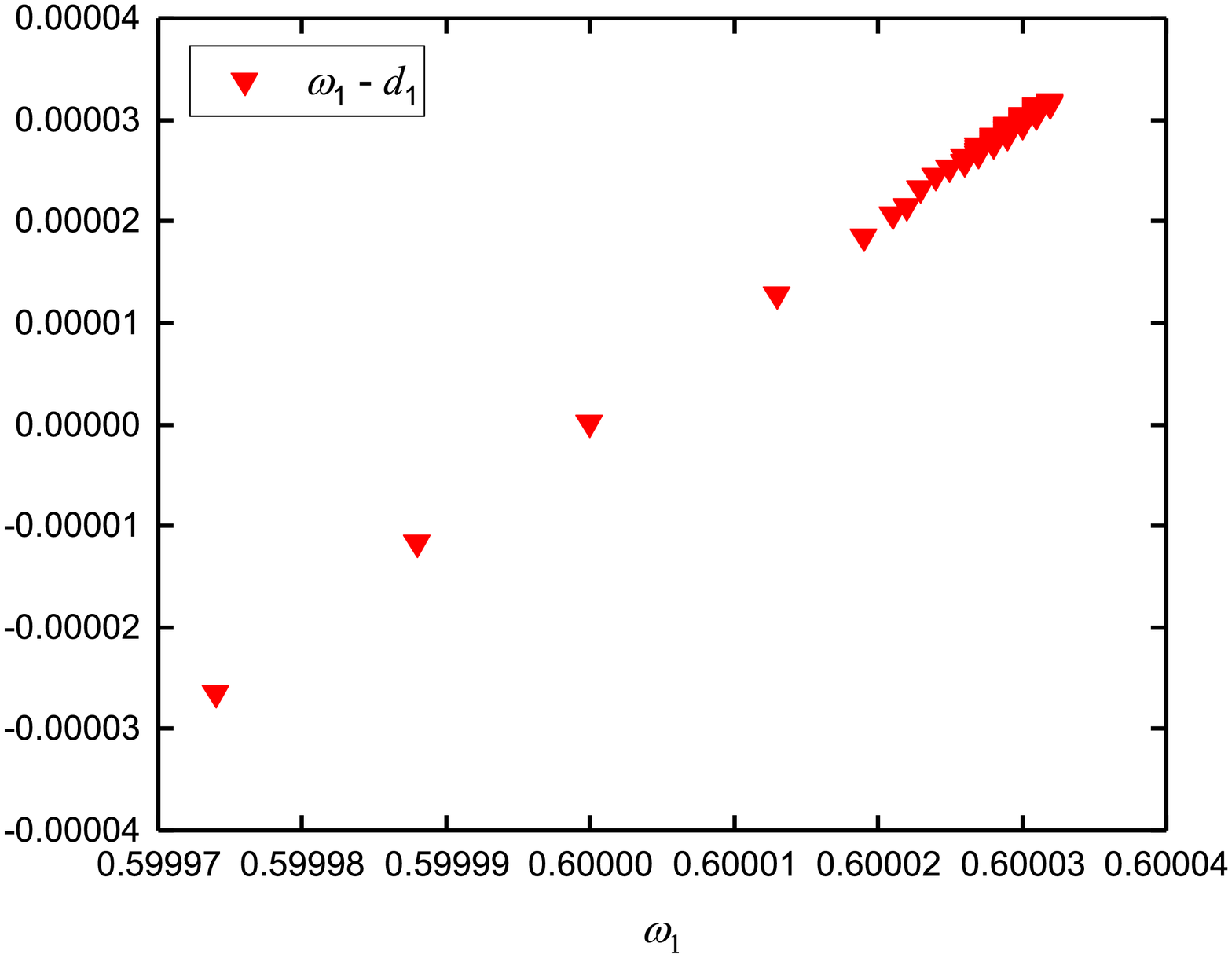}
\caption{Comparison of $\frac{f_{s}^{\mathrm{MAIM}}}{A}$ and $\omega 
- d_{s}$ for $d_{1} = 0.6$ (left: full; right: zoomed).}
\label{fig_roleOfMAIM_d1}
\end{figure}

We now explain that for the case of $k = 2n - 1, n\in \mathbb{N}^{+}$
, the quotient of $f_{s}^{\mathrm{MAIM}}$ and $A$ from Eq.(\ref
{fimFunction}) is only an approximation of the signum function but 
does not tend to $\omega - d_{s}$. In other words, when $k = 2n - 1,
n\in \mathbb{N}^{+}$, the product $f_{s}^{\mathrm{MAIM}} \cdot (
\omega - d_{s})^{k}$ in Eq.(\ref{mappingFunctionMAIM}) is entirely 
different from $A(\omega - d_{s})^{k + 1}$. To illustrate this 
result, we take the following linear advection equation with the 
periodic boundary conditions as an example
\begin{equation}\left\{
\begin{array}{l}
u_{t} + u_{x} = 0, \quad -1 \leq x \leq 1, \\
u(x, 0) = \sin (\pi x),
\end{array}\right.
\label{eq:LAE}
\end{equation}
where the output time is $t = 2.0$ and the number of cells is $N = 
800$. Without loss of generality, we employ the fifth-order 
WENO-MAIM$1$ scheme and set $\delta = 10^{-6}$ (the same holds in 
the rest of this paper), $k = 1, A = 10^{-6}$ for
$f_{s}^{\mathrm{MAIM}}$ and $m_{s} = 0.5, \epsilon_{\mathrm{A}} = 
10^{-10}$ for $f_{1,s}^{\mathrm{Ada}0}, f_{1,s}^{\mathrm{Ada}1}$. 
The comparison of $\frac{f_{s}^{\mathrm{MAIM}}}{A}$ and $\omega - 
d_{s}$, taking $d_{1} = 0.6$ as an example, is shown in 
Fig.\ref{fig_roleOfMAIM_d1}. From Fig.\ref{fig_roleOfMAIM_d1}, we 
intuitively find that $\frac{f_{s}^{\mathrm{MAIM}}}{A}$ is only an 
approximation of the signum function but does not tend to $\omega - 
d_{s}$. The same results have been achieved for $d_{0} = 0.3$ and
$d_{2} = 0.1$.

\subsubsection{Parametric study of the mapping functions}
\label{ParametricStudy}
\begin{figure}[ht]
\centering
	\includegraphics[height=0.39\textwidth]{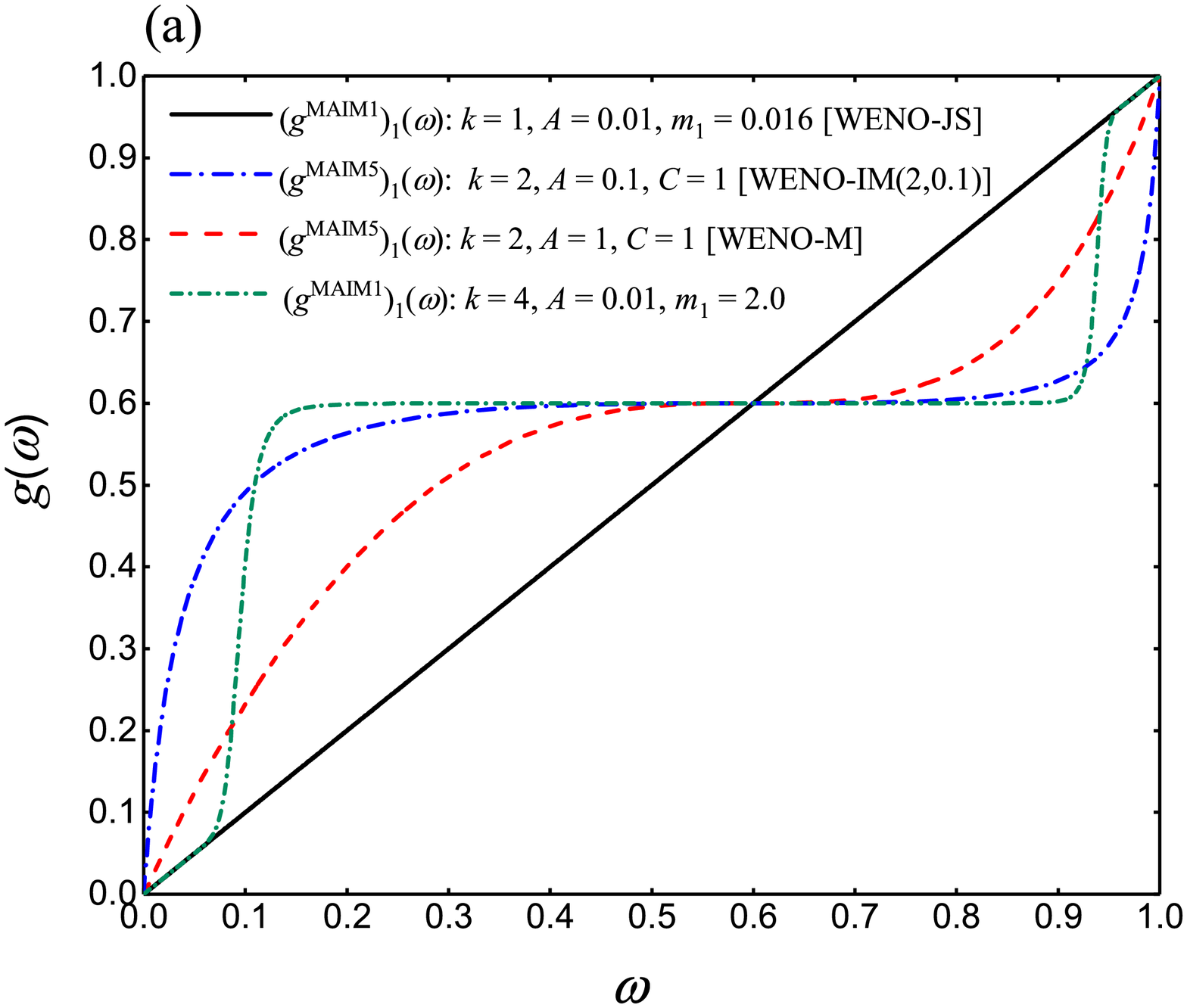}
	\quad
	\includegraphics[height=0.39\textwidth]{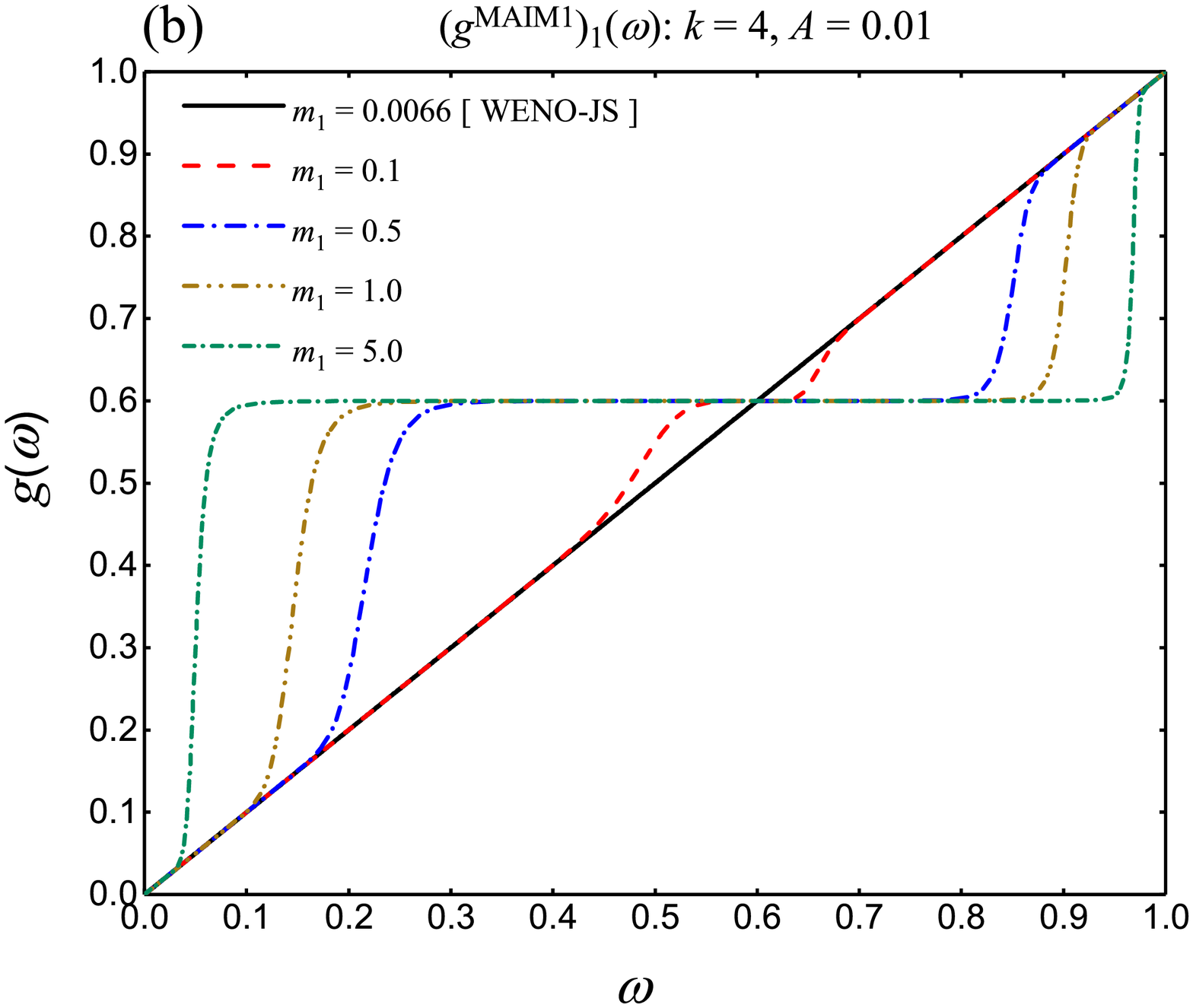} 
	\\
	\includegraphics[height=0.39\textwidth]{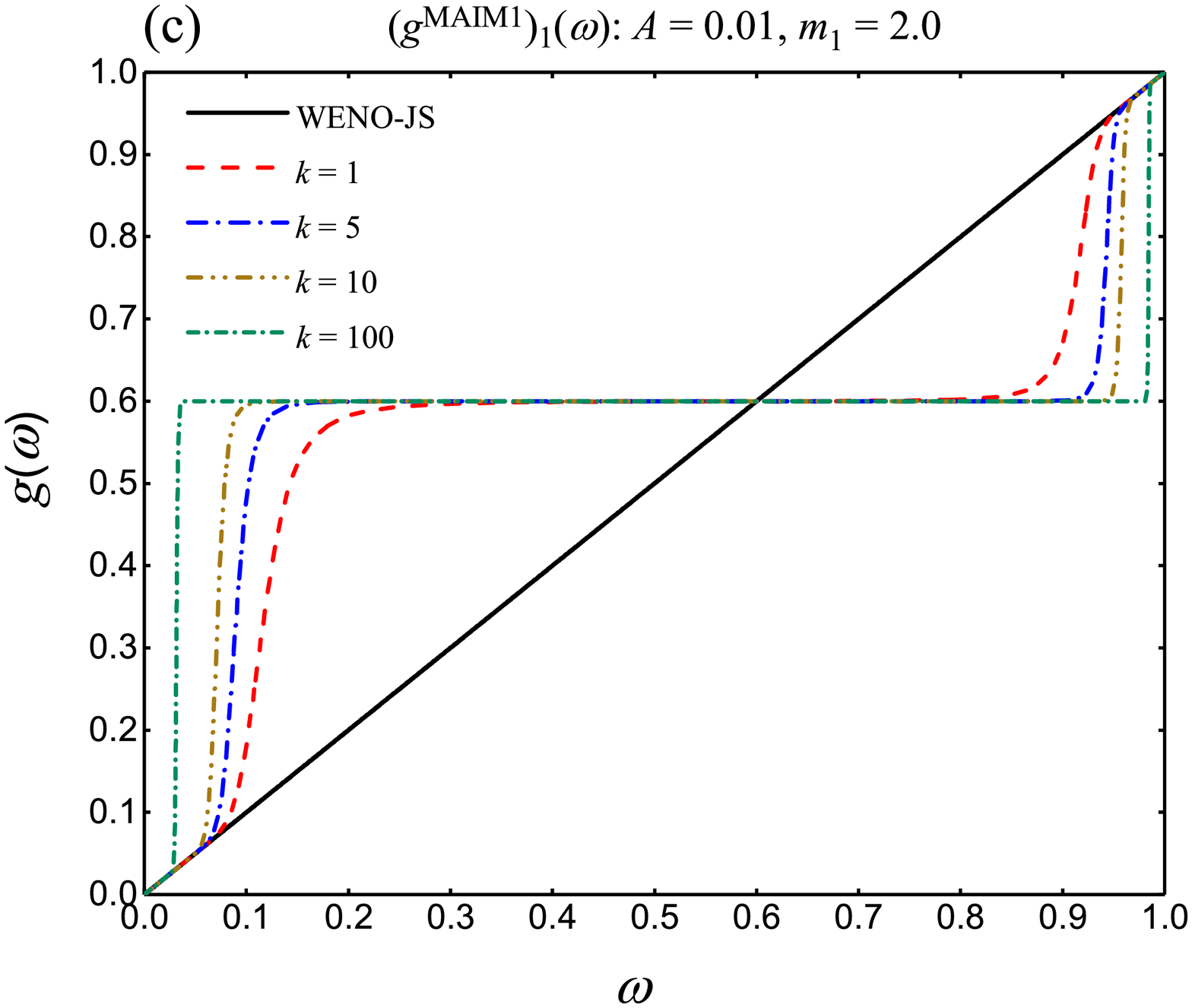}
	\quad
	\includegraphics[height=0.39\textwidth]{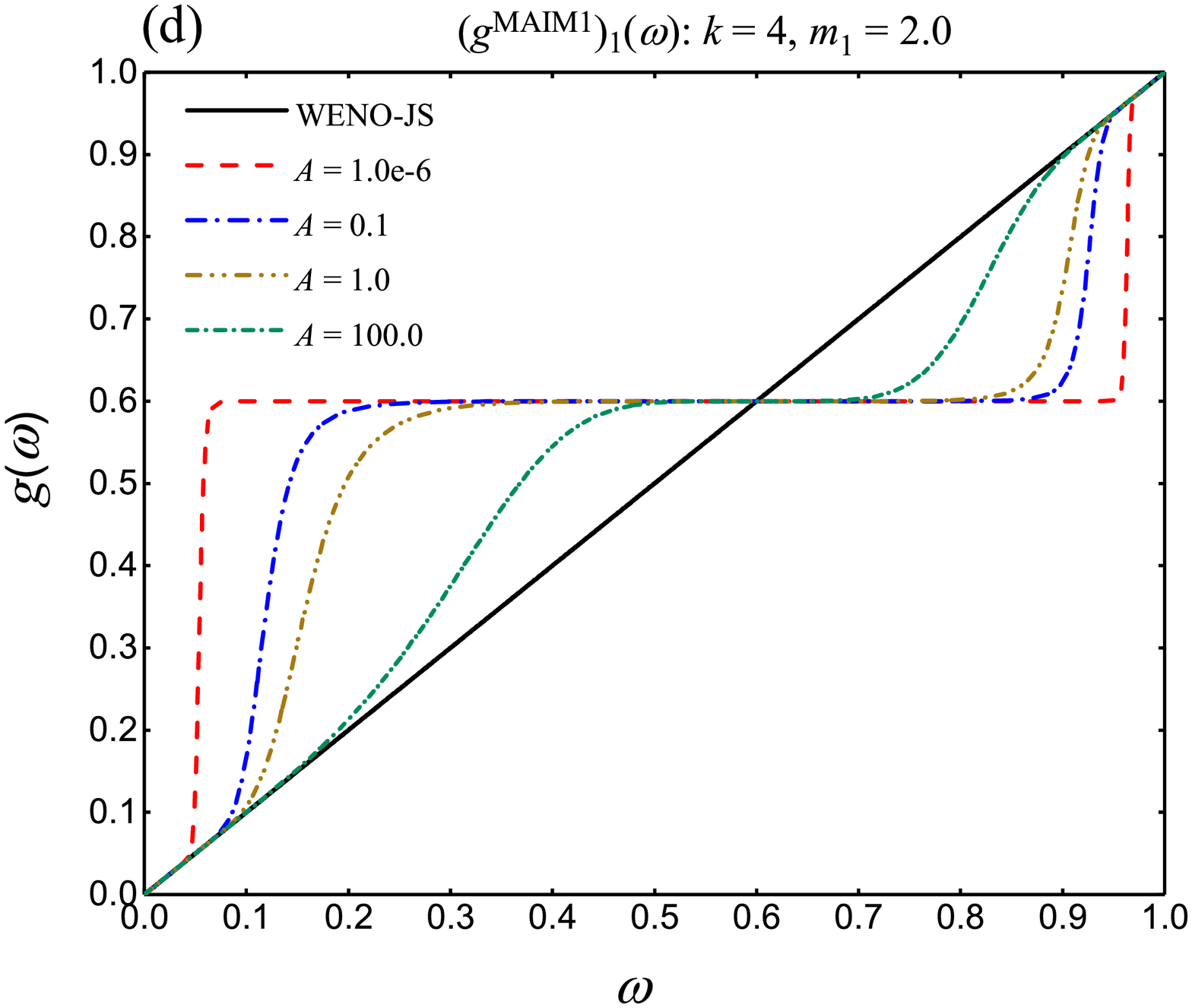}
	   \caption{New mapping functions $\big( g^{\mathrm{MAIM}1} \big)
	   _{1}(\omega)$, $\big( g^{\mathrm{MAIM}5} \big)_{1}(\omega)$, 
	   and effect of varying parameters $k, A$ and $m_{s}$ on $\big( 
	   g^{\mathrm{MAIM}1} \big)_{1}(\omega)$ for $d_{1} = 0.6$.}
	   \label{fig:parametricStudyOfMappingFunctions01}
\end{figure}

We can consider the WENO-JS scheme as a mapped WENO scheme with a 
mapping function defined as $\big( g^{\mathrm{JS}} \big)_{s}(\omega) 
= \omega$. Then, it can be easily verified that $\big( g^{\mathrm{
MAIM}1} \big)_{s}(\omega)\big\lvert_{m_{s} \to 0} = \big( g^{\mathrm{
JS}} \big)_{s}(\omega)$, $\big( g^{\mathrm{MAIM}5} \big)_{s}\big(
\omega\big)\big\lvert_{k=2,A=1,C=1} = \big( g^{\mathrm{M}} \big)_{s}(
\omega)$, $\big( g^{\mathrm{MAIM}5} \big)_{s}\big(\omega\big) \\
\big\lvert_{k=2,A=0.1,C=1} = \big( g^{\mathrm{IM}} \big)_{s}(\omega;
2, 0.1)$. We present these results directly in 
Fig.\ref{fig:parametricStudyOfMappingFunctions01}(a). For 
$\big( g^{\mathrm{MAIM}1} \big)_{s}(\omega)$, we have the following 
properties: (1) for given $k$ and $A$, decreasing $m_{s}$ will make 
the function follow the identity map more closely but to narrow the 
optimal weight interval (see
Fig.\ref{fig:parametricStudyOfMappingFunctions01}(b)), and the 
optimal weight interval stands for the interval about
$\omega = d_{s}$ over which the mapping process attempts to use the 
corresponding optimal weight; (2) for given $A$ and $m_{s}$, the 
optimal weight interval is widened by increasing $k$ as more 
derivatives vanish at $\omega = d_{s}$ (see
Fig.\ref{fig:parametricStudyOfMappingFunctions01}(c)); (3) for given 
$k$ and $m_{s}$, the optimal weight interval is narrowed by 
increasing $A$ (see
Fig.\ref{fig:parametricStudyOfMappingFunctions01}(d)).

\begin{figure}[ht]
\centering
	\includegraphics[height=0.39\textwidth]{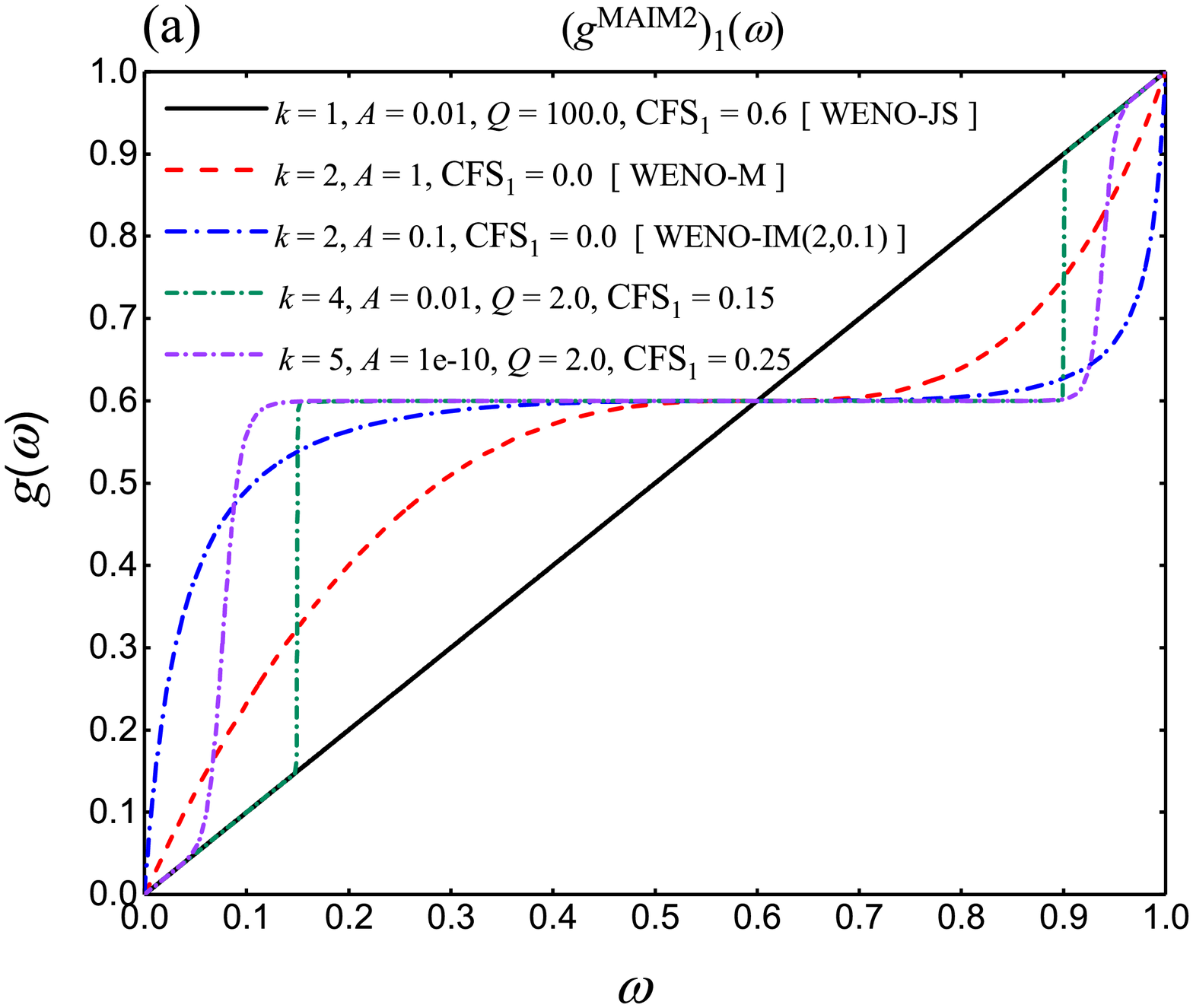}
	\quad
	\includegraphics[height=0.39\textwidth]{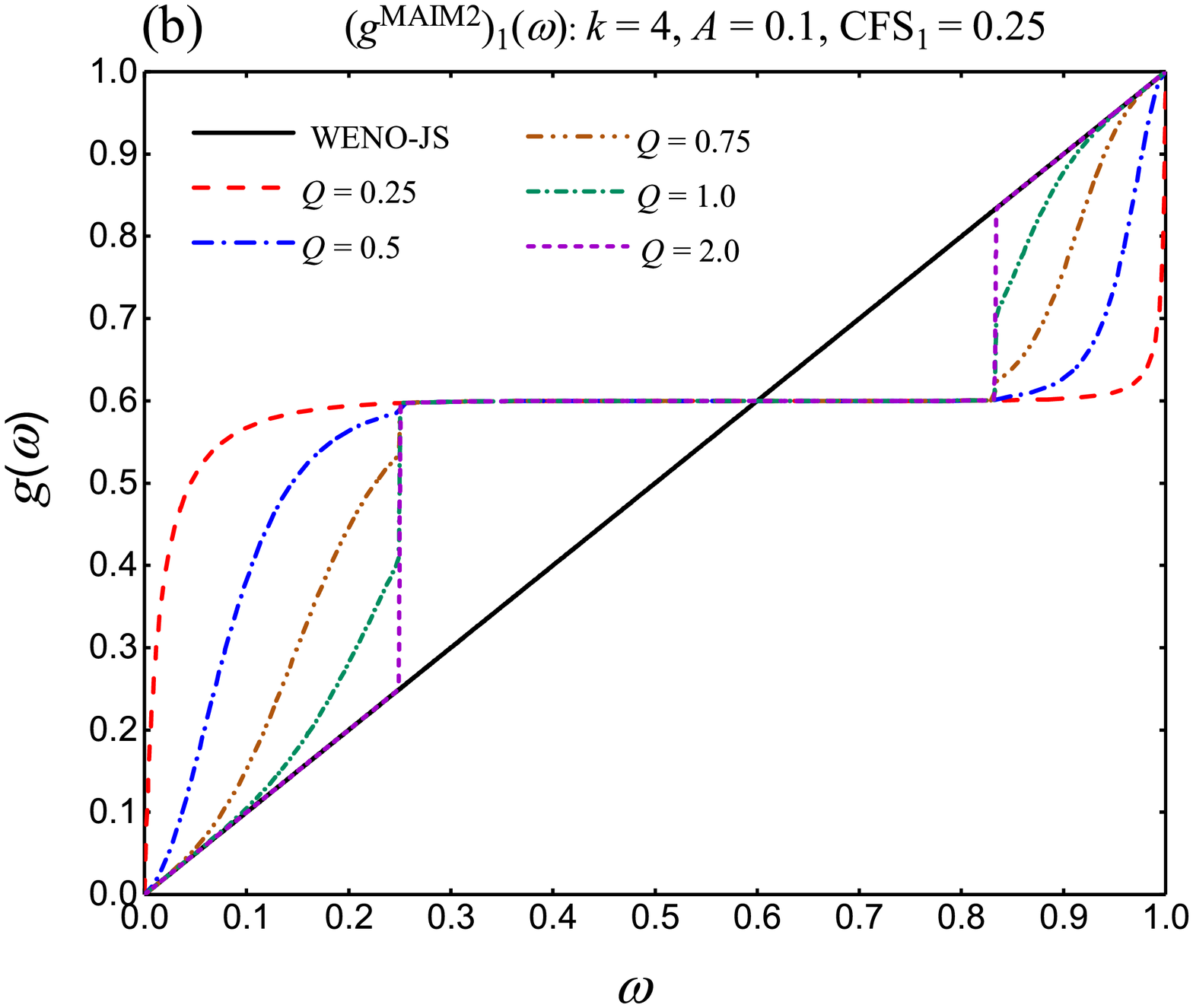} 
	\\
	\includegraphics[height=0.39\textwidth]{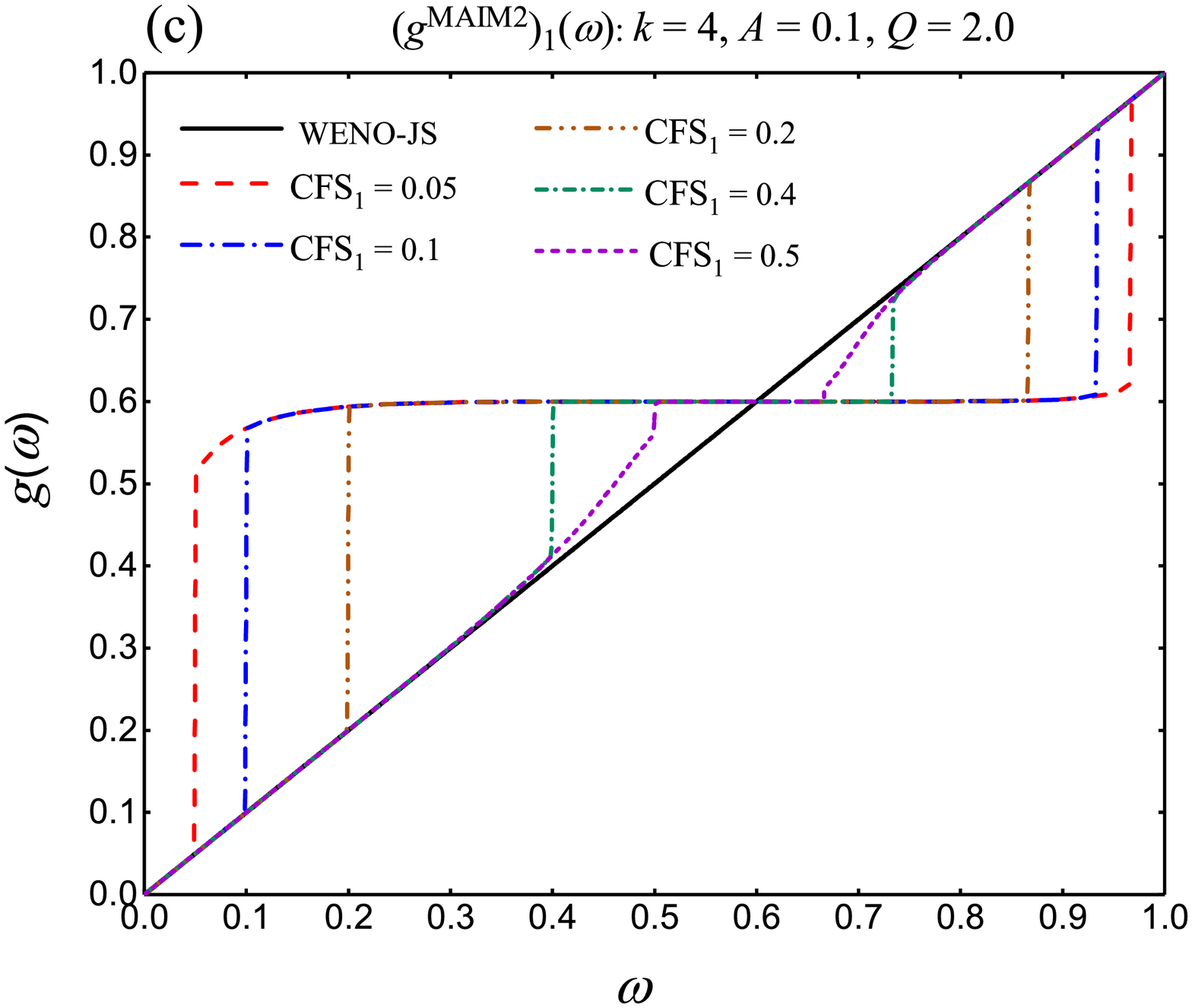}
	\quad
	\includegraphics[height=0.39\textwidth]{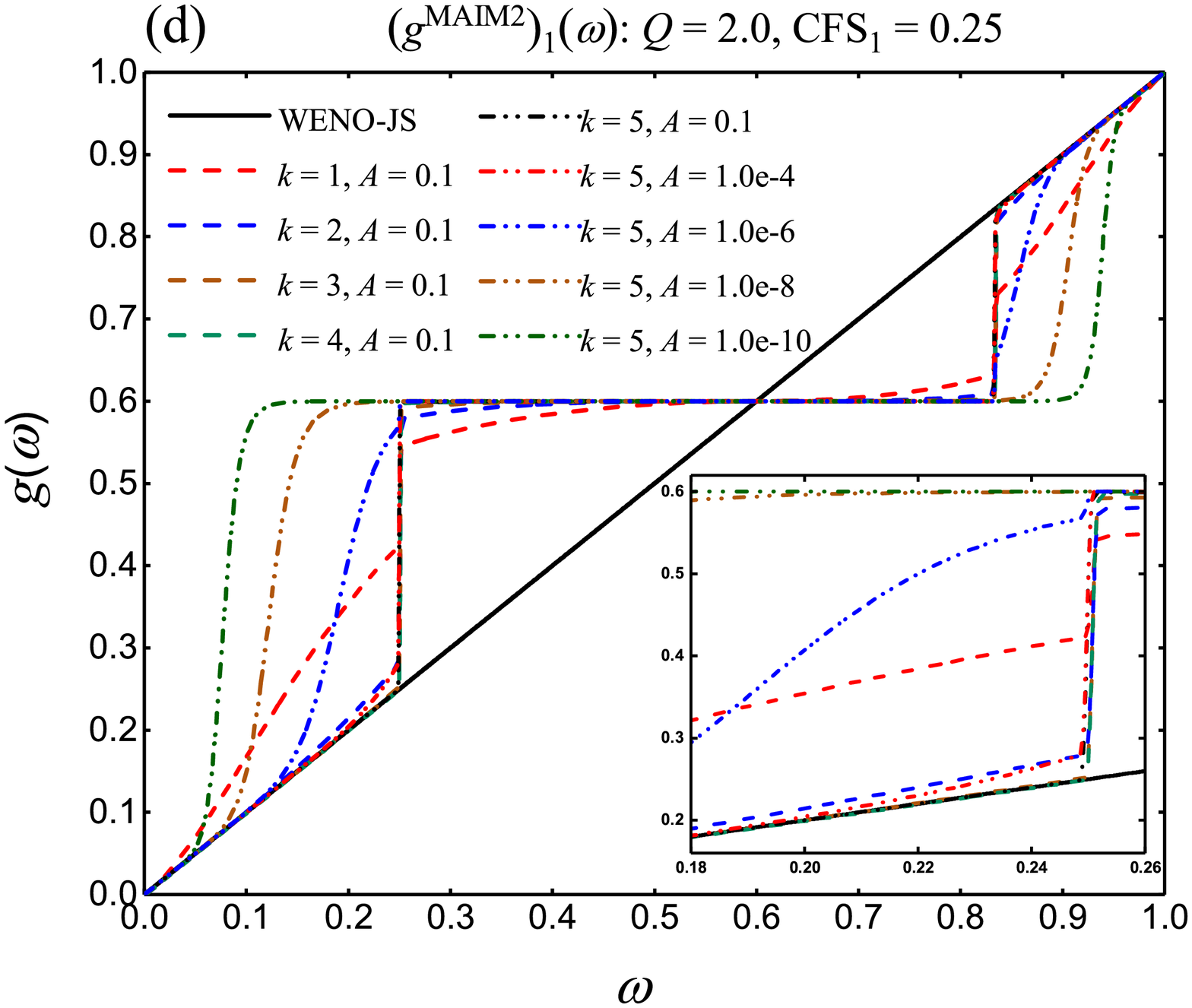}
	   \caption{Effect of varying parameters $k, A, Q$ and $\mathrm{
	   CFS}_{s}$ on $\big( g^{\mathrm{MAIM}2} \big)_{1}(\omega)$ 
	   for $d_{1} = 0.6$.}
	   \label{fig:parametricStudyOfMappingFunctions02}
\end{figure}

For $\big( g^{\mathrm{MAIM}2} \big)_{s}(\omega)$, it is trivial to 
show that $\big( g^{\mathrm{MAIM}2} \big)_{s}\big(\omega\big)\big
\lvert_{k = 1, A = 0.01, Q = 100, \mathrm{CFS}_{s} = 1} = \big( g^{
\mathrm{JS}} \big)_{s}(\omega)$, $\big( g^{\mathrm{MAIM}2} \big)_{s}
\big(\omega\big)\big\lvert_{k = 2, A = 1, \mathrm{CFS}_{s} = 0} \\
= \big( g^{\mathrm{M}} \big)_{s}(\omega)$, $\big( g^{\mathrm{MAIM}2} 
\big)_{s}\big(\omega\big)\big\lvert_{k = 2, A = 0.1, \mathrm{CFS}_{s}
= 0} = \big( g^{\mathrm{IM}} \big)_{s}(\omega; 2, 0.1)$. We 
illustrate these results directly in Fig.\ref
{fig:parametricStudyOfMappingFunctions02}(a). Similarly, for $\big( g
^{\mathrm{MAIM}2} \big)_{s}(\omega)$, we have the following 
properties: (1) for given $k, A$ and $\mathrm{CFS}_{s}$, increasing
$Q$ makes the function follow the identity map more closely and 
narrow the optimal weight interval down to the minimum optimal 
weight interval controlled by the parameter $\mathrm{CFS}_{s}$ (see 
Fig.\ref{fig:parametricStudyOfMappingFunctions02}(b)); (2) for given 
$k, A$ and $Q$, the optimal weight interval is widened by decreasing
$\mathrm{CFS}_{s}$ (see Fig.\ref
{fig:parametricStudyOfMappingFunctions02}(c)); (3) the effects of 
parameters $k$ and $A$ are the same as in $\big( g^{\mathrm{MAIM}1} 
\big)_{s}(\omega)$, but one must consider the effect of the 
parameter $\mathrm{CFS}_{s}$ (see Fig.\ref
{fig:parametricStudyOfMappingFunctions02}(d)). 

The study of the mapping functions $\big( g^{\mathrm{MAIM}3} \big)_{s
}(\omega)$ and $\big( g^{\mathrm{MAIM}4} \big)_{s}(\omega)$ is 
slightly more complicated as these functions are not independent 
functions of $\omega$: one must obtain the relationship of $\big( g^{
\mathrm{MAIM}3} \big)_{s}(\omega) \sim \omega$ or $\big( g^{\mathrm{
MAIM}4} \big)_{s}(\omega) \sim \omega$ by specific numerical 
examples. We found that if the problem is non-smooth, these mapping 
functions will not maintain the monotonicity property. To illustrate 
this result, we take the linear advection equation $u_t + u_x = 0$
as an example with the following initial condition consisting of two 
constant states separated by sharp discontinuities at $x = 0, \pm 1$

\begin{equation}
u(x, 0) = \left\{
\begin{array}{ll}
1, & x \in [-1, 0), \\
0, & x \in [0, 1].
\end{array}
\right.
\label{eq:LAE:IC2}
\end{equation}

\begin{figure}[ht]
\centering
	\includegraphics[height=0.39\textwidth]{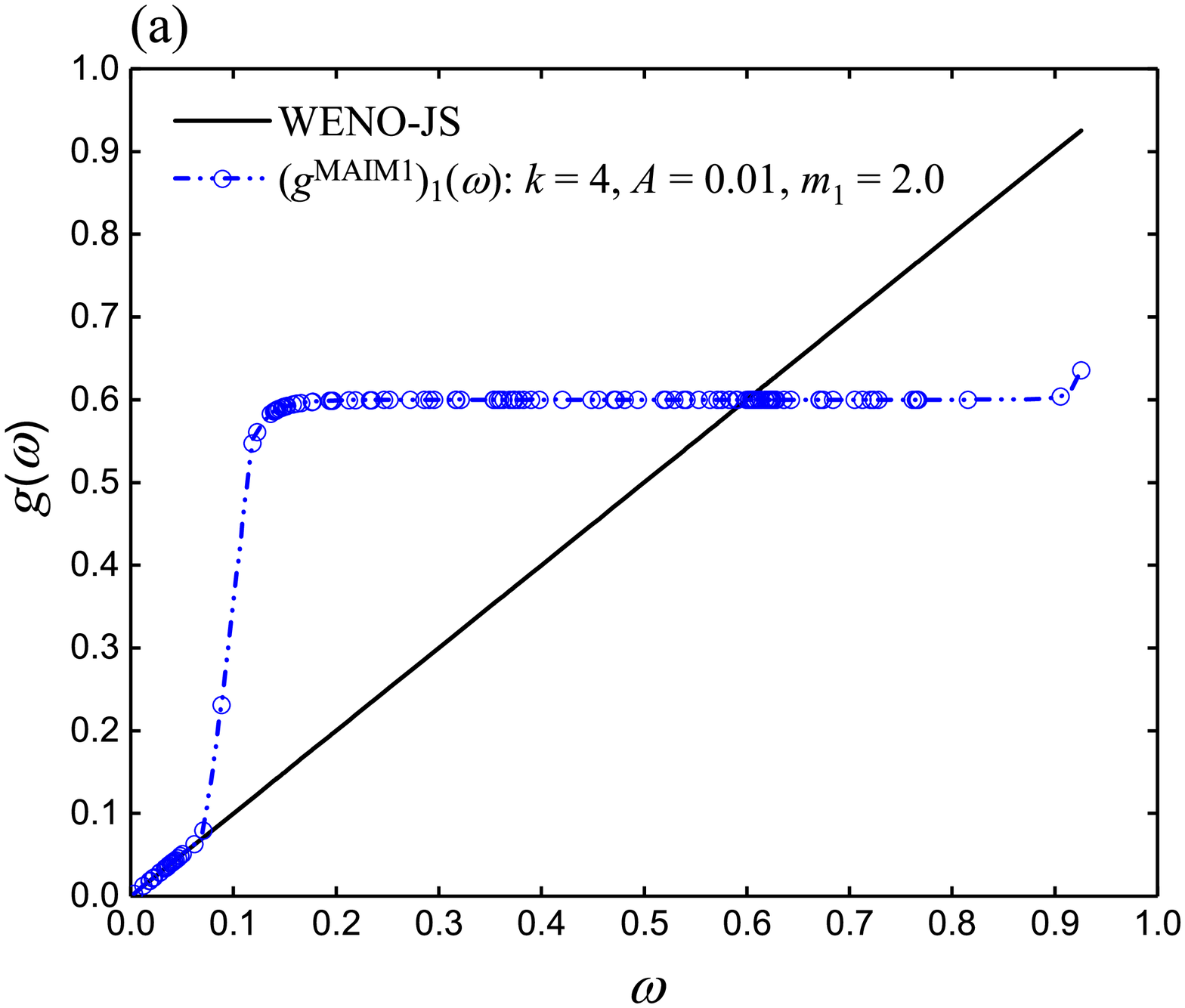}
	\quad
	\includegraphics[height=0.39\textwidth]{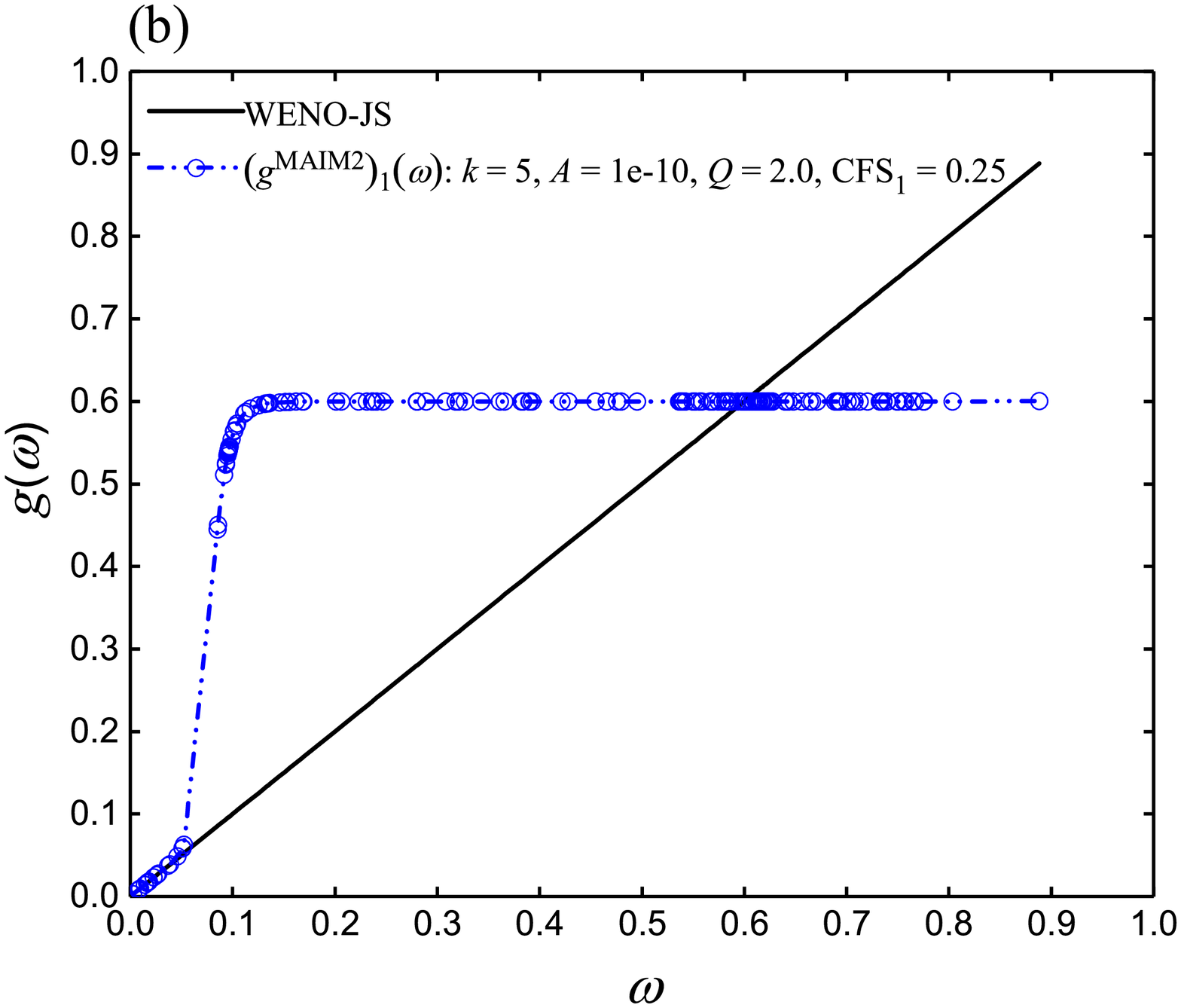} 
	\\
	\includegraphics[height=0.39\textwidth]{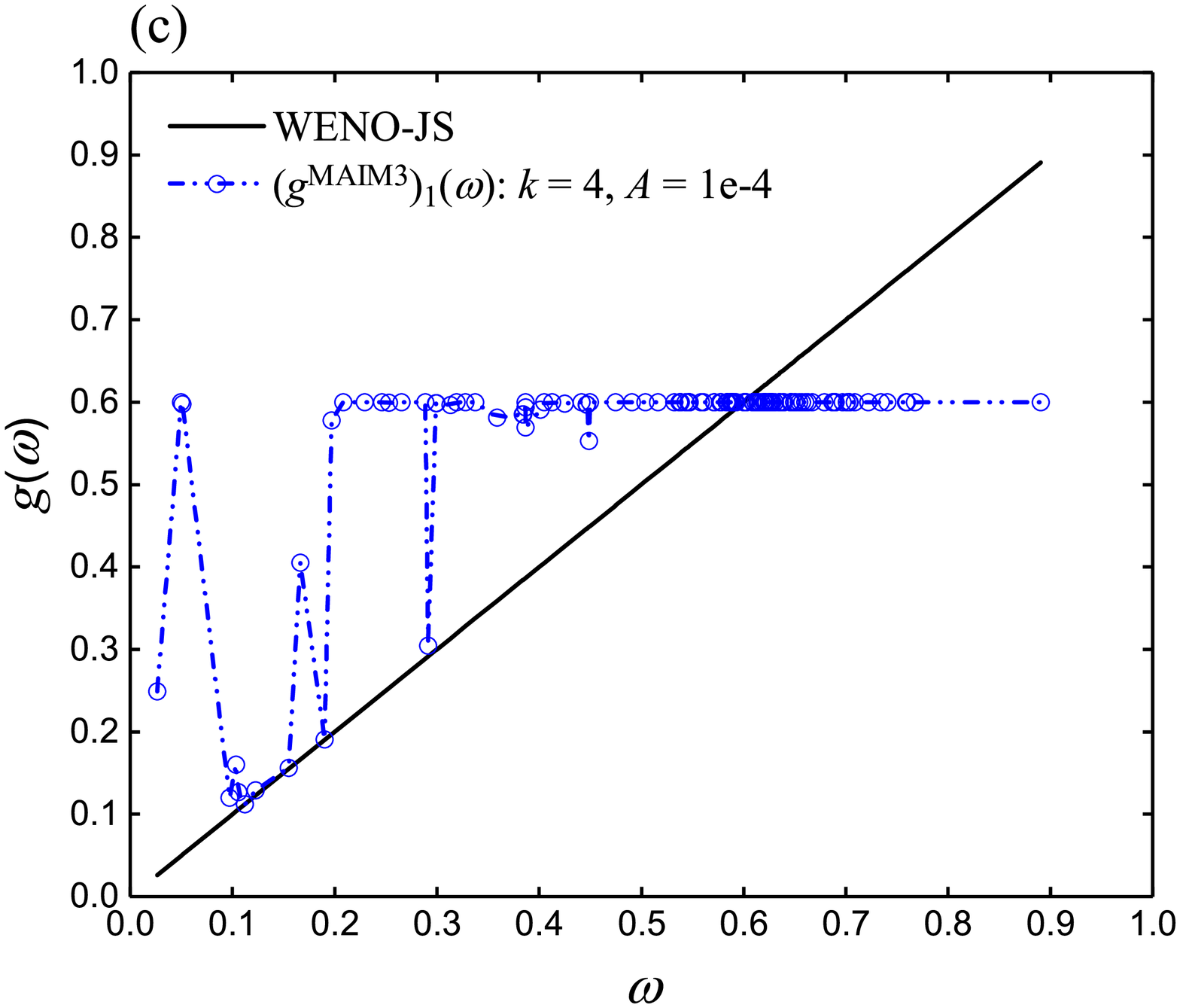}
	\quad
	\includegraphics[height=0.39\textwidth]{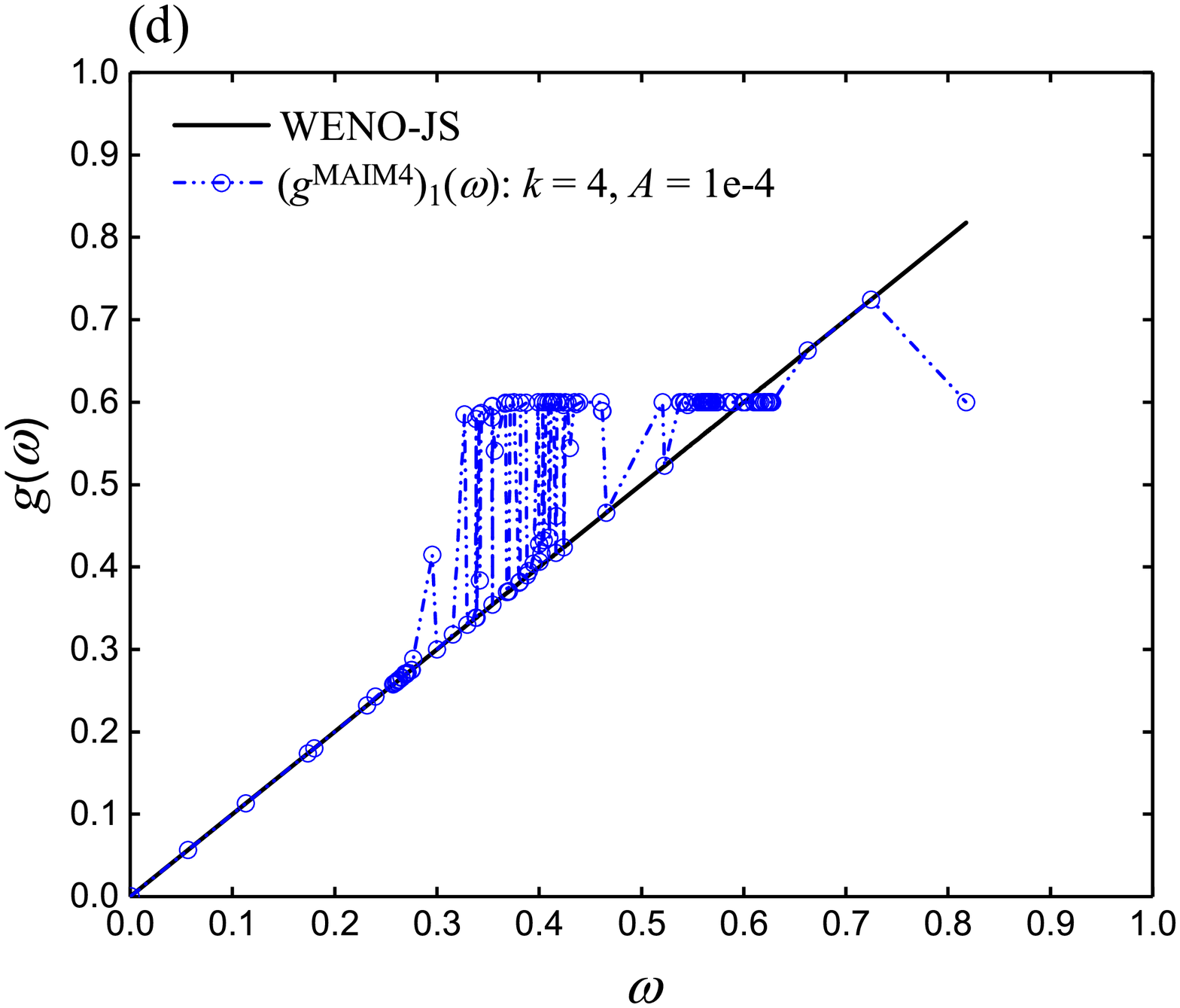}
	   \caption{Comparison of the mapping functions $\big( g^{\mathrm
	   {MAIM}1} \big)_{1}(\omega)$, $\big( g^{\mathrm{MAIM}2} \big)_{
	   1}(\omega)$, $\big( g^{\mathrm{MAIM}3} \big)_{1}(\omega)$, $
	   \big( g^{\mathrm{MAIM}4} \big)_{1}(\omega)$ and $\big( g^{
	   \mathrm{AIM}} \big)_{1}(\omega; 4, 2, 1\mathrm{e}4)$ for $d_{1
	   } = 0.6$ for linear advection equation (\ref{eq:LAE}) with 
	   initial condition (\ref{eq:LAE:IC2}).}
	   \label{fig:parametricStudyOfMappingFunctions03}
\end{figure}

The results are shown in Fig.\ref
{fig:parametricStudyOfMappingFunctions03} at the resolution $N = 200$
cells and output time $t = 200$. In the calculations, the CFL number 
is chosen to be $0.1$ and the periodic boundary conditions in two 
directions are used. From Fig.\ref
{fig:parametricStudyOfMappingFunctions03}, we see that the mapping 
functions $\big( g^{\mathrm{MAIM}3} \big)_{s}(\omega)$ and $\big( g^{
\mathrm{MAIM}4} \big)_{s}(\omega)$ are non-monotonic with respect to
$\omega$ for non-smooth problems. However, the monotonicity of the 
mapping functions $\big( g^{\mathrm{MAIM}1} \big)_{s}(\omega)$ and $
\big( g^{\mathrm{MAIM}2} \big)_{s}(\omega)$ with respect to $\omega$ 
can be maintained very well.

\subsubsection{The rate of convergence}
\label{subsubsecConvergence_WENO-MAIM}
WENO schemes up to $17$th-order of accuracy have been outlined by 
Gerolymos et al. in \cite{veryHighOrderWENO}. Notationally, let $n_{
\mathrm{cp}} \in \mathbb{N}$ denote the order of the critical point; 
that is, if $n_{\mathrm{cp}} = k_{\mathrm{crit}}$, then $f' = \cdots 
= f^{(k_{\mathrm{crit}})} = 0$ and $f^{(k_{\mathrm{crit}} + 1)} \neq 
0$. Now, we present Theorem \ref{theoremWENO-MAIM}, which will show 
that the WENO-MAIM$i$ schemes can recover the optimal convergence 
rates by setting the parameter $k$ for different values of $n_{
\mathrm{cp}}$ in smooth regions.
\begin{theorem}
Let $\lceil x \rceil$ be a ceiling function of $x$, for $n_{\mathrm{cp}} < r - 1$, the WENO-MAIM$i$ schemes can achieve 
the optimal $(2r - 1)$th-order of accuracy if the modified adaptive 
improved mapping functions $(g^{\mathrm{MAIM}i})_{s}(\omega)$ 
defined in Eq.(\ref{mappingFunctionMAIM}) are applied to the 
original weights in the WENO-JS scheme with $k \geq k^{\mathrm
{MAIM}}$, where
\begin{equation}
k^{\mathrm{MAIM}} = \Bigg\lceil \dfrac{r}{r - 1 - n_{\mathrm{cp}}} - 
2 \Bigg\rceil + \dfrac{1 + \Big( -1 \Big)^{\Bigg\lceil \dfrac{r}{r - 
1 - n_{\mathrm{cp}}} -2 \Bigg\rceil}}{2}.
\end{equation}
\label{theoremWENO-MAIM}
\end{theorem}

Before giving the proof of Theorem \ref{theoremWENO-MAIM}, we state 
the following two lemmas.
\begin{lemma}
When $n_{\mathrm{cp}} < r - 1$, the weights $\omega_{s}^{\mathrm{JS
}}$ in the $(2r - 1)$th-order WENO-JS scheme satisfy
\begin{equation}
\omega_{s}^{\mathrm{JS}} - d_{s} = O(\Delta x^{r - 1 - n_{\mathrm{cp}
}}), \quad r = 2, 3, \cdots, 9,
\label{Eq1_lemma_ncpWENO-JS}
\end{equation} 
then, the convergence order is
\begin{equation*}
r_{\mathrm{c}} = \left\{
\begin{array}{ll}
2r - 1, & \mathrm{if} \quad n_{\mathrm{cp}} = 0, \\
2r - 2 - n_{\mathrm{cp}}, & \mathrm{if} \quad n_{\mathrm{cp}} = 1, 2
, \cdots, r - 1.
\end{array}\right.
\end{equation*} 
\label{lemma_ncpWENO-JS}
\end{lemma}

The above lemma is a direct extension of Lemma 2 in \cite{WENO-IM}, 
and we only extend the range of  $r$ from $r = 2, 3, \cdots, 6$ to $r
 = 2, 3, \cdots, 9$. As the proof of the lemma can be found on page 
 565 in \cite{WENO-M}, we simply state it explicitly here.

\begin{lemma}
For $r = 2,3,\cdots,9$, the sufficient condition for the $(2r - 1)$
th-order WENO scheme to achieve the optimal order of accuracy is
\begin{equation*}
\omega_{s} - d_{s} = O(\Delta x^{r}), \quad s = 0, 1, \cdots, r-1.
\end{equation*} 
\label{lemma_SufficientCondition}
\end{lemma}

\begin{remark}
The proof of Lemma \ref{lemma_SufficientCondition} in this paper is 
the finite volume version of Lemma 1 in \cite{WENO-IM}, which is a 
natural extension of the one in \cite{WENO-M}. For more details
about the proof within the finite volume framework, the reference 
\cite{WENO-accuracy_fv-WENO} is refered to.
\label{remark_proof_lemma_SufficientCondition}
\end{remark}
\textbf{Proof of Theorem \ref{theoremWENO-MAIM}.} 

(1) If $k = 2n, n\in \mathbb{N}^{+}$, from Theorem \ref
{theorem_maimMappingFunction}, Corollary \ref
{corollary:theorem_maimMappingFunction01} and Corollary \ref
{corollary:theorem_maimMappingFunction02}, we have $\big( g^{\mathrm{
MAIM}i} \big)_{s}(d_{s}) = d_{s}$ and $\big( g^{\mathrm{MAIM}i} \big)
'_{s}(d_{s}) = \big( g^{\mathrm{MAIM}i} \big)''_{s}(d_{s}) = \cdots =
 \big( g^{\mathrm{MAIM}i} \big)^{(k)}_{s}(d_{s}) = 0, \big( g^{
 \mathrm{MAIM}i} \big)^{(k + 1)}_{s}(d_{s}) \neq 0.$ Thus, 
 evaluation at $\omega_{s}^{\mathrm{JS}}$ of the Taylor series 
 approximation of the $\big( g^{\mathrm{MAIM}i} \big)_{s}(\omega)$  
 about $d_{s}$ yields
\begin{equation*}
\begin{aligned}
\alpha_{s}^{\mathrm{MAIM}i} &= \big( g^{\mathrm{MAIM}i} \big)_{s}(d_{
s}) + \sum_{l=1}^{k}\dfrac{\big( g^{\mathrm{MAIM}i} \big)^{(l)}_{s}(d
_{s})}{l!}\Big( \omega_{s}^{\mathrm{JS}} - d_{s} \Big)^{l} + \dfrac{
\big( g^{\mathrm{MAIM}i} \big)^{(k + 1)}_{s}(d_{s})}{(k + 1)!}\Big( 
\omega_{s}^{\mathrm{JS}} - d_{s} \Big)^{k + 1} + \cdots \\
&= d_{s} + \dfrac{\big( g^{\mathrm{MAIM}i} \big)^{(k + 1)}_{s}(d_{s})
}{(k + 1)!}\Big( \omega_{s}^{\mathrm{JS}} - d_{s} \Big)^{k + 1} + 
\cdots
\end{aligned}
\end{equation*}
According to Lemma \ref{lemma_ncpWENO-JS} and Lemma \ref
{lemma_SufficientCondition}, the $(2r - 1)$th-order WENO-MAIM$i$ 
schemes can achieve the optimal $(2r - 1)$th-order of accuracy with 
requirement 
\begin{equation*} 
(r - 1 - n_{\mathrm{cp}})(k + 1) \geq r, \quad 1 \leq n_{\mathrm{cp}}
 < r - 1.
\end{equation*} 
As $k = 2n, n\in \mathbb{N}^{+}$, we can rewrite the requirement 
above as $k \geq k_{(2n)}^{\mathrm{MAIM}}$ by introducing 
\begin{equation*} 
k_{(2n)}^{\mathrm{MAIM}} = \Bigg\lceil \dfrac{r}{r - 1 - n_{\mathrm{
cp}}} - 1 \Bigg\rceil + \dfrac{1 - \Big( -1 \Big)^{\Bigg\lceil \dfrac
{r}{r - 1 - n_{\mathrm{cp}}} - 1 \Bigg\rceil}}{2}.
\end{equation*} 

(2) If $k = 2n - 1, n\in \mathbb{N}^{+}$, from Theorem \ref
{theorem_maimMappingFunction}, Corollary \ref
{corollary:theorem_maimMappingFunction01} and Corollary \ref
{corollary:theorem_maimMappingFunction02}, we have $\big( g^{\mathrm{
MAIM}i} \big)_{s}(d_{s}) = d_{s}$ and $\big( g^{\mathrm{MAIM}i} \big)
'_{s}(d_{s}) = \big( g^{\mathrm{MAIM}i} \big)''_{s}(d_{s}) = \cdots =
 \big( g^{\mathrm{MAIM}i} \big)^{(k + 1)}_{s}(d_{s}) = 0, \big( g^{
 \mathrm{MAIM}i} \big)^{(k + 2)}_{s}(d_{s}) \neq 0.$ Therefore, 
 similarly, evaluation at $\omega_{s}^{\mathrm{JS}}$ of the Taylor 
 series approximation of the $\big( g^{\mathrm{MAIM}i} \big)_{s}(
 \omega)$  about $d_{s}$ yields
\begin{equation*}
\begin{aligned}
\alpha_{s}^{\mathrm{MAIM}i} &= \big( g^{\mathrm{MAIM}i} \big)_{s}(d_{
s}) + \sum_{l=1}^{k+1}\dfrac{\big( g^{\mathrm{MAIM}i} \big)^{(l)}_{s}
(d_{s})}{l!}\Big( \omega_{s}^{\mathrm{JS}} - d_{s} \Big)^{l} + \dfrac
{\big( g^{\mathrm{MAIM}i} \big)^{(k + 2)}_{s}(d_{s})}{(k + 2)!}\Big( 
\omega_{s}^{\mathrm{JS}} - d_{s} \Big)^{k + 2} + \cdots \\
&= d_{s} + \dfrac{\big( g^{\mathrm{MAIM}i} \big)^{(k + 2)}_{s}(d_{s})
}{(k + 2)!}\Big( \omega_{s}^{\mathrm{JS}} - d_{s} \Big)^{k + 2} + 
\cdots
\end{aligned}
\end{equation*}
Furthermore, according to Lemma \ref{lemma_ncpWENO-JS} and Lemma \ref
{lemma_SufficientCondition}, the $(2r - 1)$th-order WENO-MAIM$i$ 
schemes can achieve the optimal $(2r - 1)$th-order of accuracy with 
requirement
\begin{equation*}
(r - 1 - n_{\mathrm{cp}})(k + 2) \geq r, \quad 1 \leq n_{\mathrm{cp}}
 < r - 1.
\end{equation*}
As $k = 2n - 1, n\in \mathbb{N}^{+}$, we can rewrite the requirement 
above as $k \geq k_{(2n - 1)}^{\mathrm{MAIM}}$ by introducing
\begin{equation*}
k_{(2n - 1)}^{\mathrm{MAIM}} = \Bigg\lceil \dfrac{r}{r - 1 - n_{
\mathrm{cp}}} - 2 \Bigg\rceil + \dfrac{1 + \Big( -1 \Big)^{\Bigg
\lceil \dfrac{r}{r - 1 - n_{\mathrm{cp}}} - 2 \Bigg\rceil}}{2}.
\end{equation*}

Clearly, we only need to ensure $k \geq k^{\mathrm{MAIM}}$, with $k^{
\mathrm{MAIM}} = \min \Big(k_{(2n - 1)}^{\mathrm{MAIM}}, k_{(2n)}^{
\mathrm{MAIM}} \Big)$, to achieve the optimal $(2r - 1)$th-order of
 accuracy. It is easy to check that $k_{(2n - 1)}^{\mathrm{MAIM}} 
 \leq k_{(2n)}^{\mathrm{MAIM}}$. Eventually, we obtain
\begin{equation}
k^{\mathrm{MAIM}} = \Bigg\lceil \dfrac{r}{r - 1 - n_{\mathrm{cp}}} - 
2 \Bigg\rceil + \dfrac{1 + \Big( -1 \Big)^{\Bigg\lceil \dfrac{r}{r - 
1 - n_{\mathrm{cp}}} -2 \Bigg\rceil}}{2}.
\end{equation}
$\hfill\square$ \\
\begin{remark}
For different values of $n_{\mathrm{cp}}$, we can calculate the 
detailed convergence order of the WENO-MAIM$i$ schemes using Theorem 
\ref{theoremWENO-MAIM}. Table \ref{table_convergenceOrder} is a 
comparison between WENO-JS, WENO-M($1$), WENO-M($2$), WENO-IM($k, A$
) and WENO-MAIM$i$. Note that WENO-M($n$) stands for WENO-M schemes 
with $n$ times mapping, and their rates of convergence are related 
to $r$ and $n_{\mathrm{cp}}$. A general conclusion of the rates of 
convergence of WENO-M($n$) schemes is proposed in  Lemma \ref
{Lemma_convergenceRateofWENO-M}, and we provide the proof in 
Appendix B. 
\end{remark}

\begin{lemma}
If $n$ mapping is used in the $(2r - 1)$th-order WENO-M scheme, then
for different values of $n_{\mathrm{cp}}$, the weights $\omega_{s}^{
\mathrm{M}}$ in the $(2r - 1)$th-order WENO-M scheme satisfy
\begin{equation*}
\omega_{s}^{\mathrm{M}} - d_{s} = O\Big( (\Delta x)^{3^{n}\times(r - 
1 - n_{\mathrm{cp}})} \Big), \quad r = 2, 3, \cdots, 9, \quad n_{
\mathrm{cp}} = 0, 1, \cdots, r - 1,
\end{equation*}
and the rate of convergence is
\begin{equation*}
r_{\mathrm{c}} = \left\{
\begin{array}{ll}
2r - 1,  &\mathrm{if} \quad n_{\mathrm{cp}} = 0, \cdots, \Bigg\lfloor
\dfrac{3^{n}-1}{3^{n}}r - 1 \Bigg\rfloor, \\
(3^{n}+1)(r - 1) - 3^{n}\times n_{\mathrm{cp}}, &\mathrm{if} \quad n_
{\mathrm{cp}} = \Bigg\lfloor \dfrac{3^{n}-1}{3^{n}}r - 1 \Bigg\rfloor
+ 1, \cdots, r - 1.
\end{array}\right.
\end{equation*}
where $\lfloor x \rfloor$ is a floor function of $x$.
\label{Lemma_convergenceRateofWENO-M}
\end{lemma}


%% file: article_numerics.tex
\section{Numerical experiments}
\label{secNumericalExperiments}
In this section, we present the results of the modified adaptive 
improved mapped WENO-MAIM$i$ schemes in comparison with the WENO-JS, 
WENO-M and WENO-IM($2, 0.1$) schemes, and all these schemes are 
in the finite volume version with $\epsilon = 10^{-40}$. The global 
Lax-Friedrichs numerical flux is used in the following numerical 
experiments. The $L_{1}$-, $L_{2}$- and $L_{\infty}$- norm of the 
error are calculated by comparing the numerical solution 
$(u_{h})_{j}$ with the exact solution $u_{j}^{\mathrm{exact}}$ 
according to 
\begin{equation*}
\displaystyle
\begin{array}{l}
L_{1}(h) = \displaystyle\sum\limits_{j} h \cdot \big\lvert u_{j}^{
\mathrm{exact}} - (u_{h})_{j} \big\rvert, \\
L_{2}(h) = \sqrt{\displaystyle\sum\limits_{j} h \cdot \big(u_{j}^{
\mathrm{exact}} - (u_{h})_{j} \big)^{2}}, \\
L_{\infty}(h) = \max\limits_{j} \big\lvert u_{j}^{\mathrm{exact}}-
(u_{h})_{j} \big\rvert,
\end{array}
\label{normsDefinition}
\end{equation*}
where $h = \Delta x$ is the spatial step size.

We focus on the performance of the fifth-order WENO-MAIM$1$($10, 
1\mathrm{e-}6, 0.06$), WENO-MAIM$2$($2, 0.1, 10, 1\mathrm{e-}6$), 
WENO-MAIM$3$($10, 1\mathrm{e-}6$) and WENO-MAIM$4$($1.0, 1\mathrm{e-}
6$) schemes. Without causing any confusion, we denote these four 
schemes as WENO-MAIM$1$, WENO-MAIM$2$, WENO-MAIM$3$ and WENO-MAIM$4$ 
for simplicity. Notably, the main purpose for choosing the specific 
WENO-MAIM$2$ scheme here is to verify the conclusion in subsection 
\ref{ParametricStudy} that WENO-MAIM$2$($k, A, Q, 0$) = 
WENO-IM($k, A$), through numerical examples. We show that for all 
the following examples, the WENO-MAIM$2$ scheme provides exactly the 
same results as those of the WENO-IM($2, 0.1$) scheme. 

\begin{table}[!th]
\begin{scriptsize}
\centering
\caption{Convergence orders for the WENO-JS, WENO-M(1), WENO-M(2), 
WENO-IM($k, A$) and WENO-MAIM$i$ schemes.}
\label{table_convergenceOrder}
\begin{tabular*}{\hsize}
{@{}@{\extracolsep{\fill}}lllllll@{}}
\hline
\space  & $n_{\mathrm{cp}}$ & $r_{\mathrm{c}}$-WENO-JS & $r_{\mathrm{
c}}$-WENO-M($1$) & $r_{\mathrm{c}}$-WENO-M($2$) & $r_{\mathrm{c}}$
-WENO-IM($k, A$) & $r_{\mathrm{c}}$-WENO-MAIM$i$ \\
\hline
$r = 2$  &  $0$  &  $3$   &  $3$   &  $3$   &  $3(k \geq 2)$   &  $3(
k \geq 1)$ \\
\space   &  $1$  &  $1$   &  $1$   &  $1$   &  $1$             &  $1$
\\
$r = 3$  &  $0$  &  $5$   &  $5$   &  $5$   &  $5(k \geq 2)$   &  $5(
k \geq 1)$ \\
\space   &  $1$  &  $3$   &  $5$   &  $5$   &  $5(k \geq 2)$   &  $5(
k \geq 1)$ \\
\space   &  $2$  &  $2$   &  $2$   &  $2$   &  $2$             &  $2$
\\
$r = 4$  &  $0$  &  $7$   &  $7$   &  $7$   &  $7(k \geq 2)$   &  $7(
k \geq 1)$ \\
\space   &  $1$  &  $5$   &  $7$   &  $7$   &  $7(k \geq 2)$   &  $7(
k \geq 1)$ \\
\space   &  $2$  &  $4$   &  $6$   &  $7$   &  $7(k \geq 4)$   &  $7(
k \geq 3)$ \\
\space   &  $3$  &  $3$   &  $3$   &  $3$   &  $3$             &  $3$
\\
$r = 5$  &  $0$  &  $9$   &  $9$   &  $9$   &  $9(k \geq 2)$   &  $9(
k \geq 1)$ \\
\space   &  $1$  &  $7$   &  $9$   &  $9$   &  $9(k \geq 2)$   &  $9(
k \geq 1)$ \\
\space   &  $2$  &  $6$   &  $9$   &  $9$   &  $9(k \geq 2)$   &  $9(
k \geq 1)$ \\
\space   &  $3$  &  $5$   &  $7$   &  $9$   &  $9(k \geq 4)$   &  $9(
k \geq 3)$ \\
\space   &  $4$  &  $4$   &  $4$   &  $4$   &  $4$             &  $4$
\\
$r = 6$  &  $0$  &  $11$  &  $11$  &  $11$  &  $11(k \geq 2)$  &  $11
(k \geq 1)$ \\
\space   &  $1$  &  $9$   &  $11$  &  $11$  &  $11(k \geq 2)$  &  $11
(k \geq 1)$ \\
\space   &  $2$  &  $8$   &  $11$  &  $11$  &  $11(k \geq 2)$  &  $11
(k \geq 1)$ \\
\space   &  $3$  &  $7$   &  $11$  &  $11$  &  $11(k \geq 2)$  &  $11
(k \geq 1)$ \\
\space   &  $4$  &  $6$   &  $8$   &  $11$  &  $11(k \geq 6)$  &  $11
(k \geq 5)$ \\
\space   &  $5$  &  $5$   &  $5$   &  $5$   &  $5$             &  $5$
\\
$r = 7$  &  $0$  &  $13$  &  $13$  &  $13$  &  $13(k \geq 2)$  &  $13
(k \geq 1)$ \\
\space   &  $1$  &  $11$  &  $13$  &  $13$  &  $13(k \geq 2)$  &  $13
(k \geq 1)$ \\
\space   &  $2$  &  $10$  &  $13$  &  $13$  &  $13(k \geq 2)$  &  $13
(k \geq 1)$ \\
\space   &  $3$  &  $9$   &  $13$  &  $13$  &  $13(k \geq 2)$  &  $13
(k \geq 1)$ \\
\space   &  $4$  &  $8$   &  $12$  &  $13$  &  $13(k \geq 4)$  &  $13
(k \geq 3)$ \\
\space   &  $5$  &  $7$   &  $9$   &  $13$  &  $13(k \geq 6)$  &  $13
(k \geq 5)$ \\
\space   &  $6$  &  $6$   &  $6$   &  $6$   &  $6$             &  $6$
\\
$r = 8$  &  $0$  &  $15$  &  $15$  &  $15$  &  $15(k \geq 2)$  &  $15
(k \geq 1)$ \\
\space   &  $1$  &  $13$  &  $15$  &  $15$  &  $15(k \geq 2)$  &  $15
(k \geq 1)$ \\
\space   &  $2$  &  $12$  &  $15$  &  $15$  &  $15(k \geq 2)$  &  $15
(k \geq 1)$ \\
\space   &  $3$  &  $11$  &  $15$  &  $15$  &  $15(k \geq 2)$  &  $15
(k \geq 1)$ \\
\space   &  $4$  &  $10$  &  $15$  &  $15$  &  $15(k \geq 2)$  &  $15
(k \geq 1)$ \\
\space   &  $5$  &  $9$   &  $13$  &  $15$  &  $15(k \geq 4)$  &  $15
(k \geq 3)$ \\
\space   &  $6$  &  $8$   &  $10$  &  $15$  &  $15(k \geq 8)$  &  $15
(k \geq 7)$ \\
\space   &  $7$  &  $7$   &  $7$   &  $7$   &  $7$             &  $7$
\\
$r = 9$  &  $0$  &  $17$  &  $17$  &  $17$  &  $17(k \geq 2)$  &  $17
(k \geq 1)$ \\
\space   &  $1$  &  $15$  &  $17$  &  $17$  &  $17(k \geq 2)$  &  $17
(k \geq 1)$ \\
\space   &  $2$  &  $14$  &  $17$  &  $17$  &  $17(k \geq 2)$  &  $17
(k \geq 1)$ \\
\space   &  $3$  &  $13$  &  $17$  &  $17$  &  $17(k \geq 2)$  &  $17
(k \geq 1)$ \\
\space   &  $4$  &  $12$  &  $17$  &  $17$  &  $17(k \geq 2)$  &  $17
(k \geq 1)$ \\
\space   &  $5$  &  $11$  &  $17$  &  $17$  &  $17(k \geq 2)$  &  $17
(k \geq 1)$ \\
\space   &  $6$  &  $10$  &  $14$  &  $17$  &  $17(k \geq 4)$  &  $17
(k \geq 3)$ \\
\space   &  $7$  &  $9$   &  $11$  &  $17$  &  $17(k \geq 8)$  &  $17
(k \geq 7)$ \\
\space   &  $8$  &  $8$   &  $8$   &  $8$   &  $8$             &  $8$
\\
\hline
\end{tabular*}
\end{scriptsize}
\end{table}

\subsection{Accuracy test}
\begin{example}
\bf{(Accuracy test without any critical points \cite{WENO-IM})} \rm{
We solve the one-dimensional linear advection equation $u_{t} + u_{x}
= 0$ with the following initial condition} 
\label{LAE1}
\end{example} 
\begin{equation}
u(x, 0) = \sin ( \pi x ) 
\label{eq:LAE:IC1}
\end{equation}
and the periodic boundary conditions. Note that we consider only the 
fifth-order methods here, and to ensure that the error for the 
overall scheme is a measure of the spatial convergence only, we set 
the CFL number to be $(\Delta x)^{2/3}$. Clearly, the initial 
condition in (\ref{eq:LAE:IC1}) has no critical points. 

Table \ref{table_LAE1} shows the errors and convergence rates (in 
brackets, hereinafter the same) for the numerical solutions of 
Example \ref{LAE1} produced by various schemes at output time 
$t = 2.0$. The results of the three rows are $L_{1}$-, $L_{2}$- and 
$L_{\infty}$- norm errors and orders in turn (the same below). All 
the schemes achieve close to their designed order of accuracy. In 
terms of accuracy, the WENO-MAIM$1$, WENO-MAIM$2$, WENO-MAIM$3$ and 
WENO-MAIM$4$ schemes provide almost equally accurate numerical 
solutions as the WENO-M and WENO-IM($2,0.1$) schemes, which provide 
more accurate numerical solutions than the WENO-JS scheme.

\begin{table}[ht]
\begin{scriptsize}
\centering
\caption{Convergence properties of various schemes solving $u_{t} + u
_{x} = 0$ with initial condition $u(x, 0) = \sin (\pi x)$.}
\label{table_LAE1}
\begin{tabular*}{\hsize}
{@{}@{\extracolsep{\fill}}llllll@{}}
\hline
$N$    & $20$        & $40$        & $80$        & $160$ & $320$ \\
\hline
WENO-JS     & 2.96529e-03(-) & 9.27609e-05(4.9985) & 
2.89265e-06(5.0031) & 9.03392e-08(5.0009) & 2.82330e-09(4.9999)  \\
{}          & 2.42673e-03(-) & 7.64322e-05(4.9887) & 
2.33581e-06(5.0322) & 7.19259e-08(5.0213) & 2.23105e-09(5.0107)  \\
{}          & 2.57899e-03(-) & 9.05453e-05(4.8320) & 
2.90709e-06(4.9610) & 8.85753e-08(5.0365) & 2.72458e-09(5.0228)  \\
WENO-M      & 5.18291e-04(-) & 1.59422e-05(5.0288) & 
4.98914e-07(4.9979) & 1.56021e-08(4.9990) & 4.88356e-10(4.9977)  \\
{}          & 4.06148e-04(-) & 1.25236e-05(5.0193) & 
3.91875e-07(4.9981) & 1.22541e-08(4.9991) & 3.83568e-10(4.9976)  \\
{}          & 3.94913e-04(-) & 1.24993e-05(4.9816) & 
3.91808e-07(4.9956) & 1.22538e-08(4.9988) & 3.83541e-10(4.9977)  \\
WENO-IM(2,0.1)  
            & 5.04401e-04(-) & 1.59160e-05(4.9860) & 
4.98863e-07(4.9957) & 1.56020e-08(4.9988) & 4.88355e-10(4.9977)  \\
{}          & 3.96236e-04(-) & 1.25033e-05(4.9860) & 
4.98863e-07(4.6475) & 1.22541e-08(5.3473) & 3.83568e-10(4.9976)  \\
{}          & 3.94458e-04(-) & 1.24963e-05(4.9803) & 
3.91797e-07(4.9953) & 1.22538e-08(4.9988) & 3.83547e-10(4.9977)  \\
WENO-MAIM1  & 5.08205e-04(-) & 1.59130e-05(4.9971) & 
4.98858e-07(4.9954) & 1.56020e-08(4.9988) & 4.88355e-10(4.9977)  \\
{}          & 4.26155e-04(-) & 1.25010e-05(5.0913) & 
3.91831e-07(4.9957) & 1.22541e-08(4.9989) & 3.83568e-10(4.9976)  \\
{}          & 5.03701e-04(-) & 1.24960e-05(5.3330) & 
3.91795e-07(4.9952) & 1.22538e-08(4.9988) & 3.83543e-10(4.9977)  \\
WENO-MAIM2      & 5.04401e-04(-) & 1.59160e-05(4.9860) & 
4.98863e-07(4.9957) & 1.56020e-08(4.9988) & 4.88355e-10(4.9977)  \\
{}          & 3.96236e-04(-) & 1.25033e-05(4.9860) & 
4.98863e-07(4.6475) & 1.22541e-08(5.3473) & 3.83568e-10(4.9976)  \\
{}          & 3.94458e-04(-) & 1.24963e-05(4.9803) & 
3.91797e-07(4.9953) & 1.22538e-08(4.9988) & 3.83547e-10(4.9977)  \\
WENO-MAIM3  & 5.02844e-04(-) & 1.59130e-05(4.9818) & 
4.98858e-07(4.9954) & 1.56020e-08(4.9988) & 4.88355e-10(4.9977)  \\
{}          & 3.95138e-04(-) & 1.25010e-05(4.9822) & 
3.91831e-07(4.9957) & 1.22541e-08(4.9989) & 3.83568e-10(4.9976)  \\
{}          & 3.94406e-04(-) & 1.24960e-05(4.9801) & 
3.91795e-07(4.9952) & 1.22538e-08(4.9988) & 3.83543e-10(4.9977) \\
WENO-MAIM4  & 5.02845e-04(-) & 1.59131e-05(4.9818) & 
4.98858e-07(4.9954) & 1.56020e-08(4.9988) & 4.88355e-10(4.9977)  \\
{}          & 3.95139e-04(-) & 1.25010e-05(4.9822) & 
3.91831e-07(4.9957) & 1.22541e-08(4.9989) & 3.83568e-10(4.9976)  \\
{}          & 3.94406e-04(-) & 1.24960e-05(4.9801) & 
3.91795e-07{4.9952} & 1.22538e-08(4.9988) & 3.83540e-10(4.9977)  \\
\hline
\end{tabular*}
\end{scriptsize}
\end{table}

\begin{example}
\bf{(Accuracy test with first-order critical points \cite{WENO-M})} 
\rm{We solve the one-dimensional linear advection equation $u_{t} + u
_{x} = 0$ with the following initial condition}
\label{LAE2}
\end{example}
\begin{equation}
u(x, 0) = \sin \bigg( \pi x - \dfrac{\sin(\pi x)}{\pi} \bigg) 
\label{eq:LAE:IC4}
\end{equation}
and periodic boundary conditions. Again, the CFL number is set to be 
$(\Delta x)^{2/3}$. It is easy to verify that the particular initial 
condition in (\ref{eq:LAE:IC4}) has two first-order critical points, 
which both have a non-vanishing third derivative. 

The errors and convergence rates for the numerical solutions of 
Example \ref{LAE2} produced by various schemes at output time 
$t = 2.0$ are shown in Table \ref{table_LAE2}. We can observe that 
the WENO-MAIM$1$, WENO-MAIM$2$, WENO-MAIM$3$ and WENO-MAIM$4$ 
schemes maintain the fifth-order convergence rate even in the 
presence of critical points. Moreover, in terms of accuracy, these 
four WENO schemes provide the equally accurate results as the WENO-M 
and WENO-IM($2,0.1$) schemes, and are much more accurate than the 
WENO-JS scheme, whose errors are much larger and the convergence 
rates are lower than fifth-order.

\begin{table}[ht]
\begin{scriptsize}
\centering
\caption{Convergence properties of various schemes solving $u_{t} + u
_{x} = 0$ with initial condition $u(x, 0) = \sin ( \pi x - \sin(\pi x
)/\pi )$.}
\label{table_LAE2}
\begin{tabular*}{\hsize}
{@{}@{\extracolsep{\fill}}llllll@{}}
\hline
$N$  & $20$        & $40$        & $80$       & $160$   & $320$ \\
\hline
WENO-JS      & 1.01260e-02(-) & 7.22169e-04(3.8096) & 
3.42286e-05(4.3991) & 1.58510e-06(4.4326) & 7.95517e-08(4.3165)  \\
{}           & 8.72198e-03(-) & 6.76133e-04(3.6893) & 
3.63761e-05(4.2162) & 2.29598e-06(3.9858) & 1.68304e-07(3.7700)  \\
{}           & 1.43499e-02(-) & 1.09663e-03(3.7099) & 
9.02485e-05(3.6030) & 8.24022e-06(3.4531) & 8.31702e-07(3.3085)  \\
WENO-M       & 3.70838e-03(-) & 1.45082e-04(4.6758) & 
4.80253e-06(4.9169) & 1.52120e-07(4.9805) & 4.77083e-09(4.9948)  \\
{}           & 3.36224e-03(-) & 1.39007e-04(4.5962) & 
4.52646e-06(4.9406) & 1.42463e-07(4.9897) & 4.45822e-09(4.9980)  \\
{}           & 5.43666e-03(-) & 2.18799e-04(4.6350) & 
6.81451e-06(5.0049) & 2.14545e-07(4.9893) & 6.71080e-09(4.9987)  \\
WENO-IM(2,0.1)   
             & 4.30725e-03(-) & 1.51327e-04(4.8310) & 
4.85592e-06(4.9618) & 1.52659e-07(4.9914) & 4.77654e-09(4.9982)  \\
{}           & 3.93700e-03(-) & 1.41737e-04(4.7958) & 
4.53602e-06(4.9656) & 1.42479e-07(4.9926) & 4.45805e-09(4.9982)  \\
{}           & 5.84039e-03(-) & 2.10531e-04(4.7940) & 
6.82606e-06(4.9468) & 2.14534e-07(4.9918) & 6.71079e-09(4.9986)  \\
WENO-MAIM1   & 8.07923e-03(-) & 3.32483e-04(4.6029) & 
1.01162e-05(5.0385) & 1.52910e-07(6.0478) & 4.77728e-09(5.0003)  \\
{}           & 7.08117e-03(-) & 3.36264e-04(4.3963) & 
1.49724e-05(4.4892) & 1.42515e-07(6.7150) & 4.45807e-09(4.9986)  \\
{}           & 1.03772e-02(-) & 6.62891e-04(3.9685) & 
4.48554e-05(3.8854) & 2.14522e-07(7.7080) & 6.71079e-09(4.9985)  \\
WENO-MAIM2   & 4.30725e-03(-) & 1.51327e-04(4.8310) & 
4.85592e-06(4.9618) & 1.52659e-07(4.9914) & 4.77654e-09(4.9982)  \\
{}           & 3.93700e-03(-) & 1.41737e-04(4.7958) & 
4.53602e-06(4.9656) & 1.42479e-07(4.9926) & 4.45805e-09(4.9982)  \\
{}           & 5.84039e-03(-) & 2.10531e-04(4.7940) & 
6.82606e-06(4.9468) & 2.14534e-07(4.9918) & 6.71079e-09(4.9986)  \\
WENO-MAIM3   & 4.39527e-03(-) & 1.52219e-04(4.8517) & 
4.86436e-06(4.9678) & 1.52735e-07(4.9931) & 4.77728e-09(4.9987)  \\
{}           & 4.02909e-03(-) & 1.42172e-04(4.8247) & 
4.53770e-06(4.9695) & 1.42486e-07(4.9931) & 4.45807e-09(4.9983)  \\
{}           & 5.89045e-03(-) & 2.09893e-04(4.8107) & 
6.83017e-06(4.9416) & 2.14533e-07(4.9926) & 6.71079e-09(4.9986)  \\
WENO-MAIM4   & 4.94421e-03(-) & 1.52224e-04(5.0215) & 
4.86436e-06(4.9678) & 1.52735e-07(4.9931) & 4.77728e-09(4.9987)  \\
{}           & 4.50651e-03(-) & 1.42174e-04(4.9863) & 
4.53770e-06(4.9696) & 1.42486e-07(4.9931) & 4.45807e-09(4.9983)  \\
{}           & 6.56976e-03(-) & 2.09893e-04(4.9681) & 
6.83018e-06(4.9416) & 2.14533e-07(4.9927) & 6.71079e-09(4.9986)  \\
\hline
\end{tabular*}
\end{scriptsize}
\end{table}

\subsection{Performance on calculating linear advection examples 
with discontinuities at long output times}\label{subsec:keyEx}
After extensive numerical tests, we find that the WENO-MAIM3 and 
WENO-MAIM4 schemes have a significant advantage that they not only 
can obtain high resolution, but also can prevent generating spurious 
oscillations on calculating linear advection problems with 
discontinuities when the output time is large. However, the other 
considered schemes do not have this advantage. To demonstrate this, 
we apply the considered schemes to one-dimensional linear advection 
equation $u_{t} + u_{x} = 0$ with the following two initial 
conditions on uniform grid sizes and at long output times.

Case 1. (SLP) The initial condition is given by
\begin{equation}
\begin{array}{l}
u(x, 0) = \left\{
\begin{array}{ll}
\dfrac{1}{6}\big[ G(x, \beta, z - \hat{\delta}) + 4G(x, \beta, z) + G
(x, \beta, z + \hat{\delta}) \big], & x \in [-0.8, -0.6], \\
1, & x \in [-0.4, -0.2], \\
1 - \big\lvert 10(x - 0.1) \big\rvert, & x \in [0.0, 0.2], \\
\dfrac{1}{6}\big[ F(x, \alpha, a - \hat{\delta}) + 4F(x, \alpha, a) +
 F(x, \alpha, a + \hat{\delta}) \big], & x \in [0.4, 0.6], \\
0, & \mathrm{otherwise},
\end{array}\right. 
\end{array}
\label{eq:SLP}
\end{equation}
where $G(x, \beta, z) = \mathrm{e}^{-\beta (x - z)^{2}}, F(x, \alpha
, a) = \sqrt{\max \big(1 - \alpha ^{2}(x - a)^{2}, 0 \big)}$, and 
the constants $z = -0.7, \hat{\delta} = 0.005, \beta = \dfrac{
\log 2}{36\hat{\delta} ^{2}}, a = 0.5$ and $\alpha = 10$.

Case 2. (BiCWP) The initial condition is given by
\begin{equation}
\begin{array}{l}
u(x, 0) = \left\{
\begin{array}{ll}
0,   & x \in [-1.0, -0.8] \bigcup (-0.2, 0.2] \bigcup (0.8, 1.0], \\
0.5, & x \in (-0.6, -0.4] \bigcup (0.2, 0.4]  \bigcup (0.6, 0.8], \\
1,   & x \in (-0.8, -0.6] \bigcup (-0.4, -0.2] \bigcup (0.4, 0.6].
\end{array}\right. 
\end{array}
\label{eq:LAE:BiCWP}
\end{equation}

The periodic boundary condition is used in the two directions and 
the CFL number is set to be $0.1$ for both Case 1 and Case 2. Case 1 
consists of a Gaussian, a square wave, a sharp triangle and a 
semi-ellipse, and Case 2 consists of several constant states 
separated by sharp discontinuities at $x= \pm 0.8, \pm 0.6, \pm 0.4, 
\pm 0.2$. For simplicity, we call Case 1 SLP as it is a 
\textit{\textbf{L}inear} \textit{\textbf{P}roblem} first proposed by 
\textit{\textbf{S}hu} et al. in \cite{WENO-JS}, and we call Case 2 
BiCWP as the profile of the exact solution for this \textit{\textbf{
P}roblem} looks like the \textit{\textbf{B}reach \textbf{i}n \textbf{
C}ity \textbf{W}all}. 

Firstly, we calculate both SLP and BiCWP by all considered WENO 
schemes using the uniform grid sizes of $N = 200, 400, 800$ with a 
long output time $t=2000$. Table \ref{table_advan:01} has shown the 
errors and convergence rates, and we find that: (1) in terms of 
accuracy, when solving SLP, the WENO-MAIM3 scheme provides the most 
accurate solutions among all considered WENO schemes at $N=200,400,
800$, and when solving BiCWP, the WENO-MAIM3 and WENO-MAIM4 schemes 
provide more accurate solutions than other considered WENO schemes 
at $N=800$; (2) as expected, the WENO-MAIM$i(i=1,2,3,4)$ schemes, as 
well as the WENO-IM(2, 0.1) scheme, show more accurate solutions 
than the WENO-JS and WENO-M schemes on solving both SLP and BiCWP; 
(3) in terms of convergence rates, the $L_{1}$ and $L_{2}$ orders of 
the WENO-MAIM$i(i=1,2,3,4)$ and WENO-IM(2, 0.1) schemes are 
significantly higher than those of the WENO-JS and WENO-M schemes; 
(4) furthermore, on solutions of SLP, the $L_{1}$ and $L_{2}$ 
orders of the WENO-MAIM$i(i=1,2,3,4)$ and WENO-IM(2, 0.1) schemes 
are approximately $1.0$ and $0.5 \sim 0.7$ respectively, and their
$L_{\infty}$ orders are all positive; (5) on solutions of BiCWP, 
the $L_{1}$ orders are also approximately $1.0$ for the WENO-MAIM1, 
WENO-MAIM3 and WENO-MAIM4 schemes while only about $0.8$ for the 
WENO-IM(2, 0.1) and WENO-MAIM2 schemes, and the $L_{2}$ orders are 
about $0.45$ for the WENO-MAIM1, WENO-MAIM3 and WENO-MAIM4 schemes 
while slightly lower with the value of about $0.4$ for the 
WENO-IM(2, 0.1) and WENO-MAIM2 schemes; (6) on solutions of 
BiCWP, the $L_{\infty}$ orders become negative for the WENO-IM(2, 
0.1), WENO-MAIM1 and WENO-MAIM2 schemes while maintain positive for 
the WENO-MAIM3 and WENO-MAIM4 schemes.

\begin{figure}[ht]
\centering
\includegraphics[height=0.32\textwidth]
{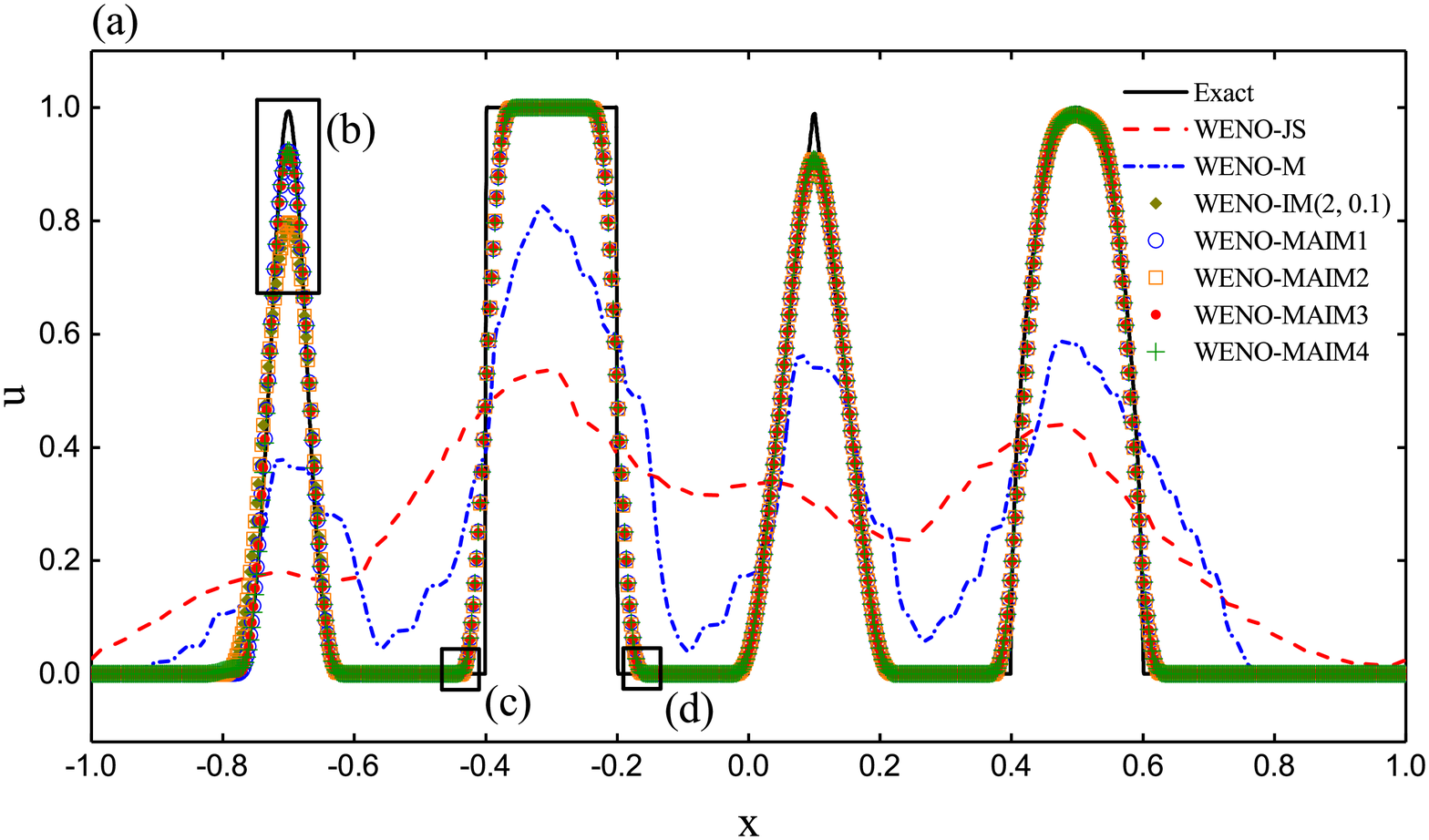}
\includegraphics[height=0.32\textwidth]
{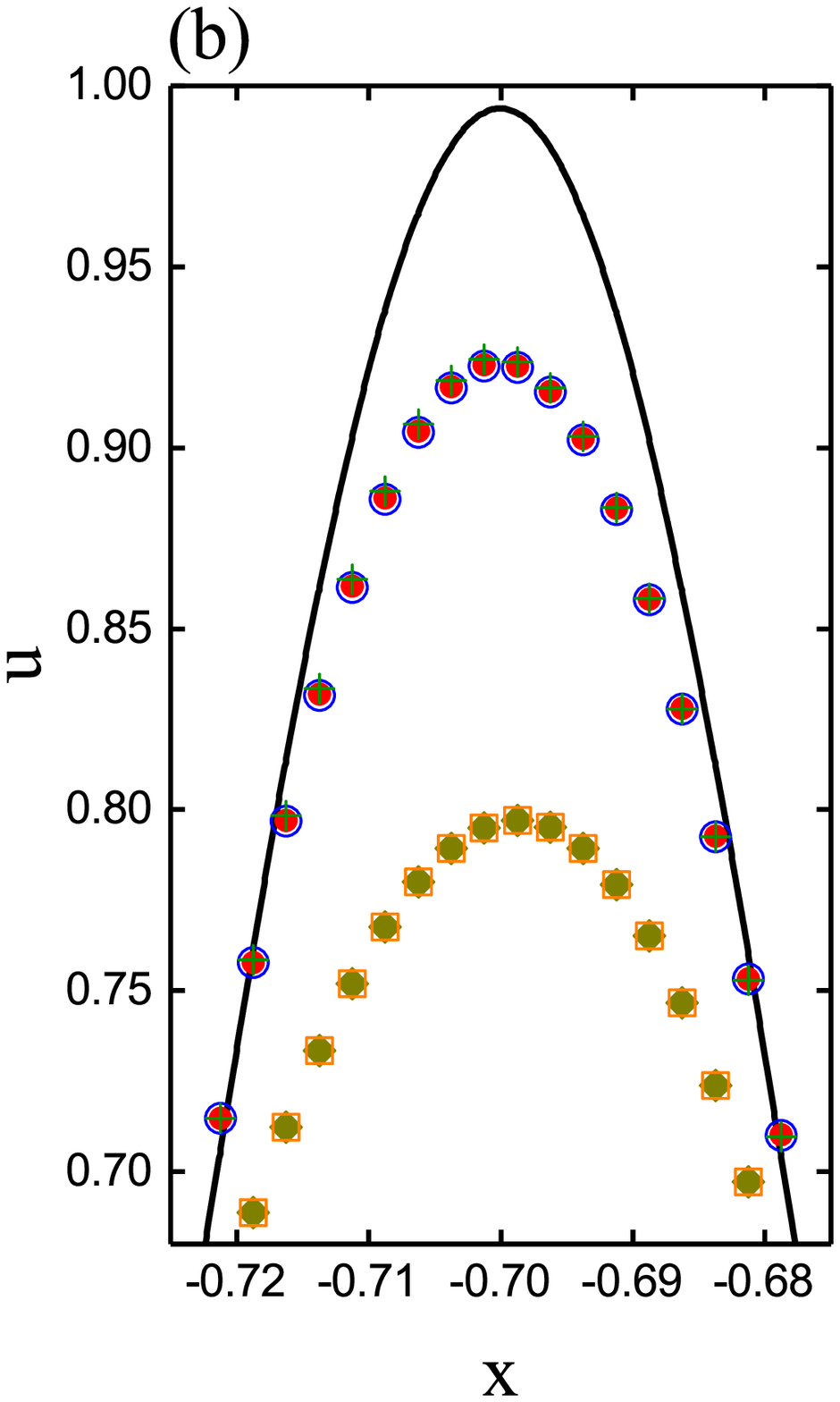}
\includegraphics[height=0.32\textwidth]
{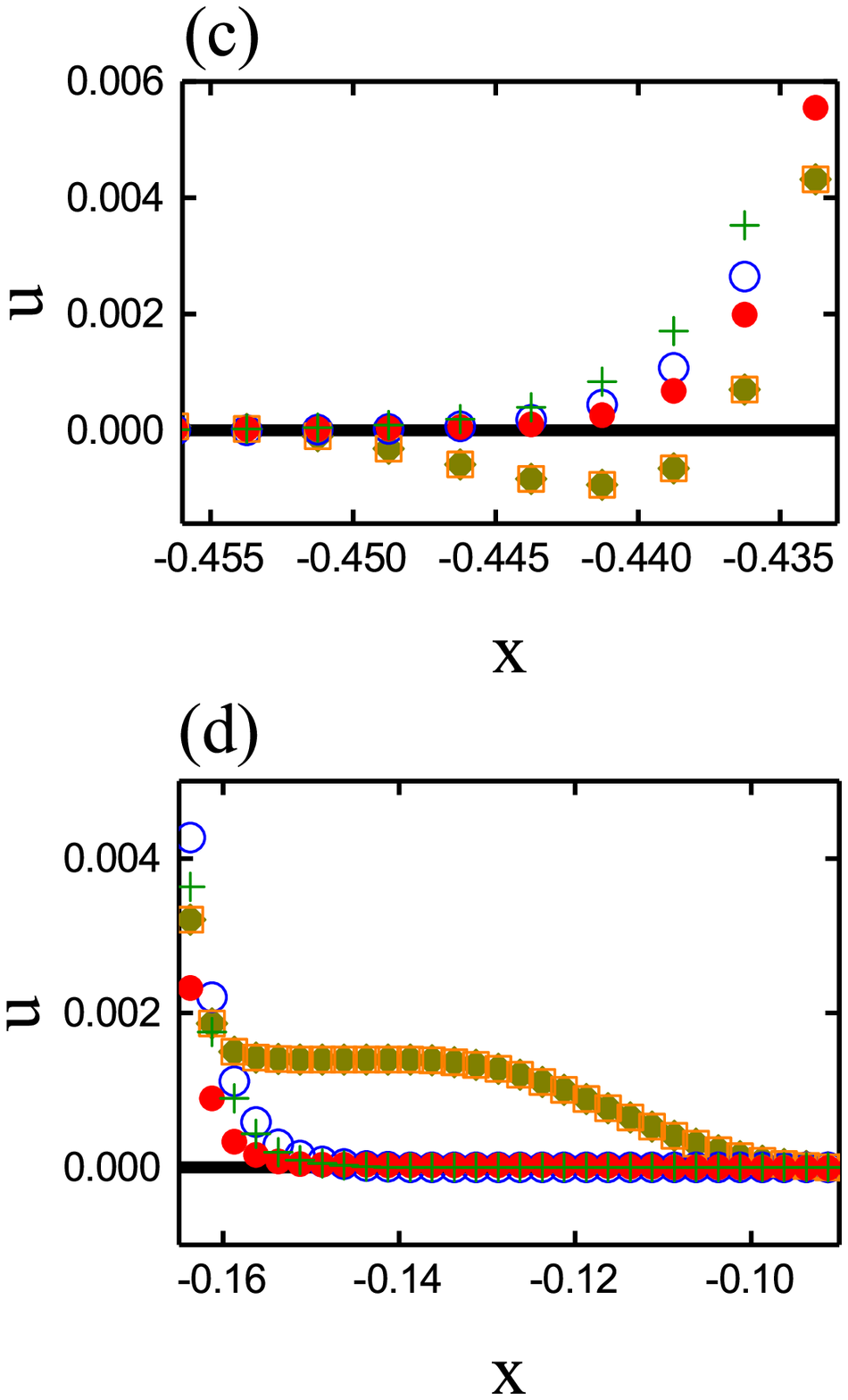}
\caption{Performance of the fifth-order WENO-JS, WENO-M, WENO-IM($2,
0.1$) and WENO-MAIM$i$($i = 1, 2, 3, 4$) schemes for the SLP with 
$N=800$ at long output time $t=2000$.}
\label{fig:SLP:800}
\end{figure}

\begin{figure}[ht]
\centering
\includegraphics[height=0.32\textwidth]
{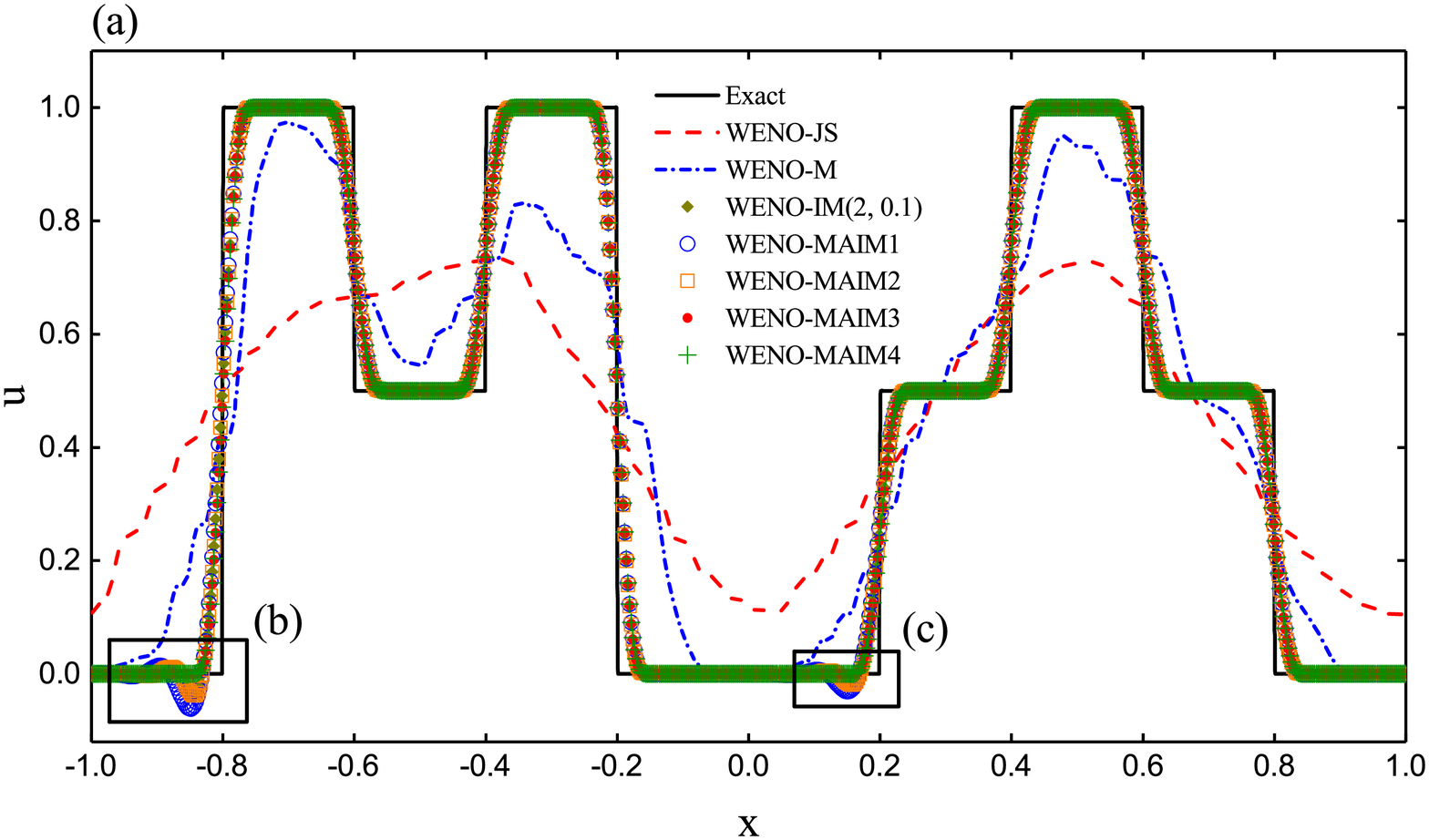}
\includegraphics[height=0.32\textwidth]
{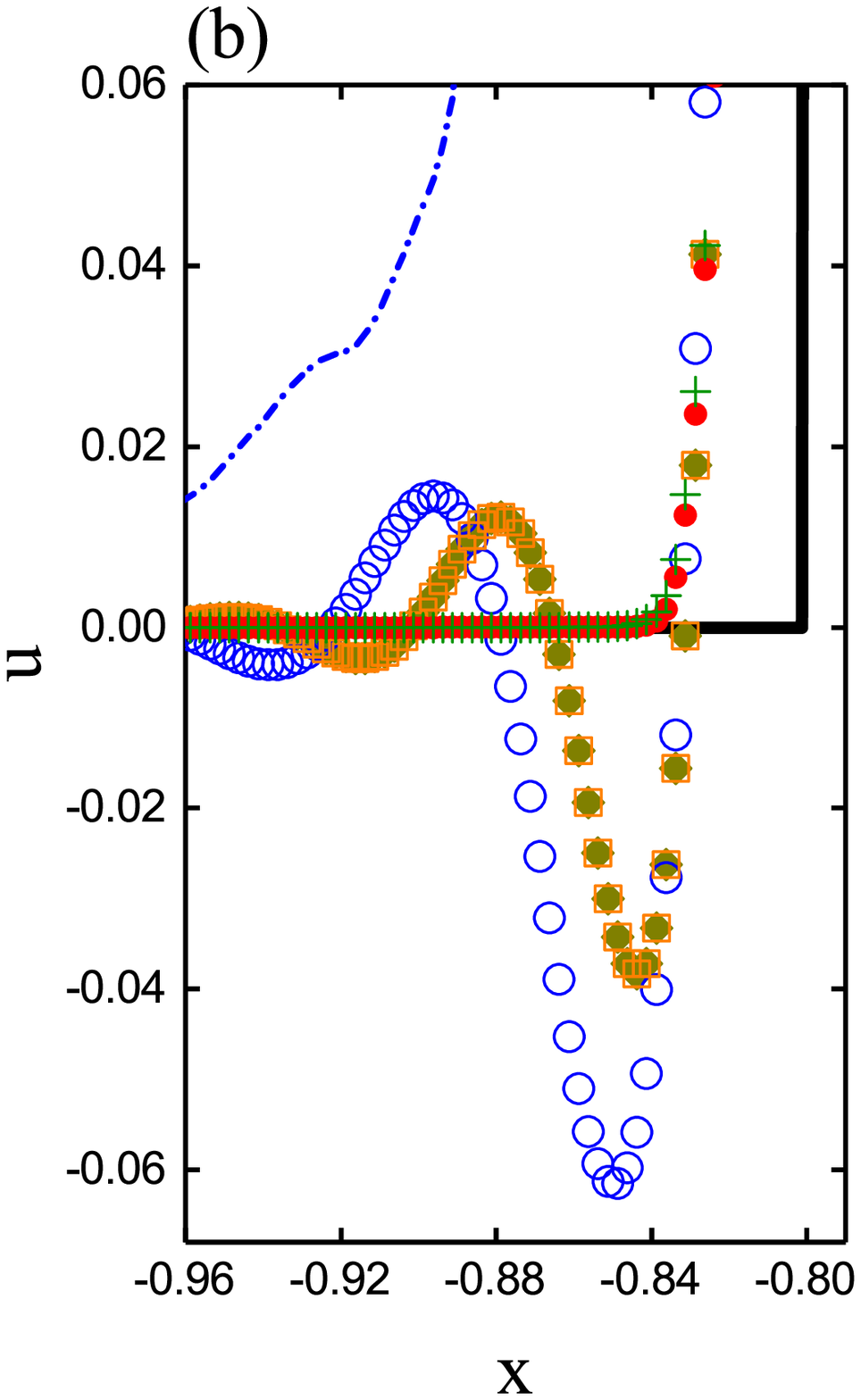}
\includegraphics[height=0.32\textwidth]
{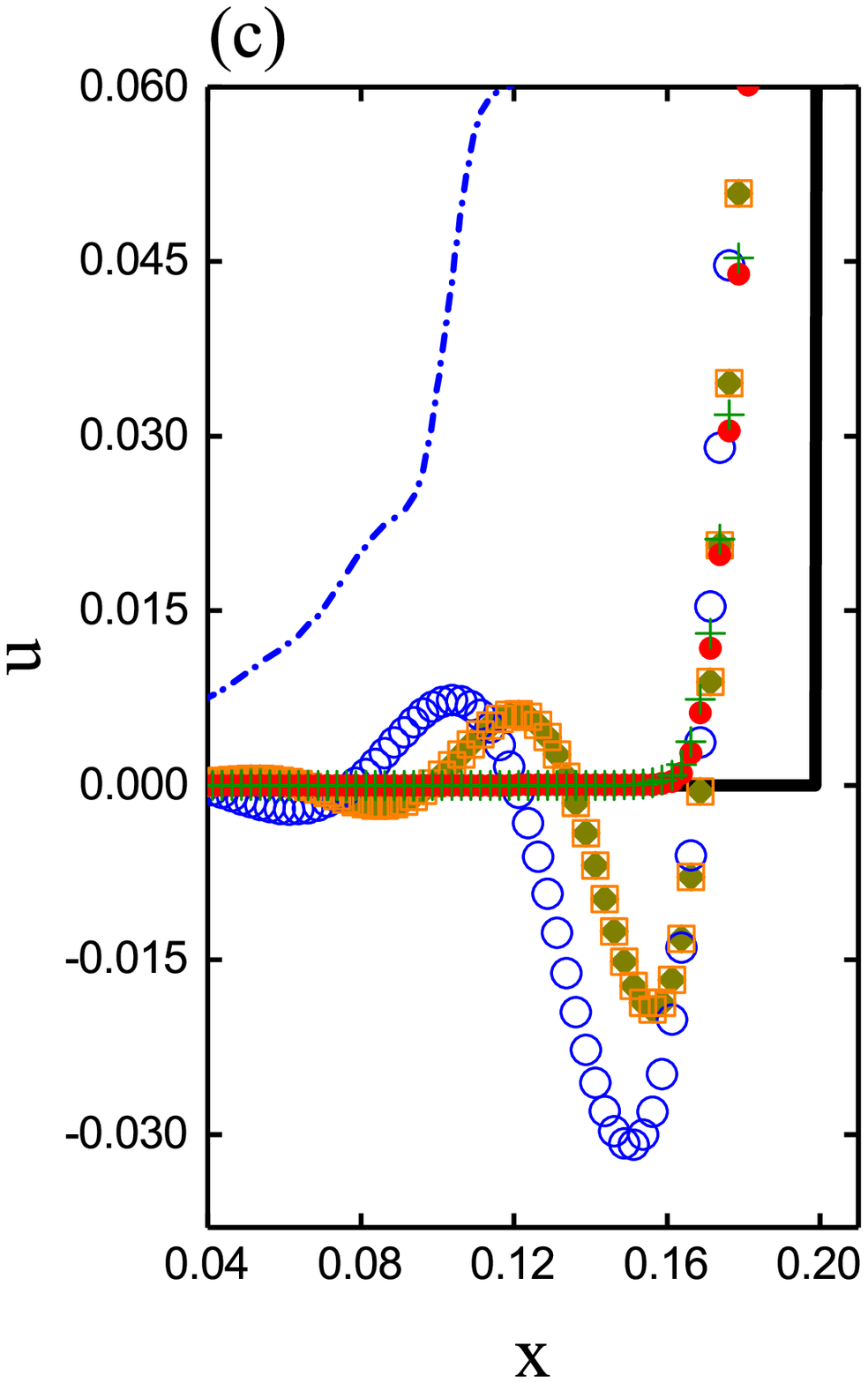}
\caption{Performance of the fifth-order WENO-JS, WENO-M, WENO-IM($2,
0.1$) and WENO-MAIM$i$($i = 1, 2, 3, 4$) schemes for the BiCWP with 
$N=800$ at long output time $t=2000$.}
\label{fig:BiCWP:800}
\end{figure}

Fig. \ref{fig:SLP:800} and Fig. \ref{fig:BiCWP:800} have shown the 
comparisons of various schemes on solving SLP and BiCWP when 
$t = 2000$ and $N=800$. We can observe that: (1) when solving SLP, 
the WENO-MAIM3 and WENO-MAIM4 schemes do not generate spurious 
oscillations and provide the numerical results with a significantly 
higher resolution than those of the WENO-JS and WENO-M schemes; (2) 
in general, the WENO-MAIM1, WENO-MAIM2 and WENO-IM(2, 0.1) schemes 
also show very high resolution and almost do not generate spurious 
oscillations when solving SLP, but if we take a closer look, we can 
see that the WENO-MAIM2 and WENO-IM(2, 0.1) schemes give a 
significantly lower resolution of the Gaussian (see Fig. 
\ref{fig:SLP:800}(b)) and a slightly lower resolution of the square 
wave near $x=-0.16$ (see Fig. \ref{fig:SLP:800}(d)) than the 
WENO-MAIM1, WENO-MAIM3 and WENO-MAIM4 schemes; (3) furthermore, the 
WENO-MAIM2 and WENO-IM(2, 0.1) schemes generate a very slight 
spurious oscillation, which is hard to be noticed without a very 
close look, around the square wave near $x=-0.435$ (see 
Fig. \ref{fig:SLP:800}(c)), while the other considered schemes do 
not generate this spurious oscillation; (4) when solving BiCWP, the 
spurious oscillations generated by the WENO-IM(2, 0.1), WENO-MAIM1 
and WENO-MAIM2 schemes become very evident, while the WENO-MAIM3 
and WENO-MAIM4 schemes still provide the numerical results with 
much sharper resolutions than those of the WENO-JS and WENO-M 
schemes but do not produce spurious oscillations.

\begin{table}[ht]
\begin{scriptsize}
\centering
\caption{Convergence properties of various schemes solving SLP and 
BiCWP at output time $t=2000$.}
\label{table_advan:01}
\begin{tabular*}{\hsize}
{@{}@{\extracolsep{\fill}}lllllll@{}}
\hline
\space  &\multicolumn{3}{l}{SLP}  &\multicolumn{3}{l}{BiCWP}  \\
\cline{2-4}  \cline{5-7}
$N$ & $200$   & $400$  & $800$ & $200$    & $400$     & $800$\\
\hline
WENO-JS    & 6.12899e-01(-) & 5.99215e-01(0.0326) & 5.50158e-01
(0.1232)   
           & 5.89672e-01(-) & 5.56639e-01(0.0832) & 4.72439e-01
(0.2366)    \\
{}         & 5.08726e-01(-) & 5.01160e-01(0.0216) & 4.67585e-01
(0.1000)   
           & 4.70933e-01(-) & 4.41787e-01(0.0922) & 3.78432e-01
(0.2233)    \\
{}         & 7.99265e-01(-) & 8.20493e-01(-0.0378)& 8.14650e-01
(0.0103)   
           & 6.41175e-01(-) & 5.94616e-01(0.1088) & 5.73614e-01
(0.0519)    \\
WENO-M     & 3.81597e-01(-) & 3.25323e-01(0.2302) & 3.48528e-01
(-0.0994)   
           & 3.27647e-01(-) & 2.64334e-01(0.3098) & 2.57390e-01
(0.0384)    \\
{}         & 3.59205e-01(-) & 3.12970e-01(0.1988) & 3.24373e-01
(-0.0516)   
           & 2.73948e-01(-) & 2.60726e-01(0.0714) & 2.47361e-01
(0.0759)    \\
{}         & 6.89414e-01(-) & 6.75473e-01(0.0295) & 6.25645e-01
(0.1106)   
           & 5.12247e-01(-) & 5.60199e-01(-0.1291)& 5.83690e-01
(-0.0593)    \\
WENO-IM(2,0.1)  
           & 2.17411e-01(-) & 1.12590e-01(0.9493) & 5.18367e-02
(1.1190)   
           & 1.96196e-01(-) & 1.12264e-01(0.8054) & 6.48339e-02
(0.7921)    \\
{}         & 2.30000e-01(-) & 1.64458e-01(0.4839) & 9.98968e-02
(0.7192)   
           & 2.07227e-01(-) & 1.54544e-01(0.4232) & 1.16534e-01
(0.4073)    \\
{}         & 5.69864e-01(-) & 4.82180e-01(0.2410) & 4.73102e-01
(0.0274)   
           & 4.98939e-01(-) & 4.68309e-01(0.0914) & 4.91291e-01
(-0.0691)    \\
WENO-MAIM1 & 2.18238e-01(-) & 1.09902e-01(0.9897) & 4.41601e-02
(1.3154)   
           & 2.04996e-01(-) & 1.33104e-01(0.6230) & 6.77607e-02
(0.9740)    \\
{}         & 2.29151e-01(-) & 1.51024e-01(0.6015) & 9.35506e-02
(0.6910)   
           & 2.07725e-01(-) & 1.62892e-01(0.3508) & 1.18338e-01
(0.4610)    \\
{}         & 5.63682e-01(-) & 4.94657e-01(0.1885) & 4.72393e-01
(0.0664)   
           & 4.93792e-01(-) & 4.96724e-01(-0.0085)& 5.13570e-01
(-0.0481)    \\
WENO-MAIM2 & 2.17411e-01(-) & 1.12590e-01(0.9493) & 5.18367e-02
(1.1190)   
           & 1.96196e-01(-) & 1.12264e-01(0.8054) & 6.48339e-02
(0.7921)    \\
{}         & 2.30000e-01(-) & 1.64458e-01(0.4839) & 9.98968e-02
(0.7192)   
           & 2.07227e-01(-) & 1.54544e-01(0.4232) & 1.16534e-01
(0.4073)    \\
{}         & 5.69864e-01(-) & 4.82180e-01(0.2410) & 4.73102e-01
(0.0274)   
           & 4.98939e-01(-) & 4.68309e-01(0.0914) & 4.91291e-01
(-0.0691)    \\
WENO-MAIM3 & 2.17339e-01(-) & 9.91687e-02(1.1320) & 4.37214e-02
(1.1815)   
           & 1.78226e-01(-) & 1.19271e-01(0.5795) & 6.25800e-02
(0.9305)    \\
{}         & 2.28723e-01(-) & 1.45461e-01(0.6530) & 9.31160e-02
(0.6435)   
           & 1.97298e-01(-) & 1.57392e-01(0.3260) & 1.15580e-01
(0.4455)    \\
{}         & 5.63600e-01(-) & 4.79774e-01(0.2323) & 4.71185e-01
(0.0261)   
           & 5.01513e-01(-) & 4.75377e-01(0.0772) & 4.71183e-01
(0.0128)    \\
WENO-MAIM4 & 2.18548e-01(-) & 1.03499e-01(1.0783) & 4.39609e-02
(1.2353)   
           & 2.05283e-01(-) & 1.28349e-01(0.6775) & 6.33284e-02
(1.0191)    \\
{}         & 2.30043e-01(-) & 1.48300e-01(0.6334) & 9.34567e-02
(0.6661)   
           & 2.07890e-01(-) & 1.60928e-01(0.3694) & 1.16112e-01
(0.4709)    \\
{}         & 5.65659e-01(-) & 4.91295e-01(0.2033) & 4.71240e-01
(0.0601)   
           & 4.90417e-01(-) & 4.89146e-01(0.0037) & 4.71240e-01
(0.0538)    \\
\hline
\end{tabular*}
\end{scriptsize}
\end{table}

In addition, we calculate both SLP and BiCWP with a shortened but 
still long output time $t=200$ and larger uniform grid sizes of 
$N=1600, 3200, 6400$. Table \ref{table_advan:02} has shown the 
errors and convergence rates, and we find that: (1) in terms of 
accuracy, the WENO-JS and WENO-M schemes generate less accurate 
solutions than other considered WENO schemes as expected; (2) the 
WENO-MAIM3 scheme provides the most accurate solutions among all 
considered WENO schemes at $N = 1600,3200,6400$ on solving SLP, and 
it also provides the most accurate solutions at $N = 3200,6400$ on 
solving BiCWP; (3) when solving both SLP and BiCWP, the WENO-MAIM4 
scheme gives more accurate solutions than other consider WENO 
schemes except the WENO-MAIM3 scheme at $N = 6400$; (4) in terms of 
convergence rates, the $L_{1}$ and $L_{2}$ orders of the WENO-MAIM$i
(i=1,2,3,4)$ and WENO-IM(2, 0.1) schemes are higher than those of 
the WENO-JS and WENO-M schemes in general; (5) on solutions of 
SLP, the $L_{1}$ and $L_{2}$ orders are approximately $0.88$ and
$0.42$ respectively for the WENO-MAIM1, WENO-MAIM3 and WENO-MAIM4 
schemes, while slightly lower with the values about $0.83$ and $0.40$
respectively for the WENO-IM(2, 0.1) and WENO-MAIM2 schemes; (6) 
similarly, on solutions of BiCWP, the $L_{1}$ and $L_{2}$ orders 
are approximately $0.83$ and $0.42$ respectively for the WENO-MAIM1, 
WENO-MAIM3 and WENO-MAIM4 schemes, while slightly lower with the 
values about $0.80$ and $0.41$ respectively for the WENO-IM(2, 0.1) 
and WENO-MAIM2 schemes; (7) on solutions of both SLP and BiCWP, the 
$L_{\infty}$ orders of the WENO-MAIM$i(1,2,3,4)$ and WENO-IM(2, 0.1) 
schemes are all negative, but the absolute values of the 
$L_{\infty}$ orders of the WENO-MAIM3 and WENO-MAIM4 schemes are 
smaller than those of the WENO-IM(2, 0.1), WENO-MAIM1 and WENO-MAIM2 
schemes. 

\begin{figure}[ht]
\centering
\includegraphics[height=0.32\textwidth]
{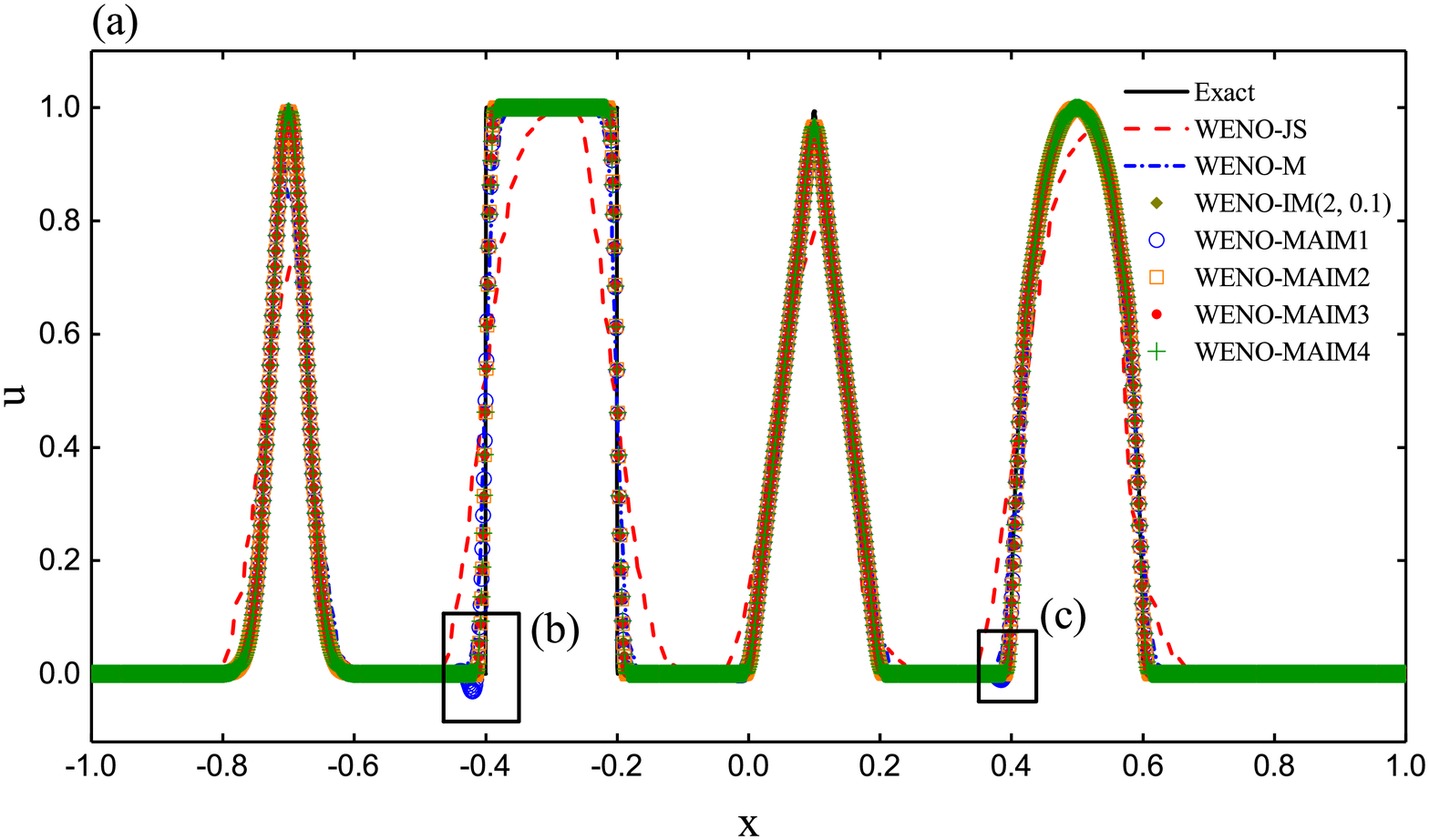}
\includegraphics[height=0.32\textwidth]
{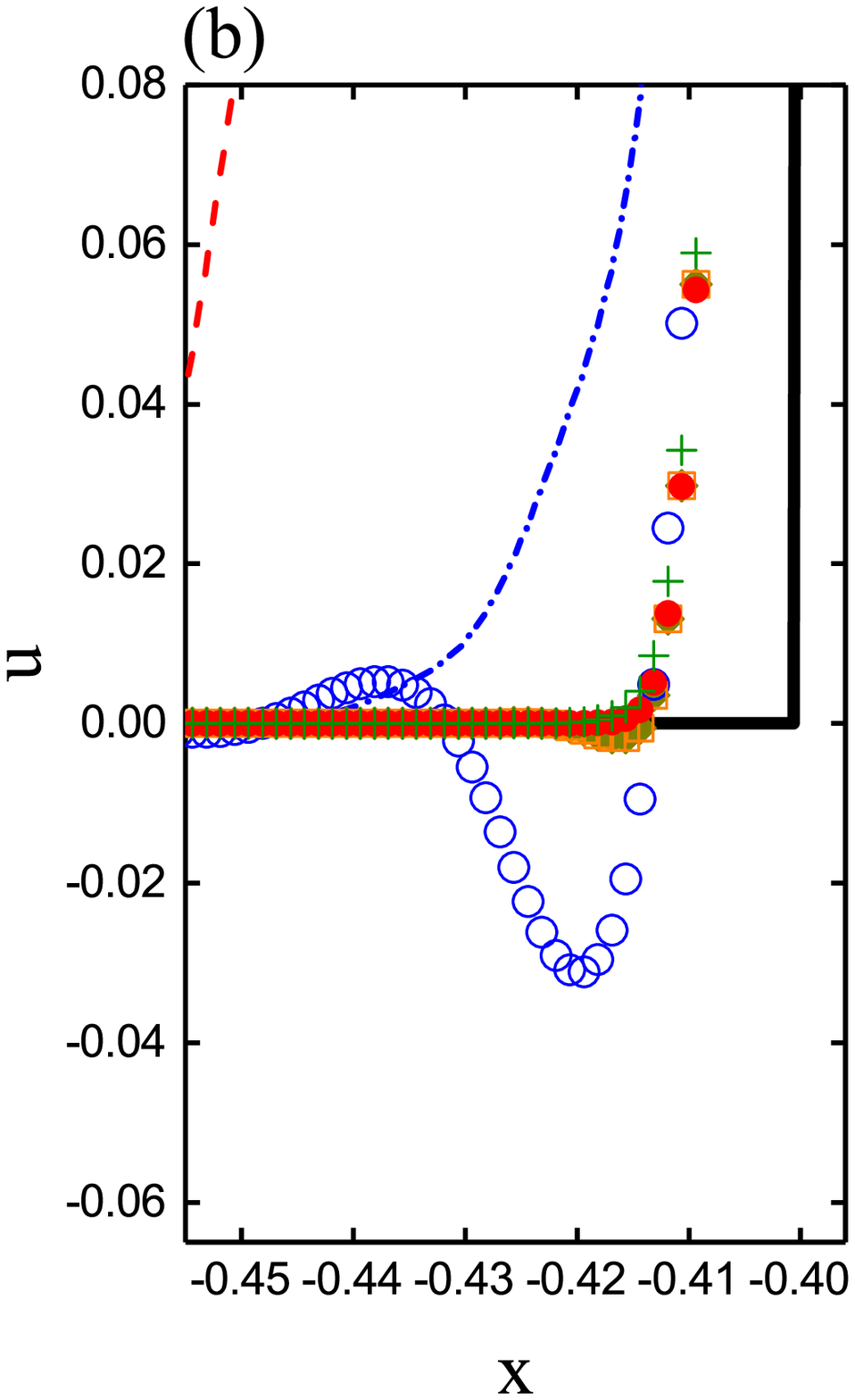}
\includegraphics[height=0.32\textwidth]
{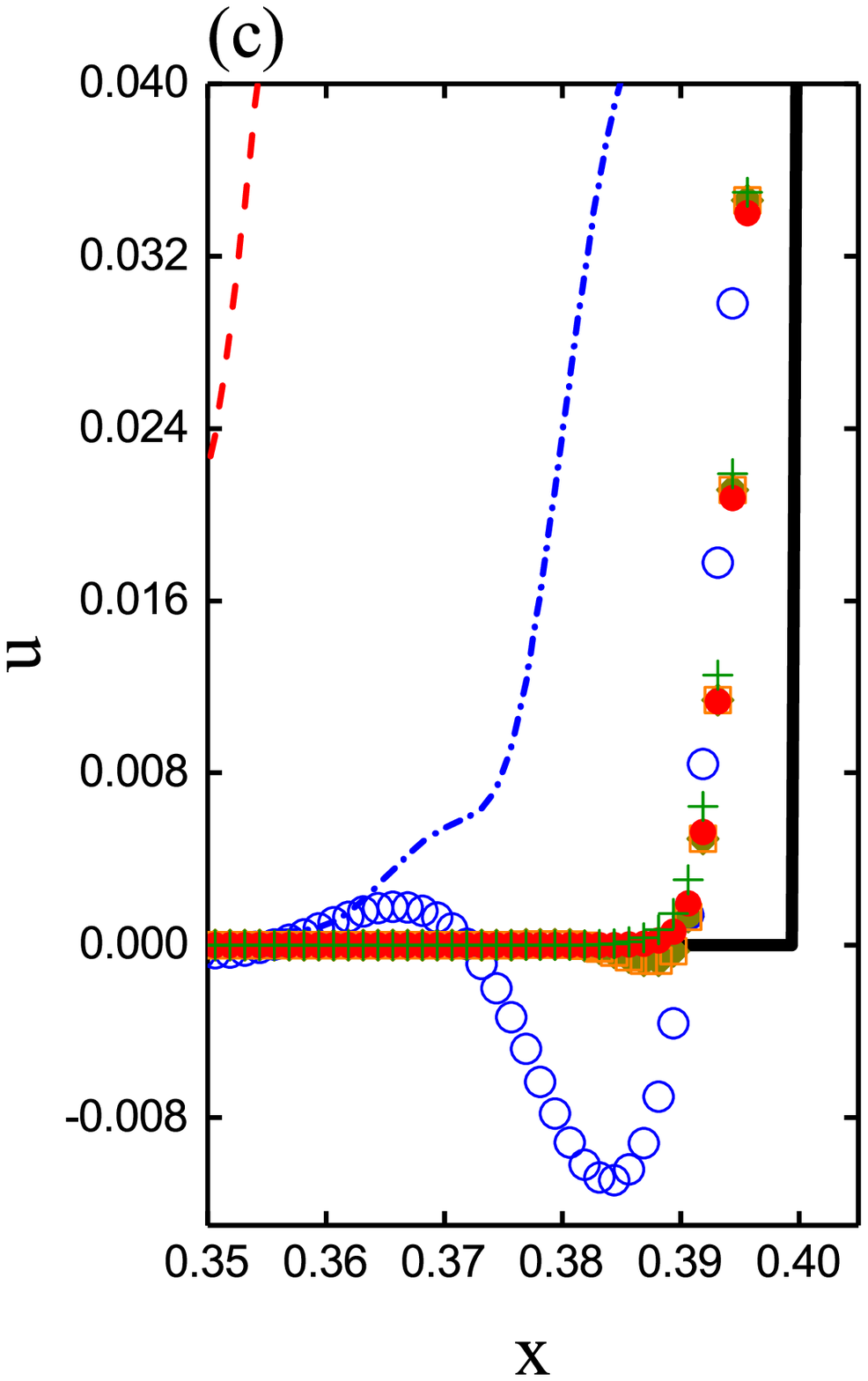}
\caption{Performance of the fifth-order WENO-JS, WENO-M, WENO-IM($2,
0.1$) and WENO-MAIM$i$($i = 1, 2, 3, 4$) schemes for the SLP with 
$N=1600$ at long output time $t=200$.}
\label{fig:SLP:1600}
\end{figure}

\begin{figure}[ht]
\centering
\includegraphics[height=0.32\textwidth]
{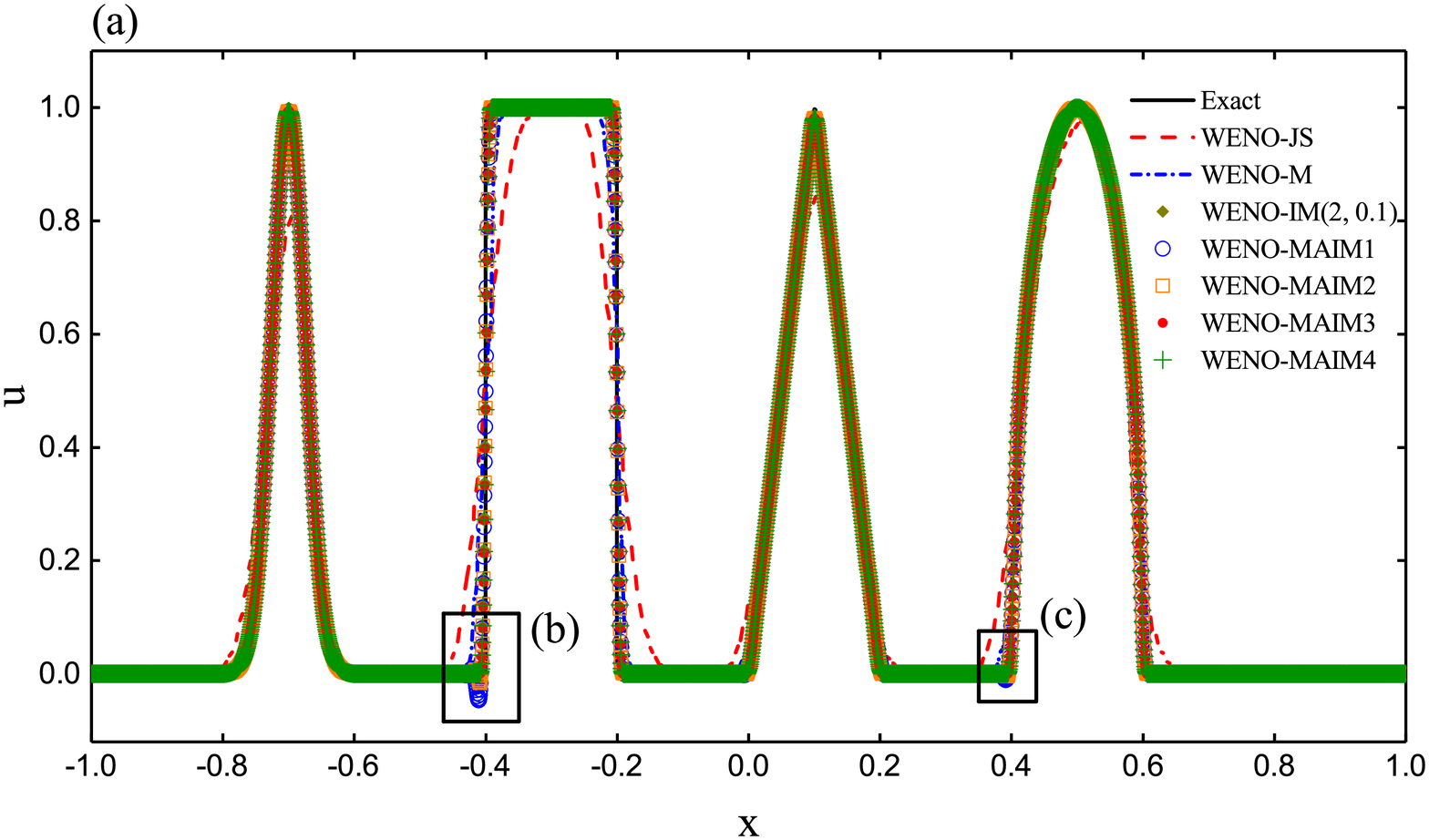}
\includegraphics[height=0.32\textwidth]
{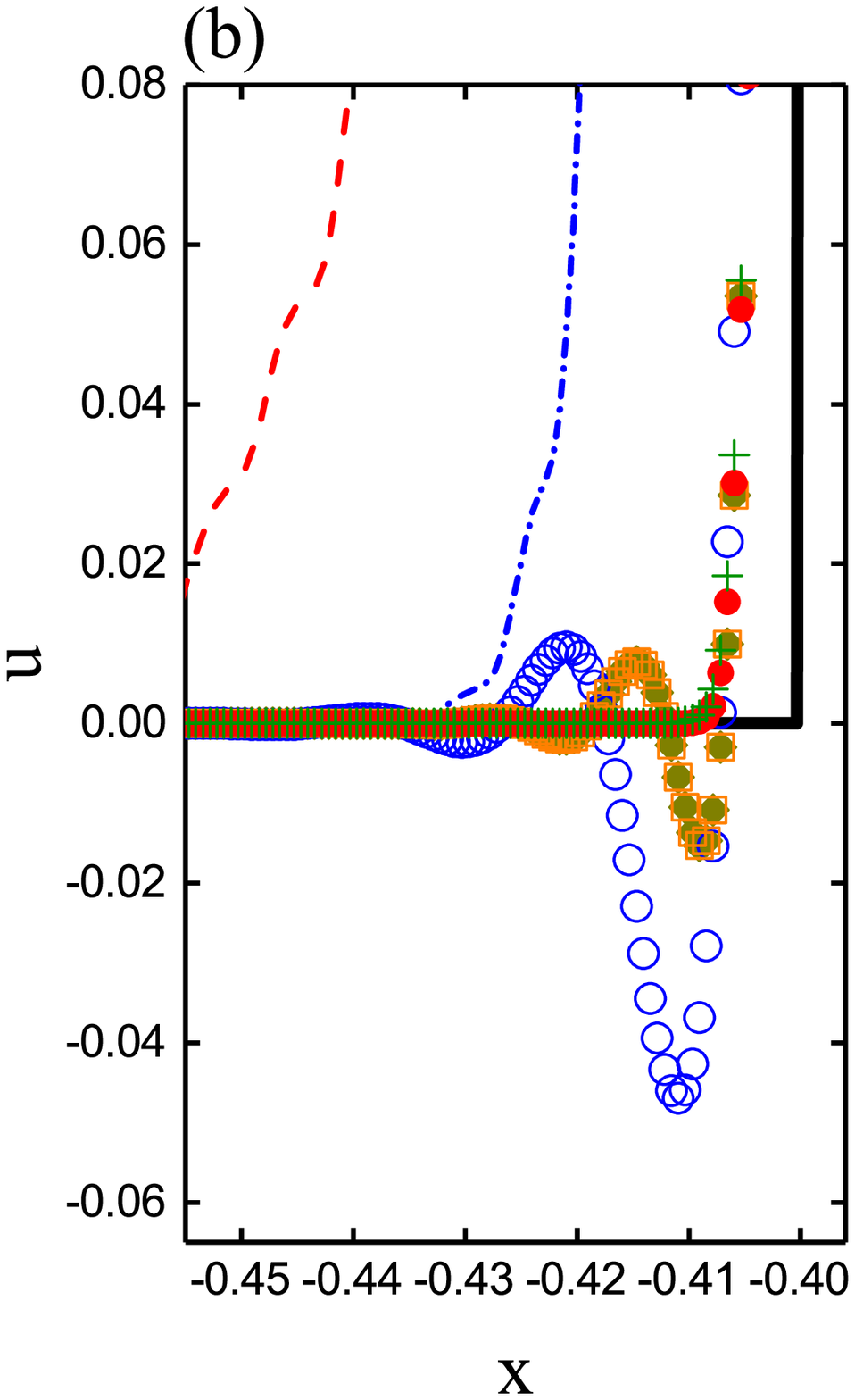}
\includegraphics[height=0.32\textwidth]
{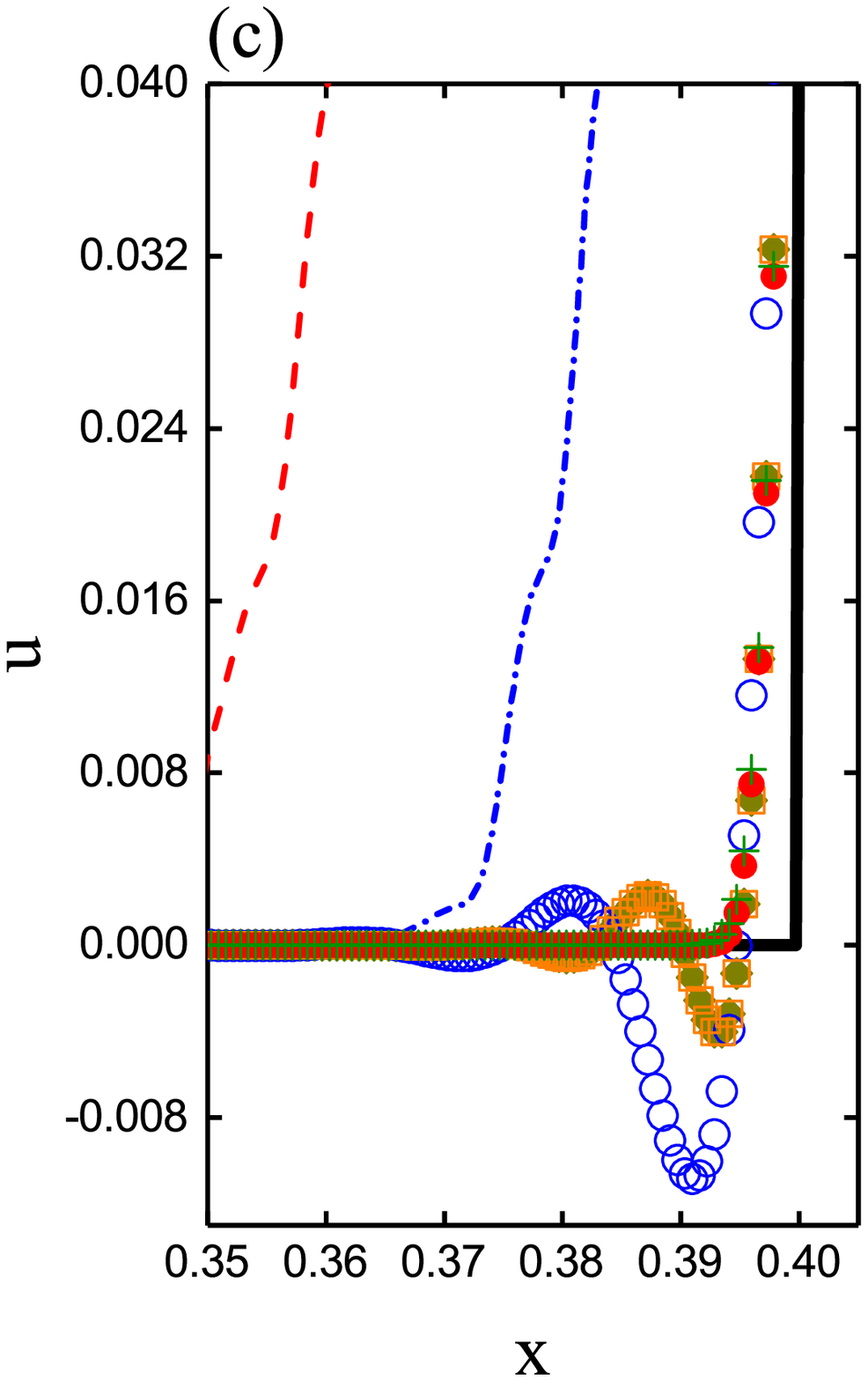}
\caption{Performance of the fifth-order WENO-JS, WENO-M, WENO-IM($2,
0.1$) and WENO-MAIM$i$($i = 1, 2, 3, 4$) schemes for the SLP with 
$N=3200$ at long output time $t=200$.}
\label{fig:SLP:3200}
\end{figure}

\begin{figure}[ht]
\centering
\includegraphics[height=0.32\textwidth]
{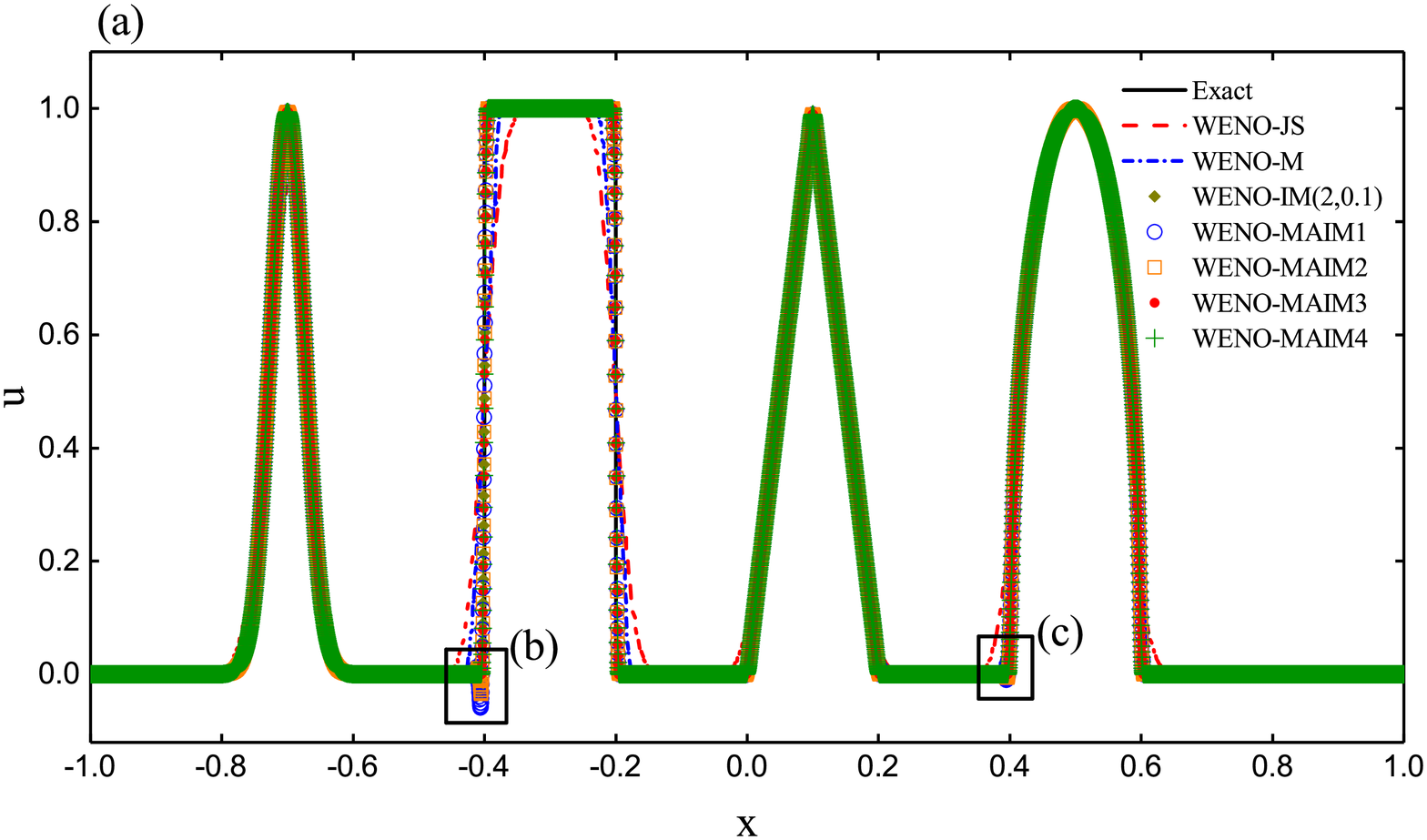}
\includegraphics[height=0.32\textwidth]
{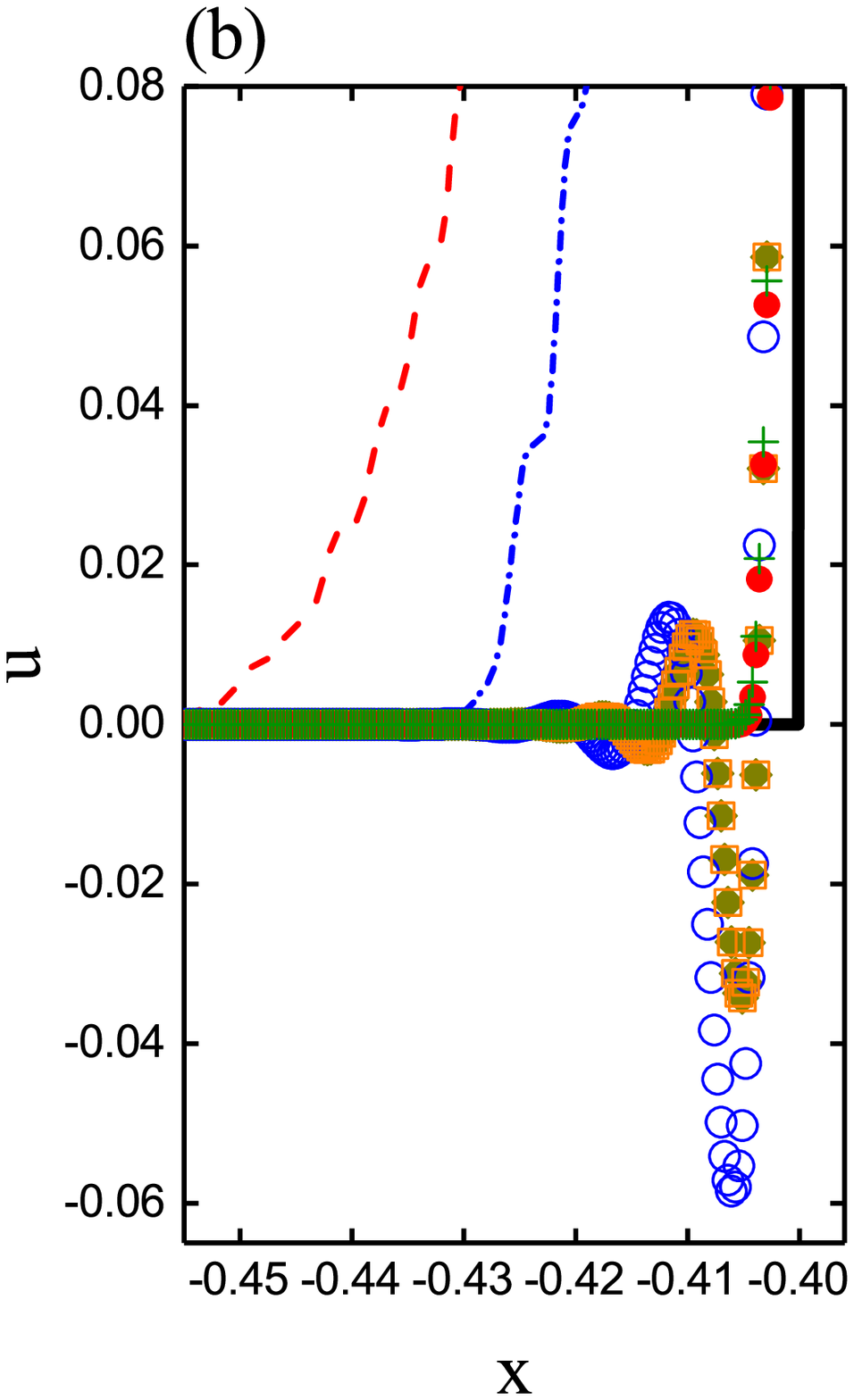}
\includegraphics[height=0.32\textwidth]
{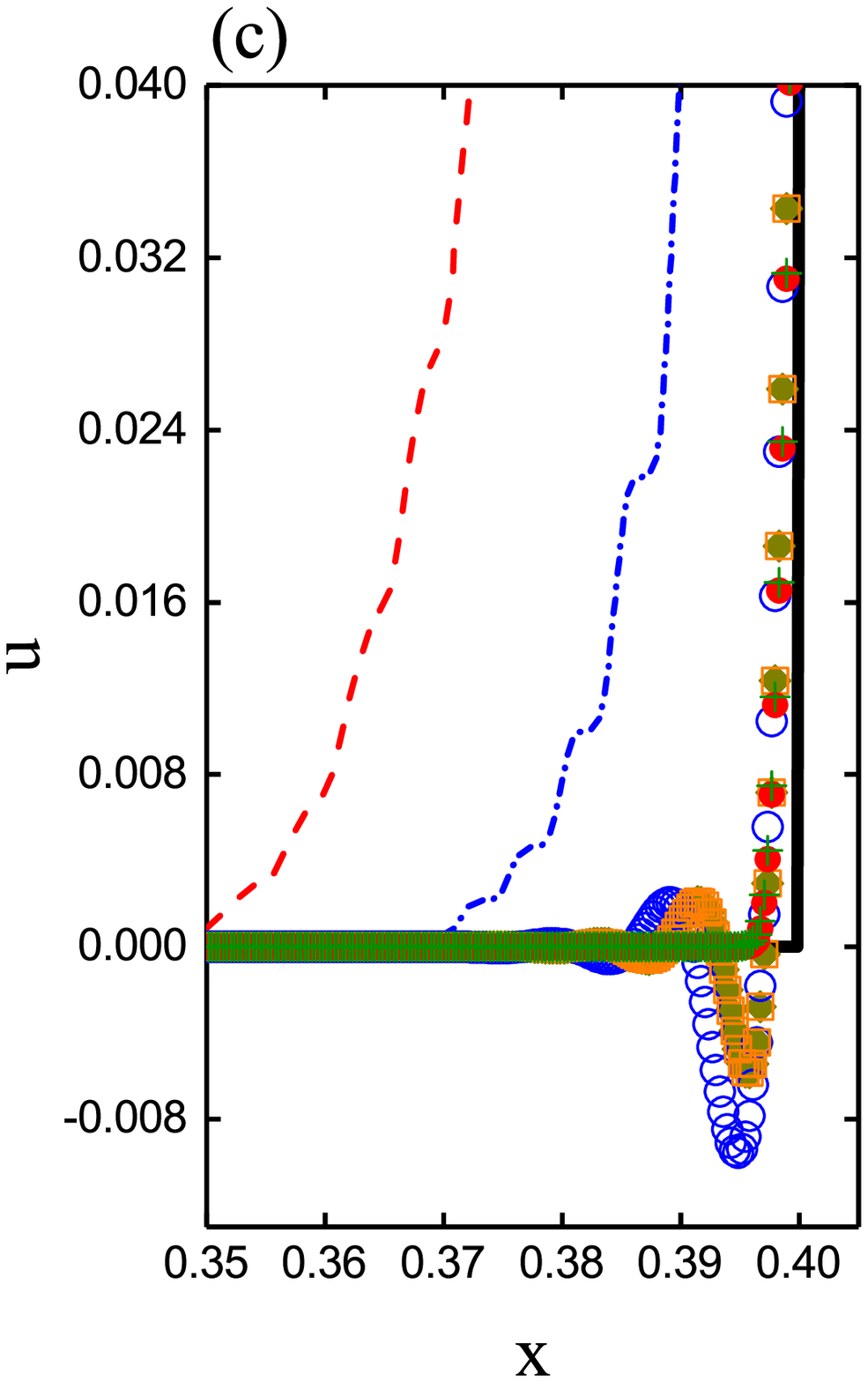}
\caption{Performance of the fifth-order WENO-JS, WENO-M, WENO-IM($2,
0.1$) and WENO-MAIM$i$($i = 1, 2, 3, 4$) schemes for the SLP with 
$N=6400$ at long output time $t=200$.}
\label{fig:SLP:6400}
\end{figure}

\begin{figure}[ht]
\centering
\includegraphics[height=0.32\textwidth]
{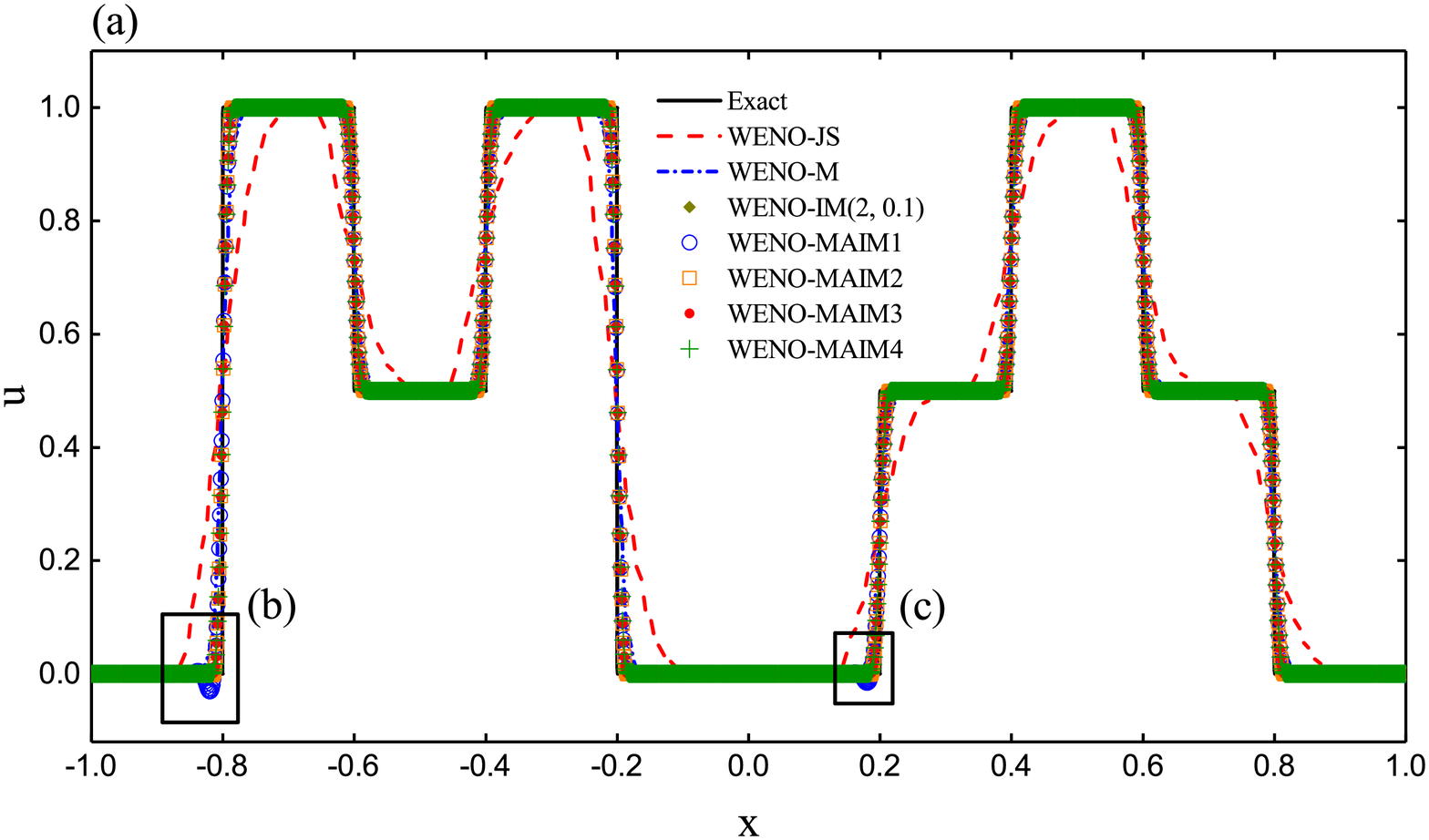}
\includegraphics[height=0.32\textwidth]
{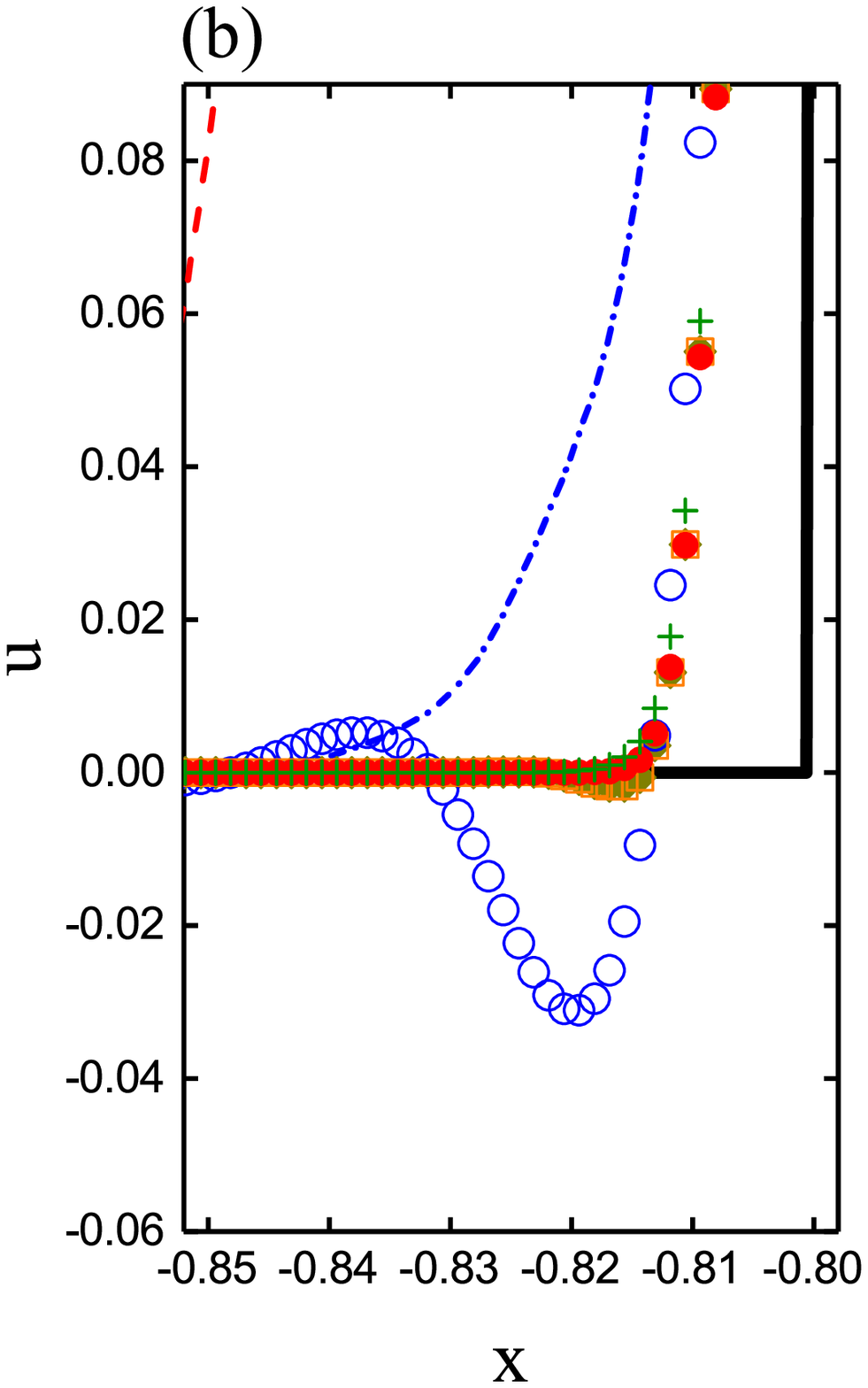}
\includegraphics[height=0.32\textwidth]
{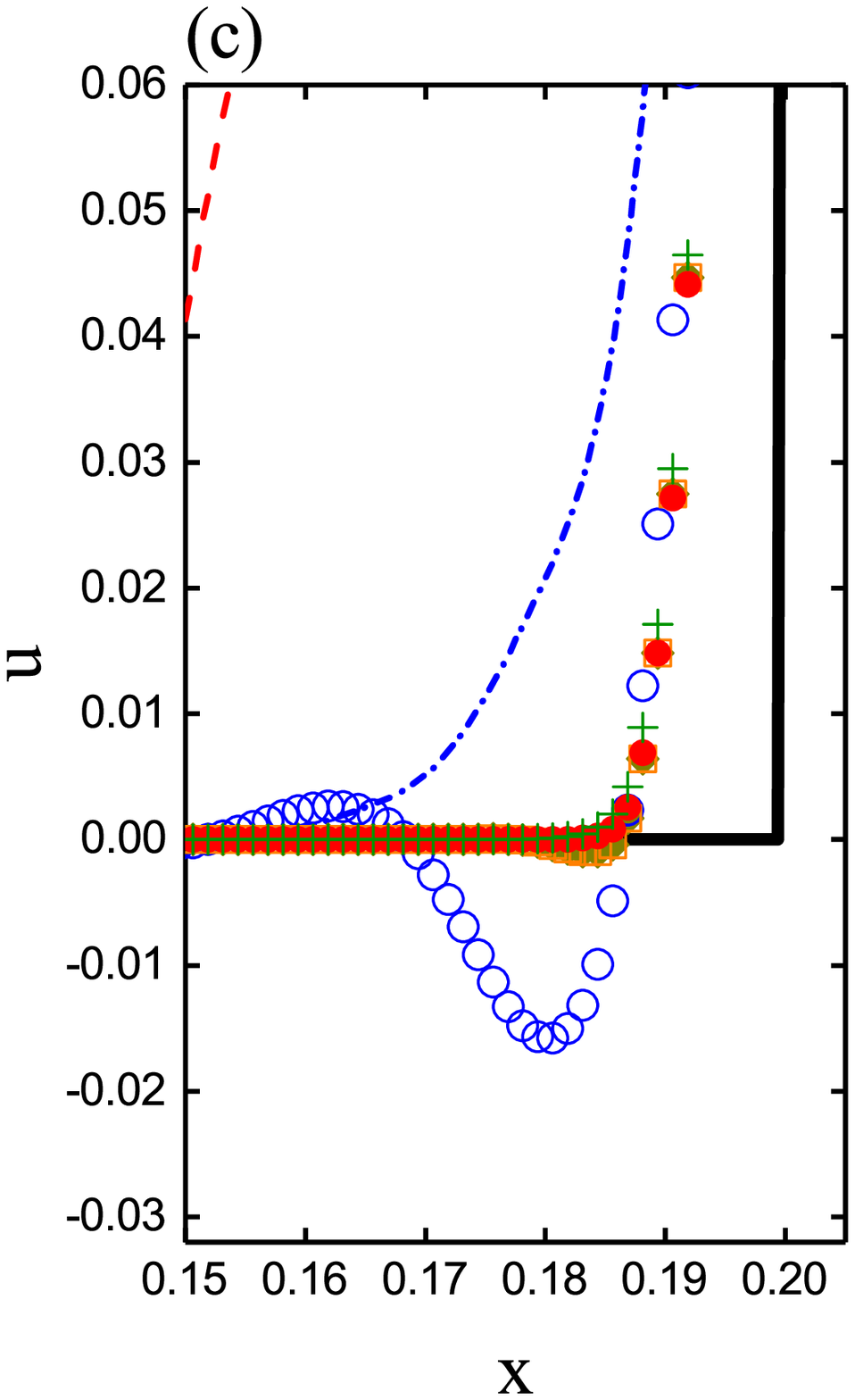}
\caption{Performance of the fifth-order WENO-JS, WENO-M, WENO-IM($2,
0.1$) and WENO-MAIM$i$($i = 1, 2, 3, 4$) schemes for the BiCWP with 
$N=1600$ at long output time $t=200$.}
\label{fig:BiCWP:1600}
\end{figure}

\begin{figure}[ht]
\centering
\includegraphics[height=0.32\textwidth]
{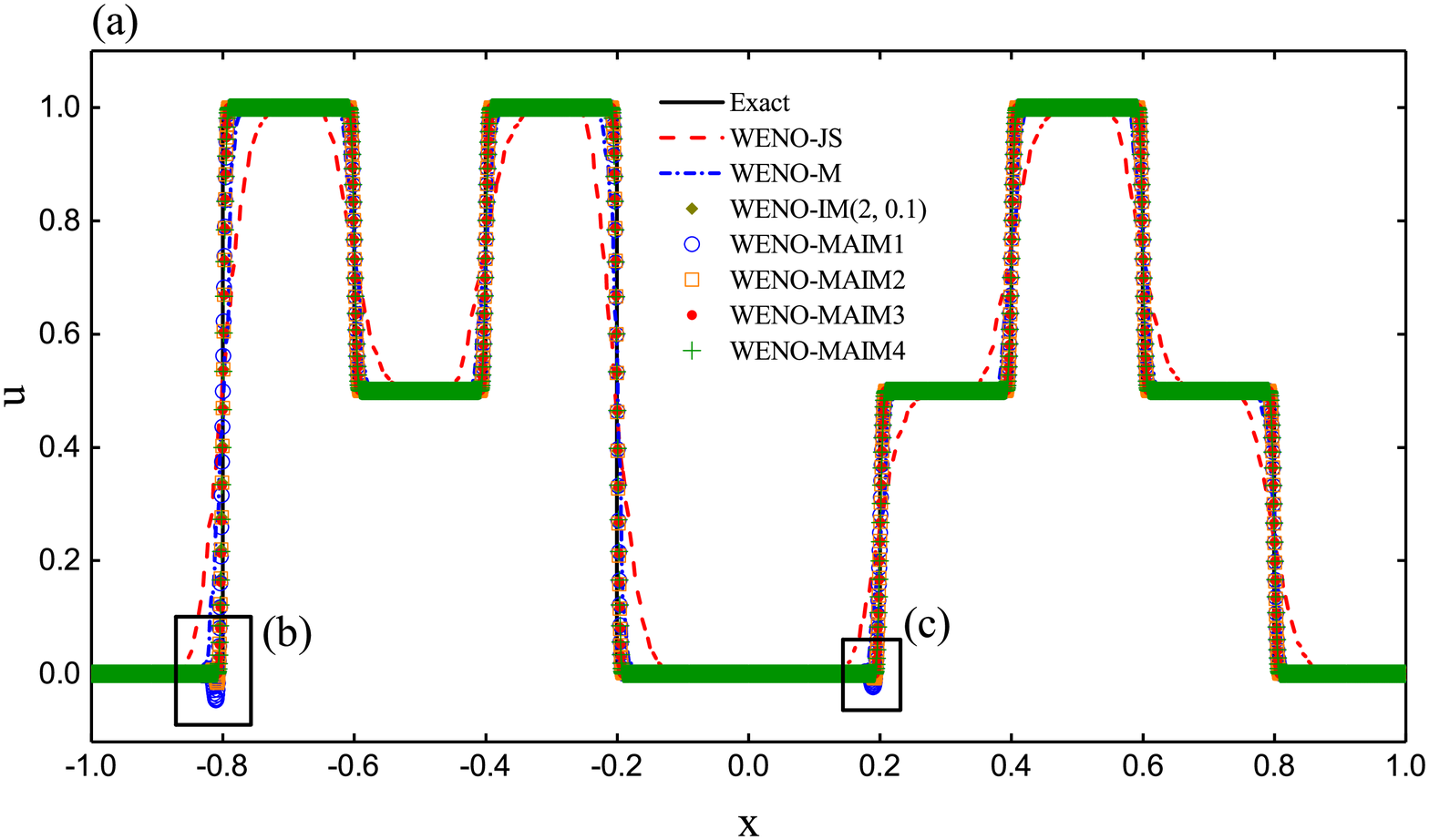}
\includegraphics[height=0.32\textwidth]
{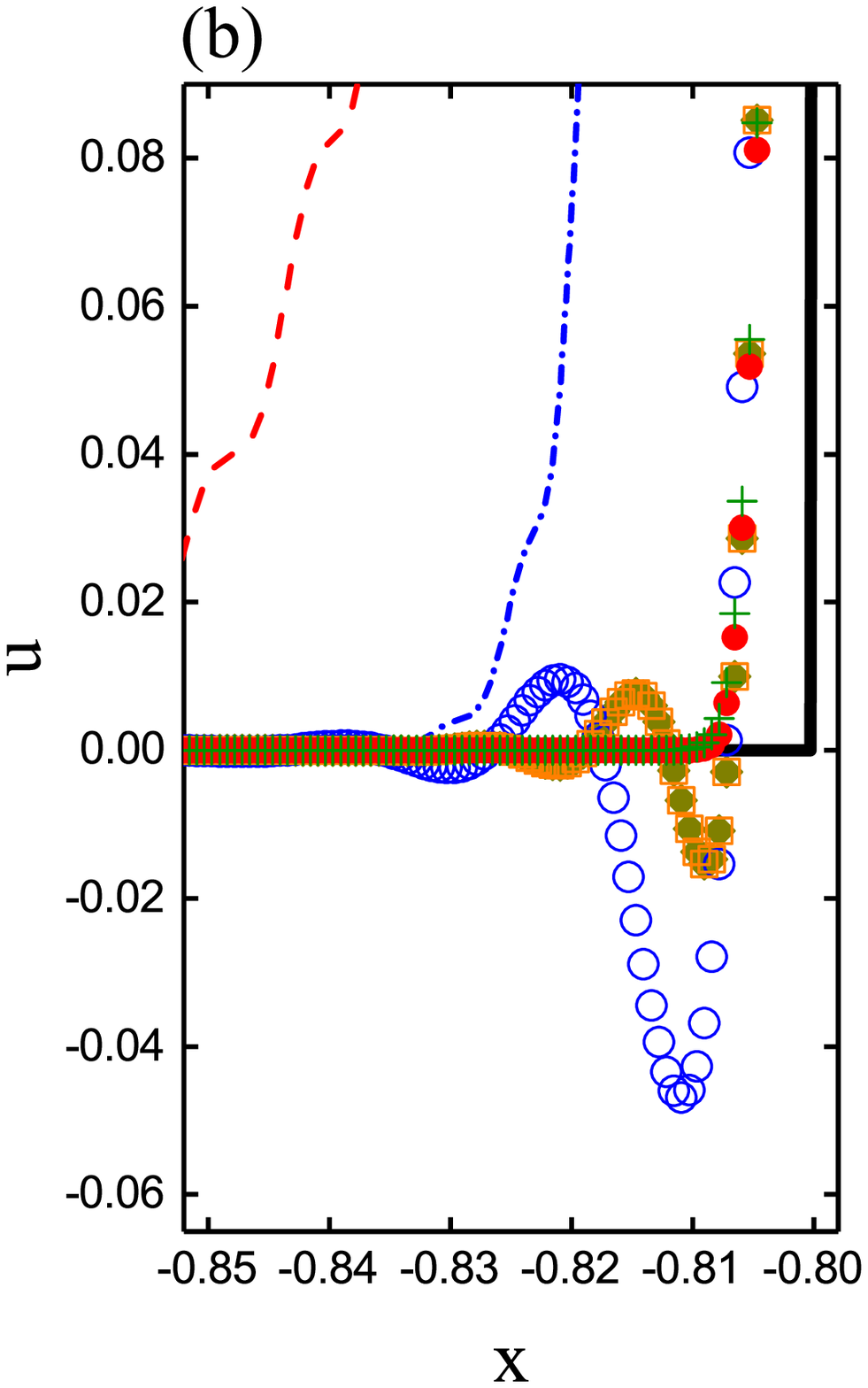}
\includegraphics[height=0.32\textwidth]
{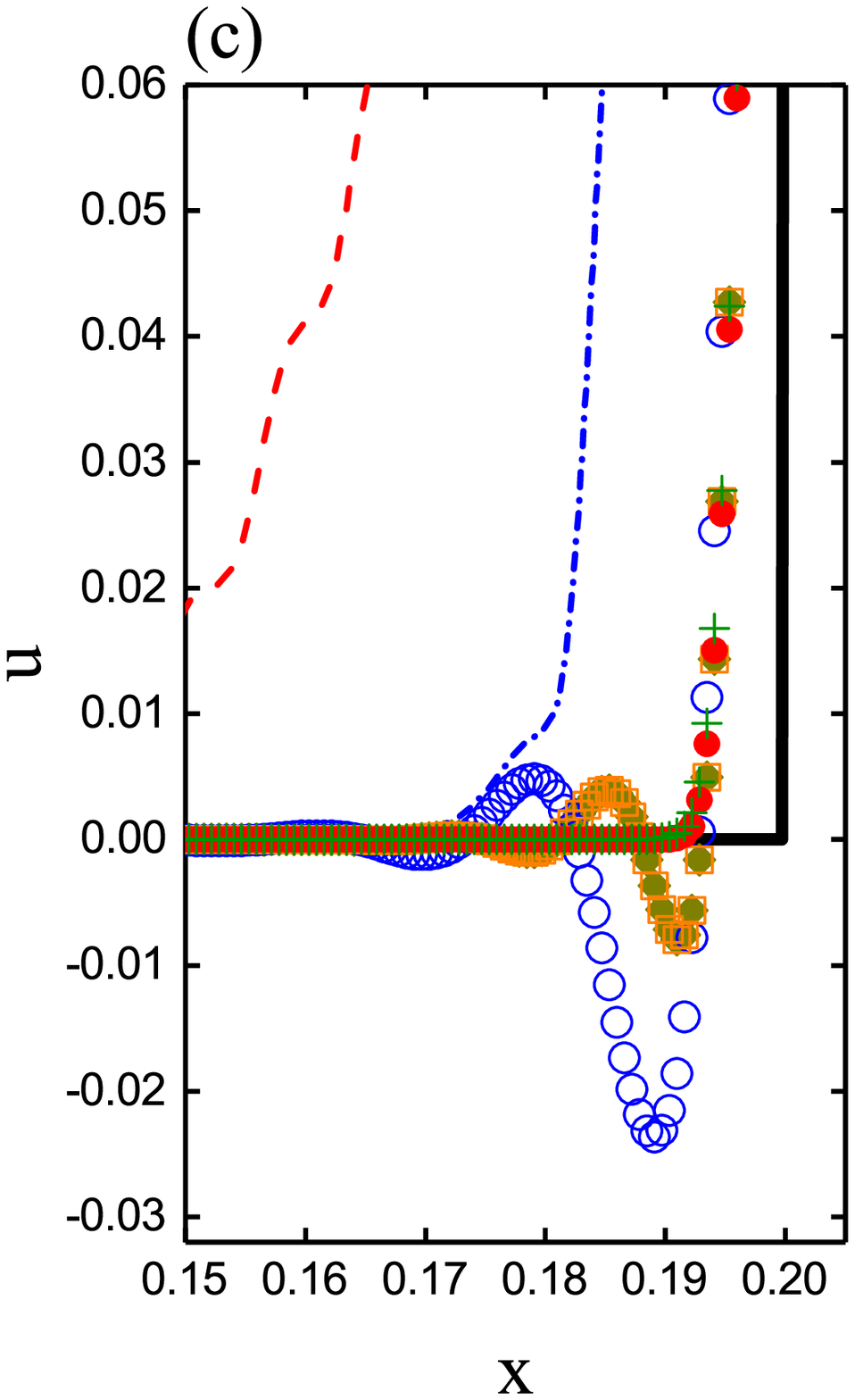}
\caption{Performance of the fifth-order WENO-JS, WENO-M, WENO-IM($2,
0.1$) and WENO-MAIM$i$($i = 1, 2, 3, 4$) schemes for the BiCWP with 
$N=3200$ at long output time $t=200$.}
\label{fig:BiCWP:3200}
\end{figure}

\begin{figure}[ht]
\centering
\includegraphics[height=0.32\textwidth]
{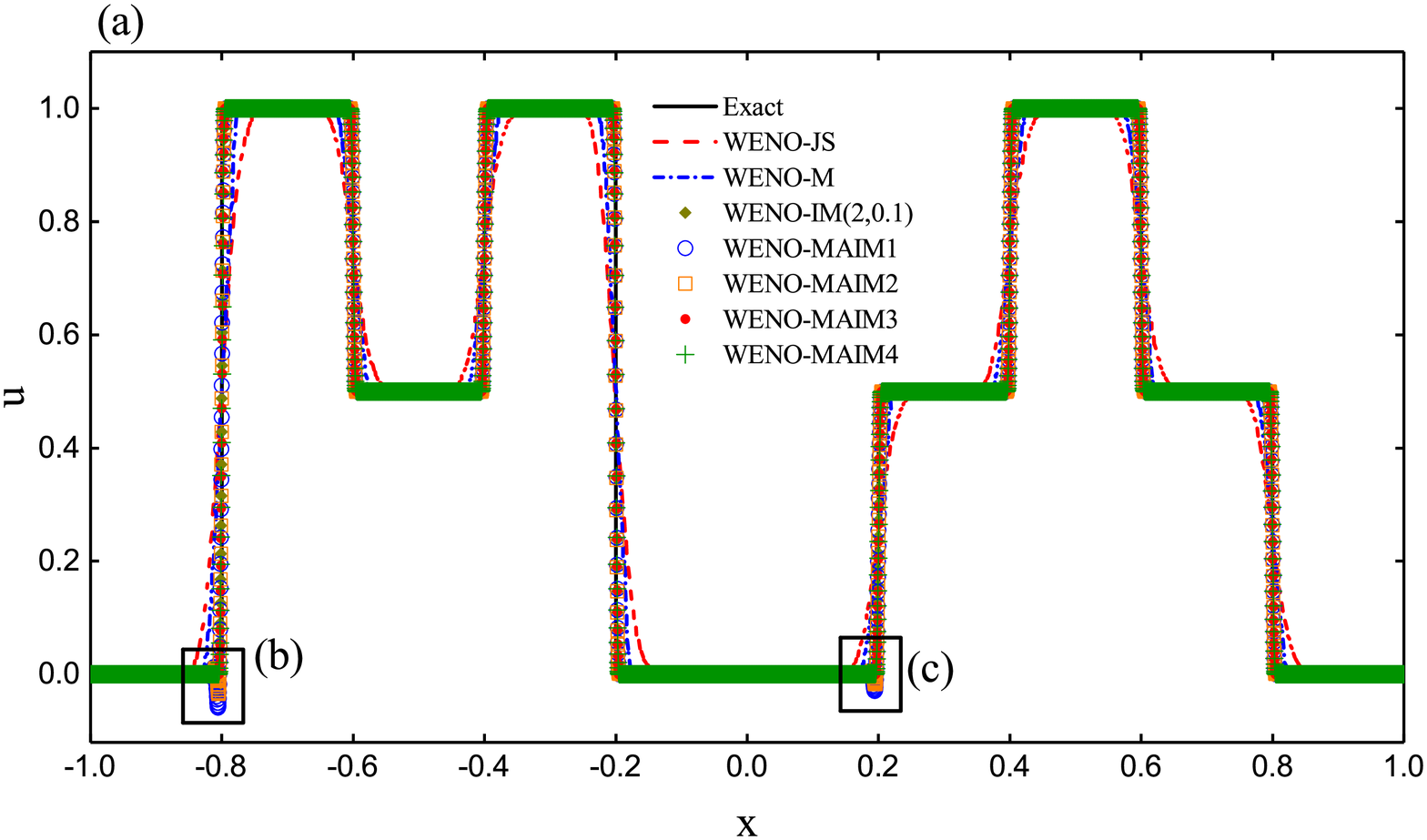}
\includegraphics[height=0.32\textwidth]
{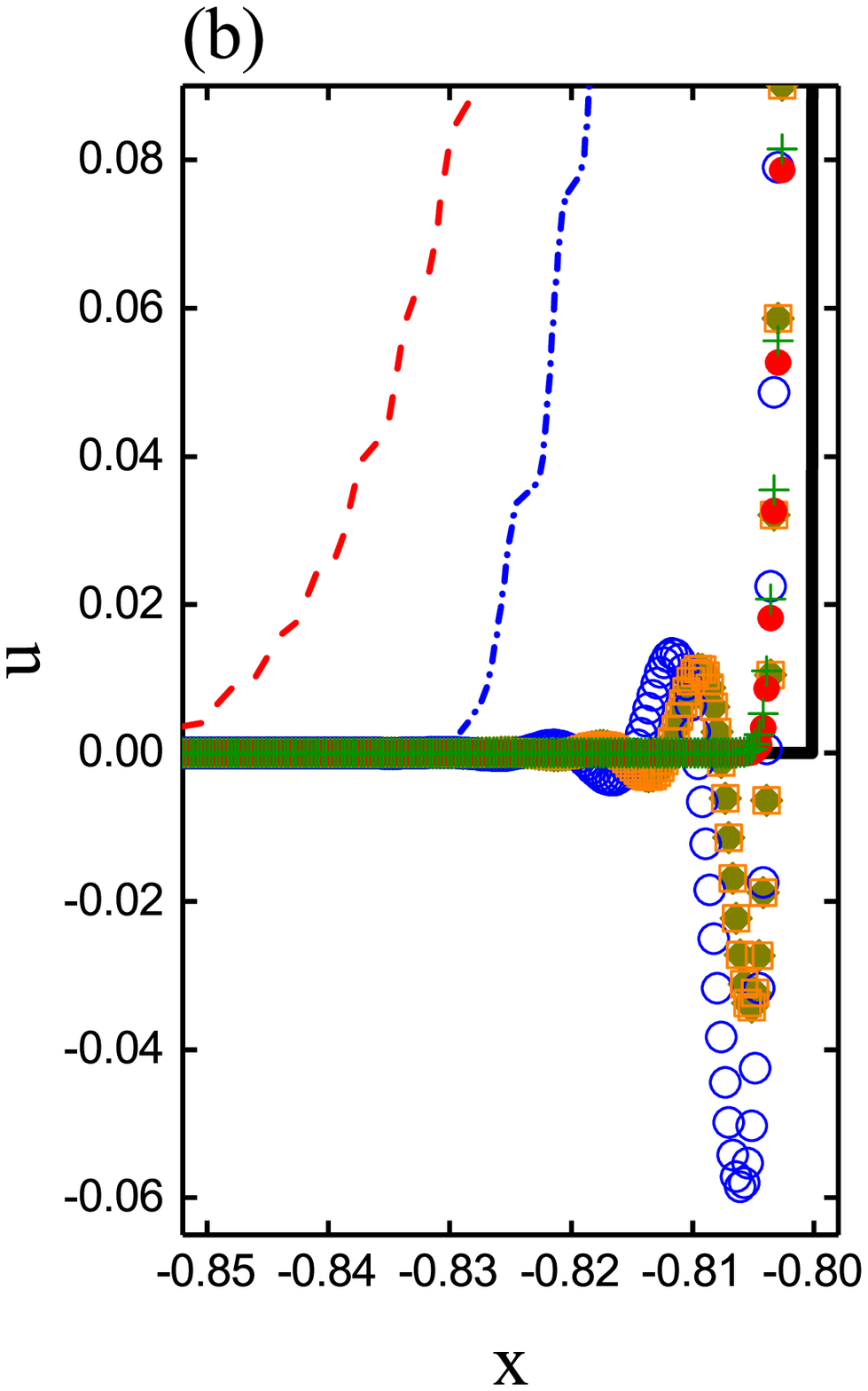}
\includegraphics[height=0.32\textwidth]
{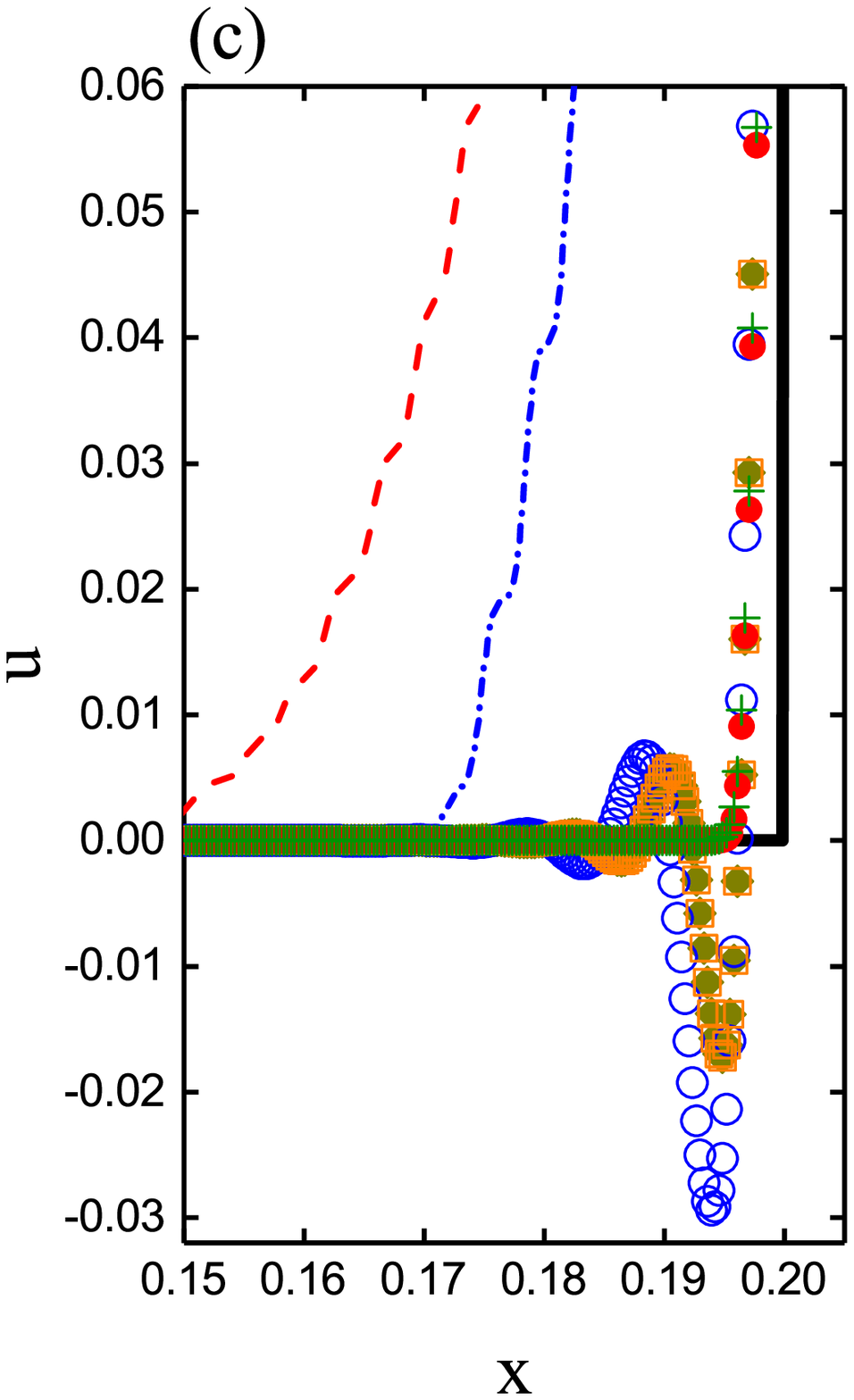}
\caption{Performance of the fifth-order WENO-JS, WENO-M, WENO-IM($2,
0.1$) and WENO-MAIM$i$($i = 1, 2, 3, 4$) schemes for the BiCWP with 
$N=6400$ at long output time $t=200$.}
\label{fig:BiCWP:6400}
\end{figure}

In Figs.\ref{fig:SLP:1600}, \ref{fig:SLP:3200}, \ref{fig:SLP:6400} 
and Figs.\ref{fig:BiCWP:1600}, \ref{fig:BiCWP:3200}, 
\ref{fig:BiCWP:6400}, we have shown the comparison of various 
schemes when $t=200$ and $N=1600, 3200, 6400$ for SLP and BiCWP 
respectively. We observe that: (1) when solving SLP, as the grid 
number increases, the WENO-MAIM1, WENO-MAIM2 and WENO-IM(2, 0.1) 
schemes generate evident spurious oscillations, while the WENO-MAIM3 
and WENO-MAIM4 schemes do not but provide improved resolutions; (2) 
although the WENO-JS and WENO-M schemes obtain numerical results
with improved resolutions for both SLP and BiCWP when the grid 
number increases, their resolutions are significantly lower than 
those of the WENO-MAIM3 and WENO-MAIM4 schemes; (3) for both SLP and 
BiCWP, when the grid number increases, the spurious oscillations of 
the WENO-IM(2, 0.1), WENO-MAIM1 and WENO-MAIM2 schemes get closer to 
the discontinuities, and the amplitudes of these spurious 
oscillations become larger; (4) however, the solutions of the 
WENO-MAIM3 and WENO-MAIM4 schemes get closer to the exact solutions 
without generating spurious oscillations.

\begin{table}[ht]
\begin{scriptsize}
\centering
\caption{Convergence properties of various schemes solving SLP and 
BiCWP at output time $t=200$.}
\label{table_advan:02}
\begin{tabular*}{\hsize}
{@{}@{\extracolsep{\fill}}lllllll@{}}
\hline
\space  &\multicolumn{3}{l}{SLP}  &\multicolumn{3}{l}{BiCWP}  \\
\cline{2-4}  \cline{5-7}
$N$   & $1600$   & $3200$  & $6400$   & $1600$  & $3200$   & $6400$\\
\hline
WENO-JS    & 1.26804e-01(-) & 8.81339e-02(0.5248) & 6.14430e-02
(0.5204)
           & 1.47671e-01(-) & 1.15707e-01(0.3519) & 8.55015e-02
(0.4365)    \\         
{}         & 1.67490e-01(-) & 1.37054e-01(0.2893) & 1.10915e-01
(0.3053)   
           & 1.78065e-01(-) & 1.56061e-01(0.1903) & 1.32873e-01
(0.2321)    \\
{}         & 5.17972e-01(-) & 5.17039e-01(0.0026) & 5.07030e-01
(0.0282)   
           & 5.28594e-01(-) & 5.35169e-01(-0.0275)& 5.13562e-01
(0.0595)    \\
WENO-M     & 4.21095e-02(-) & 3.13290e-02(0.4266) & 2.89979e-02
(0.1116)   
           & 4.45182e-02(-) & 4.70277e-02(-0.0791)& 5.36042e-02
(-0.1888)    \\
{}         & 8.47168e-02(-) & 8.10687e-02(0.0635) & 8.58473e-02
(-0.0826)   
           & 9.35523e-02(-) & 1.01968e-01(-0.1243)& 1.09394e-01
(-0.1014)    \\
{}         & 5.06850e-01(-) & 6.10861e-01(-0.2693)& 5.04121e-01
(0.2771)   
           & 5.06372e-01(-) & 5.71919e-01(-0.1756)& 5.04359e-01
(0.1814)    \\
WENO-IM(2,0.1)  
           & 1.26519e-02(-) & 6.56498e-03(0.9465) & 3.68884e-03
(0.8316)   
           & 2.39427e-02(-) & 1.35781e-02(0.8183) & 7.78055e-03
(0.8033)    \\
{}         & 5.48733e-02(-) & 4.08544e-02(0.4256) & 3.08433e-02
(0.4055)   
           & 7.14167e-02(-) & 5.35510e-02(0.4153) & 4.03989e-02
(0.4066)    \\
{}         & 4.62587e-01(-) & 4.69803e-01(-0.0223)& 4.87134e-01
(-0.0523)   
           & 4.62587e-01(-) & 4.69803e-01(-0.0223)& 4.87134e-01
(-0.0523)    \\
WENO-MAIM1 & 1.35847e-02(-) & 7.18556e-03(0.9188) & 3.93600e-03
(0.8684)   
           & 2.55329e-02(-) & 1.44669e-02(0.8196) & 8.15452e-03
(0.8271)    \\
{}         & 5.64038e-02(-) & 4.20972e-02(0.4221) & 3.14976e-02
(0.4185)   
           & 7.29915e-02(-) & 5.48110e-02(0.4133) & 4.10730e-02
(0.4163)    \\
{}         & 4.82675e-01(-) & 4.98673e-01(-0.0470)& 5.10410e-01
(-0.0336)   
           & 4.82672e-01(-) & 4.98672e-01(-0.0470)& 5.10411e-01
(-0.0336)    \\
WENO-MAIM2 & 1.26519e-02(-) & 6.56498e-03(0.9465) & 3.68884e-03
(0.8316)   
           & 2.39427e-02(-) & 1.35781e-02(0.8183) & 7.78055e-03
(0.8033)    \\
{}         & 5.48733e-02(-) & 4.08544e-02(0.4256) & 3.08433e-02
(0.4055)   
           & 7.14167e-02(-) & 5.35510e-02(0.4153) & 4.03989e-02
(0.4066)    \\
{}         & 4.62587e-01(-) & 4.69803e-01(-0.0223)& 4.87134e-01
(-0.0523)   
           & 4.62587e-01(-) & 4.69803e-01(-0.0223)& 4.87134e-01
(-0.0523)    \\
WENO-MAIM3 & 1.23903e-02(-) & 6.43177e-03(0.9459) & 3.49736e-03
(0.8789)   
           & 2.39789e-02(-) & 1.34355e-02(0.8357) & 7.53547e-03
(0.8343)    \\
{}         & 5.48637e-02(-) & 4.08071e-02(0.4270) & 3.05283e-02
(0.4187)   
           & 7.14302e-02(-) & 5.35161e-02(0.4166) & 4.01067e-02
(0.4161)    \\
{}         & 4.62409e-01(-) & 4.66452e-01(-0.0126)& 4.70087e-01
(-0.0112)   
           & 4.62408e-01(-) & 4.66453e-01(-0.0126)& 4.70087e-01
(-0.0112)    \\
WENO-MAIM4 & 1.23934e-02(-) & 6.53848e-03(0.9225) & 3.54066e-03
(0.8849)   
           & 2.44954e-02(-) & 1.36623e-02(0.8423) & 7.63475e-03
(0.8395)    \\
{}         & 5.53035e-02(-) & 4.10626e-02(0.4295) & 3.06775e-02
(0.4206)   
           & 7.19931e-02(-) & 5.38487e-02(0.4189) & 4.03025e-02
(0.4180)    \\
{}         & 4.62634e-01(-) & 4.66590e-01(-0.0123)& 4.70166e-01
(-0.0110)   
           & 4.62641e-01(-) & 4.66590e-01(-0.0123)& 4.70166e-01
(-0.0110)    \\
\hline
\end{tabular*}
\end{scriptsize}
\end{table}

\subsection{One-dimensional Euler system}
In this subsection, we present numerical examples of the 
one-dimensional Euler system for gas dynamics in the following 
conservation form with different initial and boundary conditions
\begin{equation}
\mathbf{U}_t + \mathbf{F}\big( \mathbf{U} \big)_x = 0,
\label{1DEulerEquations}
\end{equation}
where $\mathbf{U} = \big( \rho, \rho u, E \big)^{\mathrm{T}}, \mathbf
{F}\big( \mathbf{U} \big) = \big( \rho u, \rho u^{2} + p, u(E + p) 
\big)^{\mathrm{T}}$, and $\rho, u, p$ and $E$ are the density, 
velocity, pressure and total energy, respectively. The Euler system 
(\ref{1DEulerEquations}) is closed by the equation of state for an 
ideal polytropic gas, which is given by
\begin{equation*}
p = (\gamma - 1)\Big( E - \dfrac{1}{2}\rho u^{2} \Big), \quad \gamma 
= 1.4
\end{equation*}
The finite volume version of the characteristic-wise one-dimensional 
WENO procedure is employed, and we refer to \cite{WENOoverview} for 
details. In all calculations of this subsection, the CFL number is 
set to be $0.1$ and the reference solutions are obtained by using the
WENO-JS scheme with the resolution of $N = 10000$.

\begin{example}
\bf{(Woodward-Colella interacting blast wave problem)} \rm{We solve 
the one-dimensional blast waves interaction problem of Woodward and 
Colella \cite{interactingBlastWaves-Woodward-Colella}, which has the 
following initial condition with reflective boundary conditions}
\label{Euler3}
\end{example}
\begin{equation}
\big( \rho, u, p \big)(x, 0) =\left\{
\begin{array}{ll}
(1, 0, 1000),   & x \in [0, 0.1),   \\
(1, 0, 0.01), & x \in [0.1, 0.9), \\
(1, 0, 100), & x \in [0.9, 1.0].
\end{array}\right.
 \label{initial_1DEuler3}
\end{equation}
This problem is run with $N = 400$ uniform cells till $t = 0.038$ 
by the WENO-JS, WENO-M, WENO-IM($2,0.1$) and 
WENO-MAIM$i$($i=1,2,3,4$) schemes. 

\begin{figure}[ht]
\centering
\includegraphics[height=0.42\textwidth]
{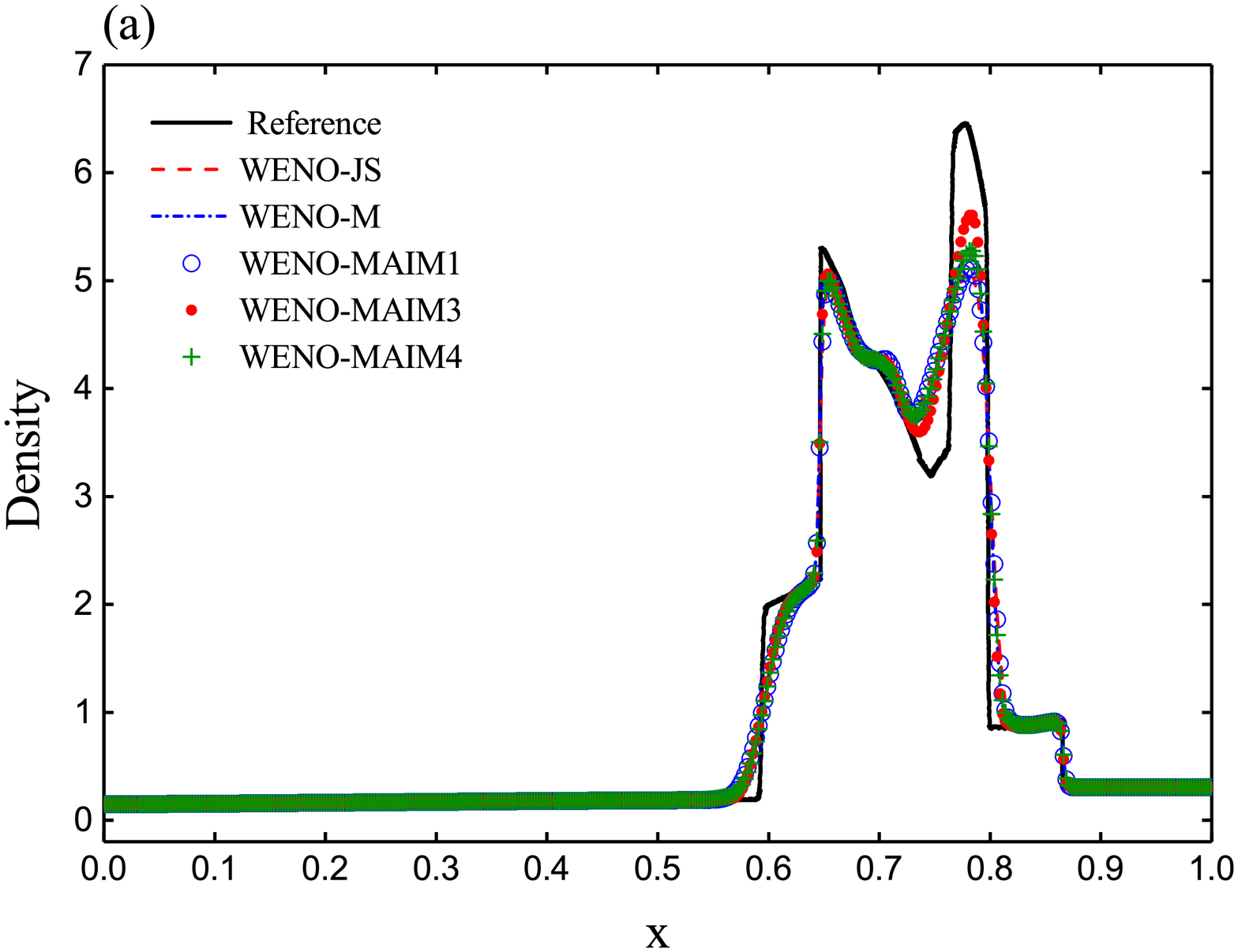}
\includegraphics[height=0.42\textwidth]
{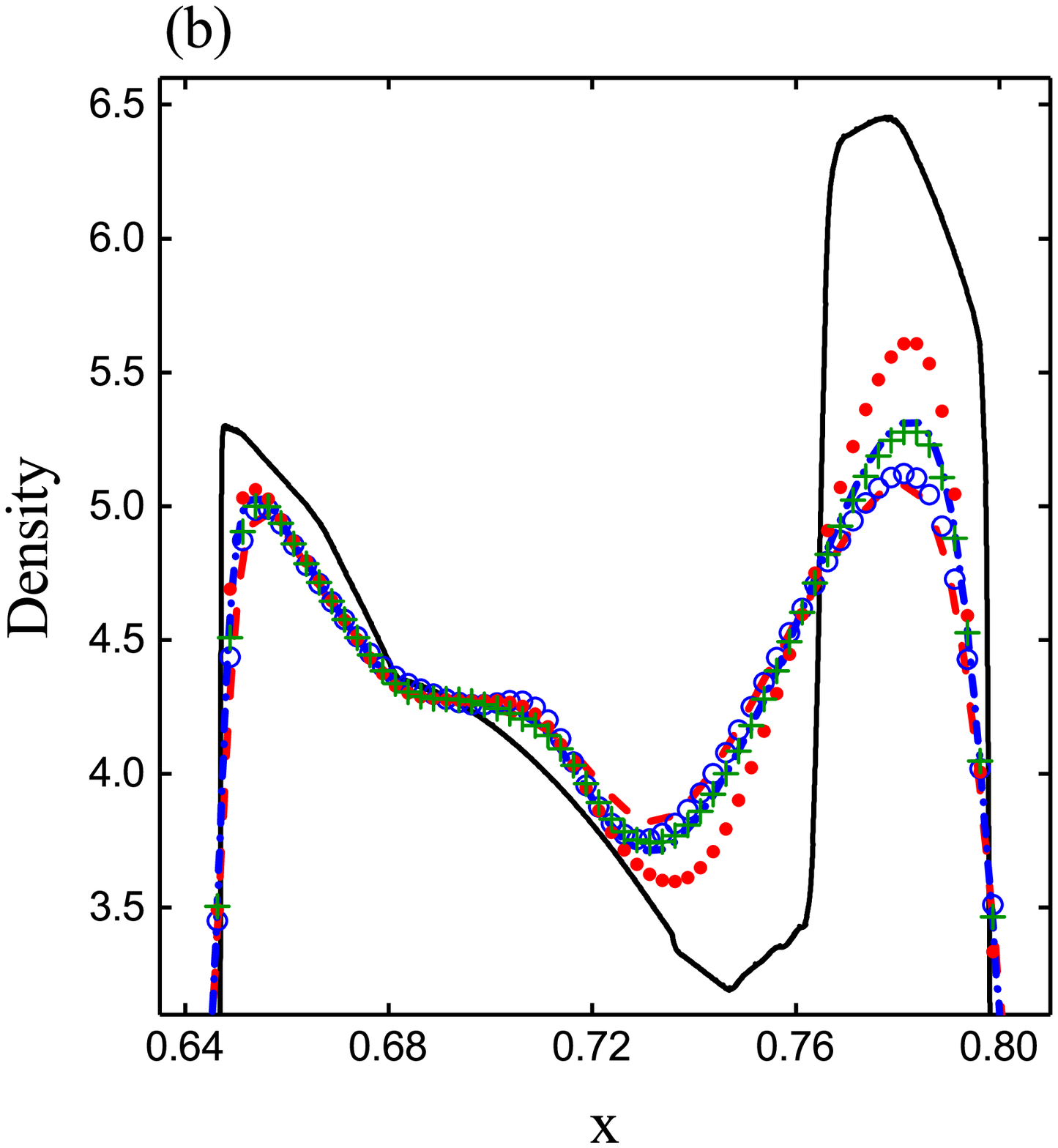}
\caption{Performance of fifth-order WENO-JS, WENO-M and WENO-MAIM$i$
($i = 1, 3, 4$) for Euler system (\ref{1DEulerEquations}) with 
initial condition (\ref{initial_1DEuler3}).}
\label{fig:exampleEuler3}
\end{figure}

Surprisingly, for this specific problem, the WENO-MAIM2 and 
WENO-IM($2, 0.1$) schemes blow up in our calculations, but the other 
schemes perform well. The results are given in 
Fig. \ref{fig:exampleEuler3}. It can be observed directly from 
Fig. \ref{fig:exampleEuler3} that the WENO-MAIM$3$ scheme provides 
the best performance, followed by the WENO-MAIM4 scheme and then the 
WENO-M scheme. The WENO-MAIM$1$ scheme performs worse than most of 
the considered schemes but slightly better than the WENO-JS scheme.

\begin{example}
\bf{(Titarev-Toro shock-entropy wave interaction problem)} \rm{We 
solve the one-dimensional shock-entropy wave interaction problem 
proposed by Titarev and Toro \cite{Titarev-Toro-1,Titarev-Toro-2,
Titarev-Toro-3} which is a more severe version of the classic 
Shu-Osher shock-entropy wave interaction problem \cite{ENO-Shu1989}.
This example is commonly used to test the ability to capture both 
shocks and short wavelength oscillations of the numerical schemes, 
and the amplitudes of the short wavelength oscillations are a 
measure of the numerical viscosity of the schemes. The initial 
condition is given by}
\label{Euler2:b}
\end{example}
\begin{equation}
\big( \rho, u, p \big)(x, 0) =\left\{
\begin{array}{ll}
(1.515695, 0.5233346, 1.80500), & x \in [-5.0, -4.5], \\
(1.0 + 0.1\sin(20\pi x), 0, 1), & x \in [-4.5, 5.0].
\end{array}\right.
 \label{initial_1DEuler2:b}
\end{equation} 
It is run with $N = 1500$ uniform cells till $t = 5.0$ and the CFL 
number is set to be $0.4$. 

\begin{figure}[ht]
\centering
\includegraphics[height=0.45\textwidth]
{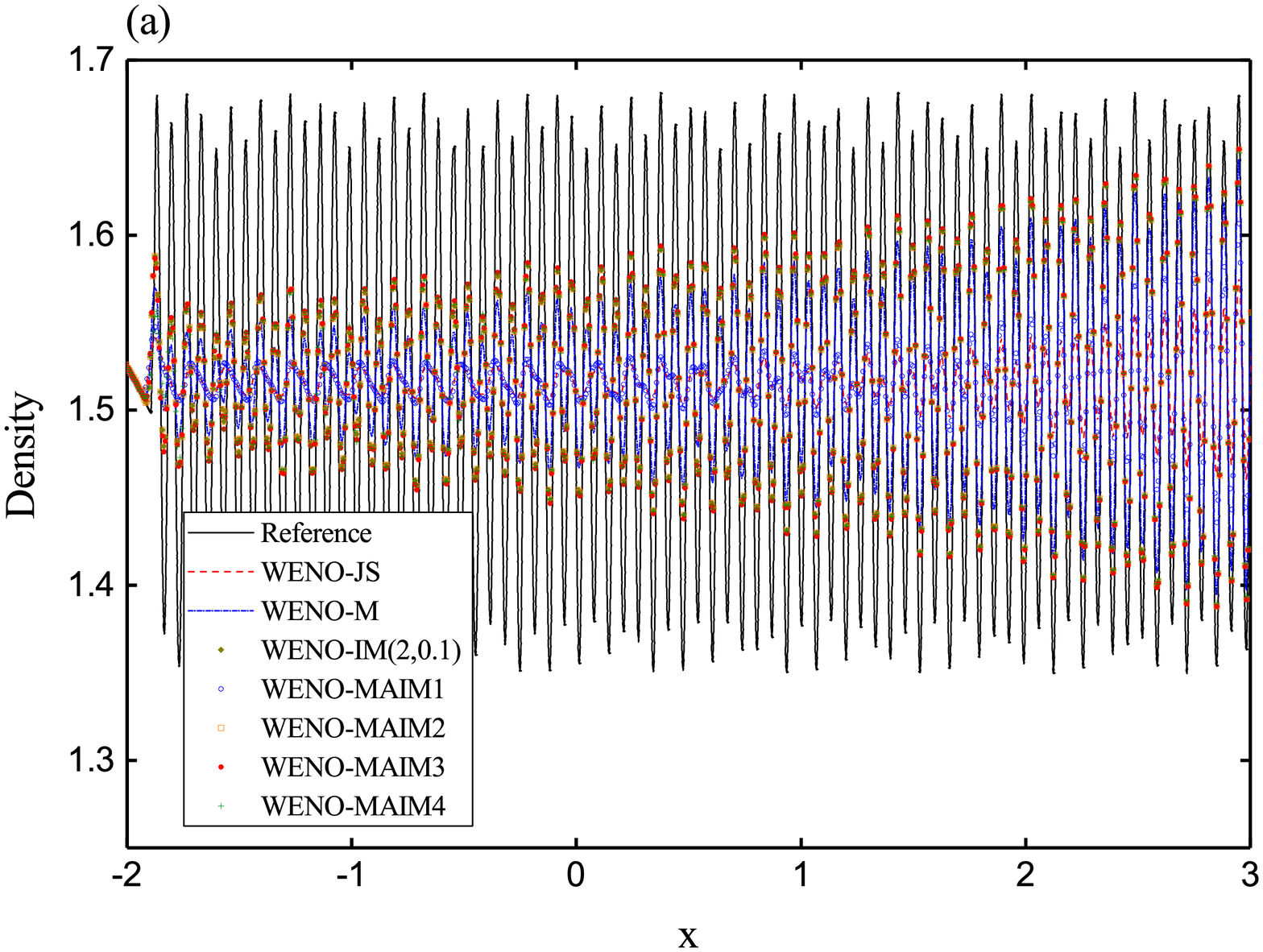}
\includegraphics[height=0.45\textwidth]
{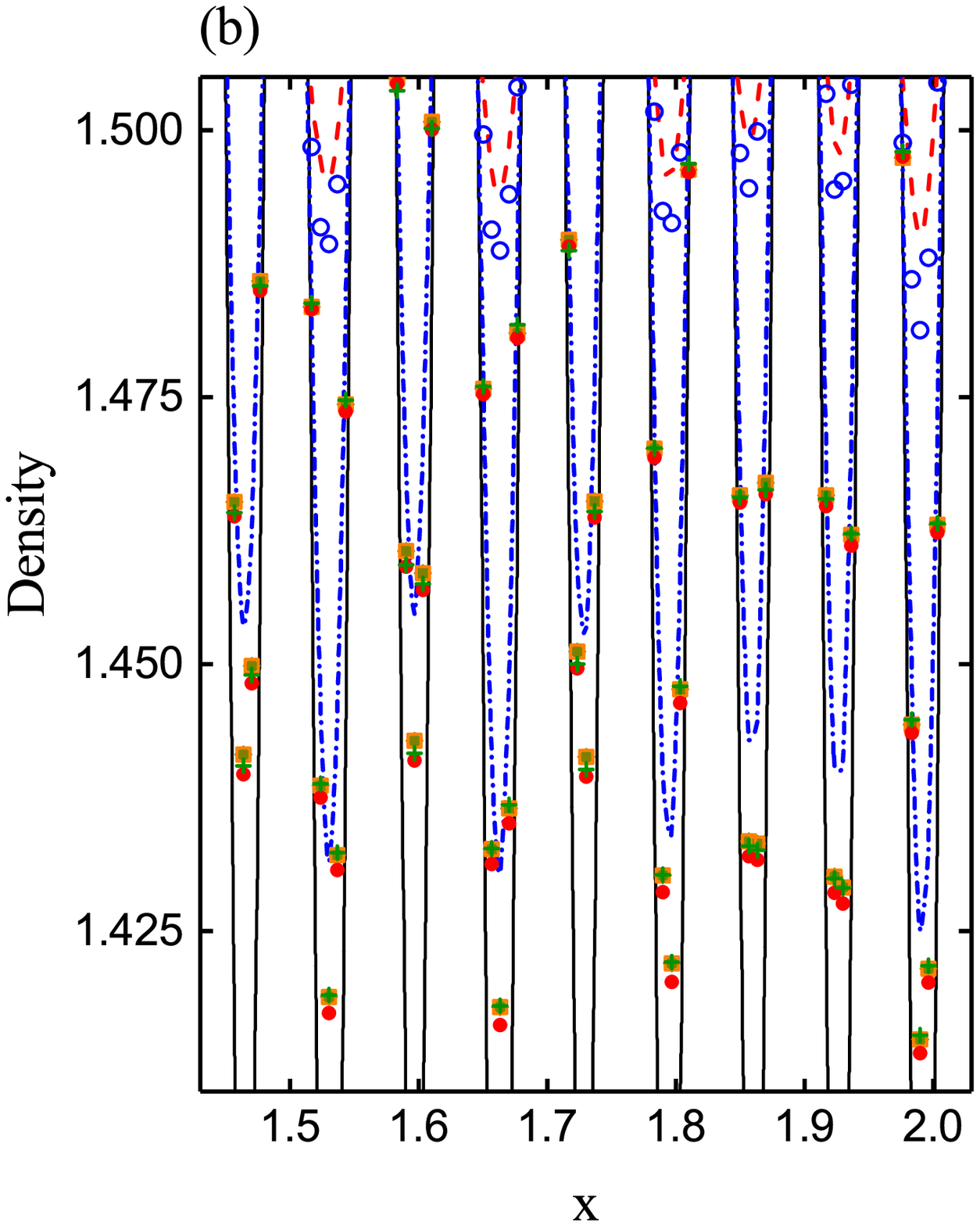}
\caption{Performance of fifth-order WENO-JS, WENO-M, WENO-IM($2,0.1$
) and WENO-MAIM$i$($i = 1, 2, 3, 4$) for Euler system 
(\ref{1DEulerEquations}) with initial condition 
(\ref{initial_1DEuler2:b}).}
\label{fig:exampleEuler2:b}
\end{figure}

The results are shown in Fig. \ref{fig:exampleEuler2:b}. It can be 
seen that the WENO-MAIM3 scheme results in slightly greater 
amplitudes compared to the WENO-MAIM2, WENO-MAIM4 and 
WENO-IM($2, 0.1$) schemes which produce higher amplitudes of the 
short wavelength oscillations compared to the WENO-M scheme. The 
WENO-MAIM1 scheme performs slightly better than the WENO-JS scheme 
which performs worst among all the considered WENO schemes for this 
case.

\subsection{Two-dimensional Euler system}
\label{subsec:examples_2D_Euler}
In this subsection, we extend the considered WENO schemes to 
calculate the two-dimensional compressible Euler system of gas 
dynamics taking the following strong conservation form
\begin{equation}
\mathbf{U}_t + \mathbf{F}\big( \mathbf{U} \big)_x +
\mathbf{G}\big( \mathbf{U} \big)_y= 0,
\label{2DEulerEquations}
\end{equation}
where $\mathbf{U} = \big(\rho, \rho u, \rho v, E\big)^{\mathrm{T}}, 
\mathbf{F}\big( \mathbf{U} \big) = \big( \rho u, \rho u^{2} + p, \rho
uv, u(E + p) \big)^{\mathrm{T}}, \mathbf{G}\big( \mathbf{U} \big) = 
\big( \rho v, \rho vu, \rho v^{2} + p, v(E + p) \big)^{\mathrm{T}}$, 
and $\rho, u, v, p$ and $E$ are the density, component of velocity 
in the $x$ and $y$ coordinate directions, pressure and total energy, 
respectively. The Euler system (\ref{2DEulerEquations}) is closed by 
the equation of state for an ideal polytropic gas given by
\begin{equation*}
p = (\gamma - 1)\Big( E - \dfrac{1}{2}\rho (u^{2} + v^{2}) \Big), 
\quad \gamma = 1.4
\end{equation*}
In all numerical examples of this subsection, the CFL number is set 
to be $0.5$. As the WENO-MAIM2 scheme gives exactly the same 
solutions as the WENO-IM(2, 0.1) scheme, we only show the results of 
the WENO-IM(2, 0.1) scheme for brevity in the following presentation.

\begin{example}
\bf{(2D Riemann problem)} 
\rm{We calculate the 2D Riemann problem proposed by 
\cite{Riemann-2D-01,Riemann2D-02}. This problem is done over a unit 
square domain $[0,1] \times[0,1]$. There are many different 
configurations for 2D Riemann problem \cite{Riemann2D-02}, and we 
focus on the one solved in \cite{WENO-PPM5,Riemann2D-02}. It 
initially involves the constant states of flow variables over each 
quadrant which is got by dividing the computational domain using 
lines $x = 0.5$ and $y = 0.5$. The initial constant states are given 
by}
\label{ex:Riemann2D}
\end{example}
\begin{equation*}
\big( \rho, u, v, p \big)(x, y, 0) = \left\{
\begin{aligned}
\begin{array}{ll}
(1.0, 0.0, -0.3, 1.0),       & 0.5 \leq x \leq 1.0, 
                               0.5 \leq y \leq 1.0, \\
(2.0, 0.0, 0.3, 1.0),        & 0.0 \leq x \leq 0.5, 
                               0.5 \leq y \leq 1.0, \\
(1.0625, 0.0, 0.8145, 0.4),  & 0.0 \leq x \leq 0.5, 
                               0.0 \leq y \leq 0.5, \\
(0.5313, 0.0, 0.4276, 0.4),  & 0.5 \leq x \leq 1.0, 
                             0.0 \leq y \leq 0.5. \\
\end{array}
\end{aligned}
\right.
\label{eq:initial_Euer2D:Riemann2D}
\end{equation*}
The transmission boundary condition is used on all boundaries. We 
evolve the initial data until time $t=0.3$ using considered WENO 
schemes with a uniform mesh size of $1200\times1200$.

\begin{figure}[ht]
\centering
  \includegraphics[height=0.29\textwidth]
  {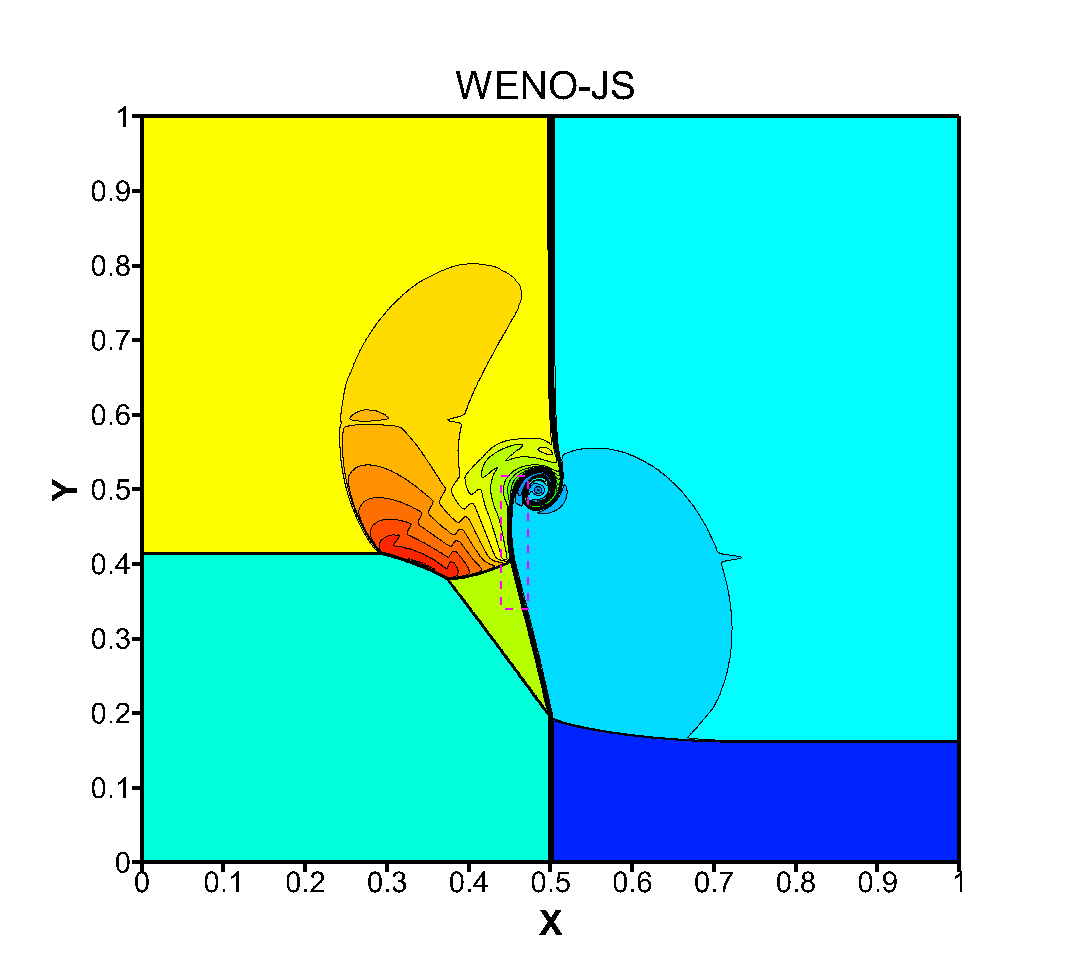}
  \includegraphics[height=0.29\textwidth]
  {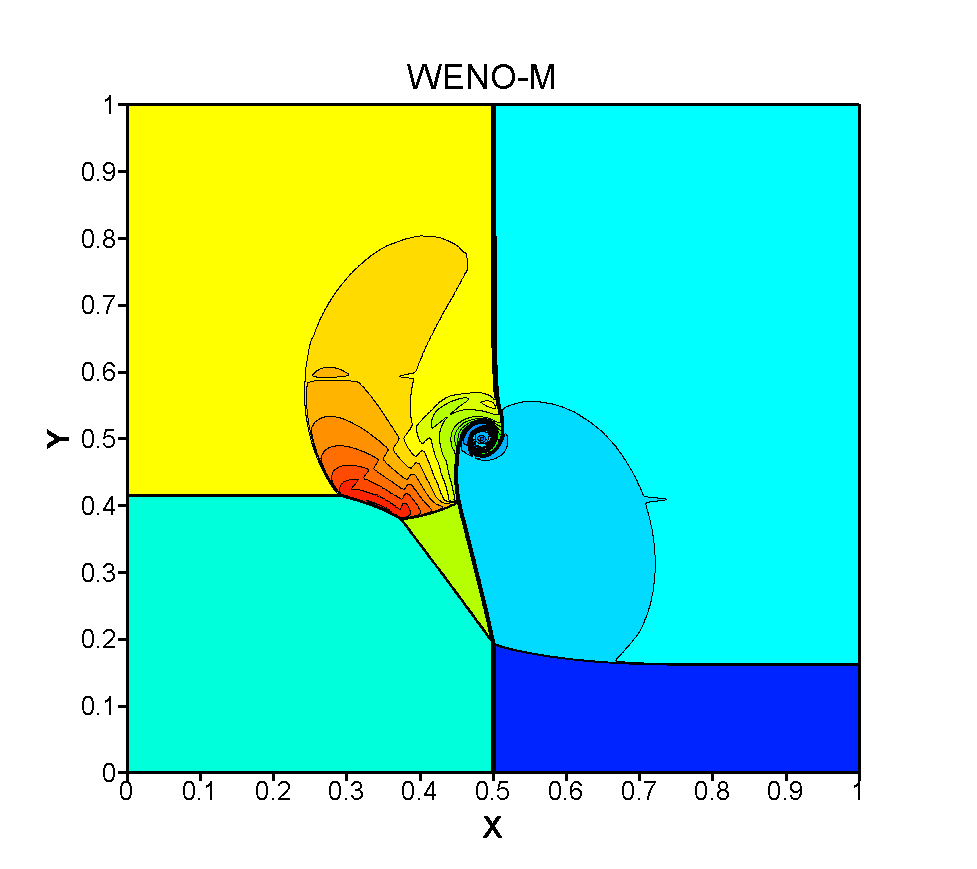}
  \includegraphics[height=0.29\textwidth]
  {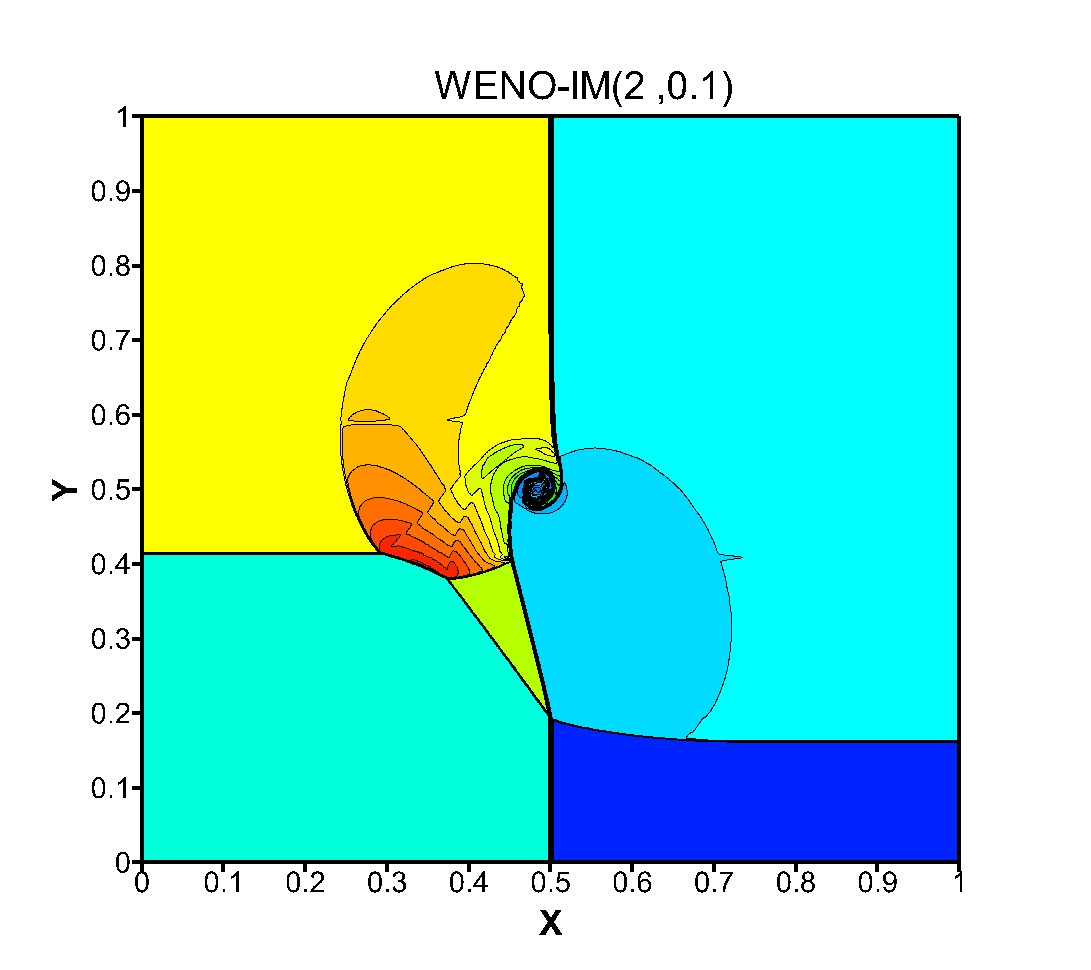}\\
  \includegraphics[height=0.29\textwidth]
  {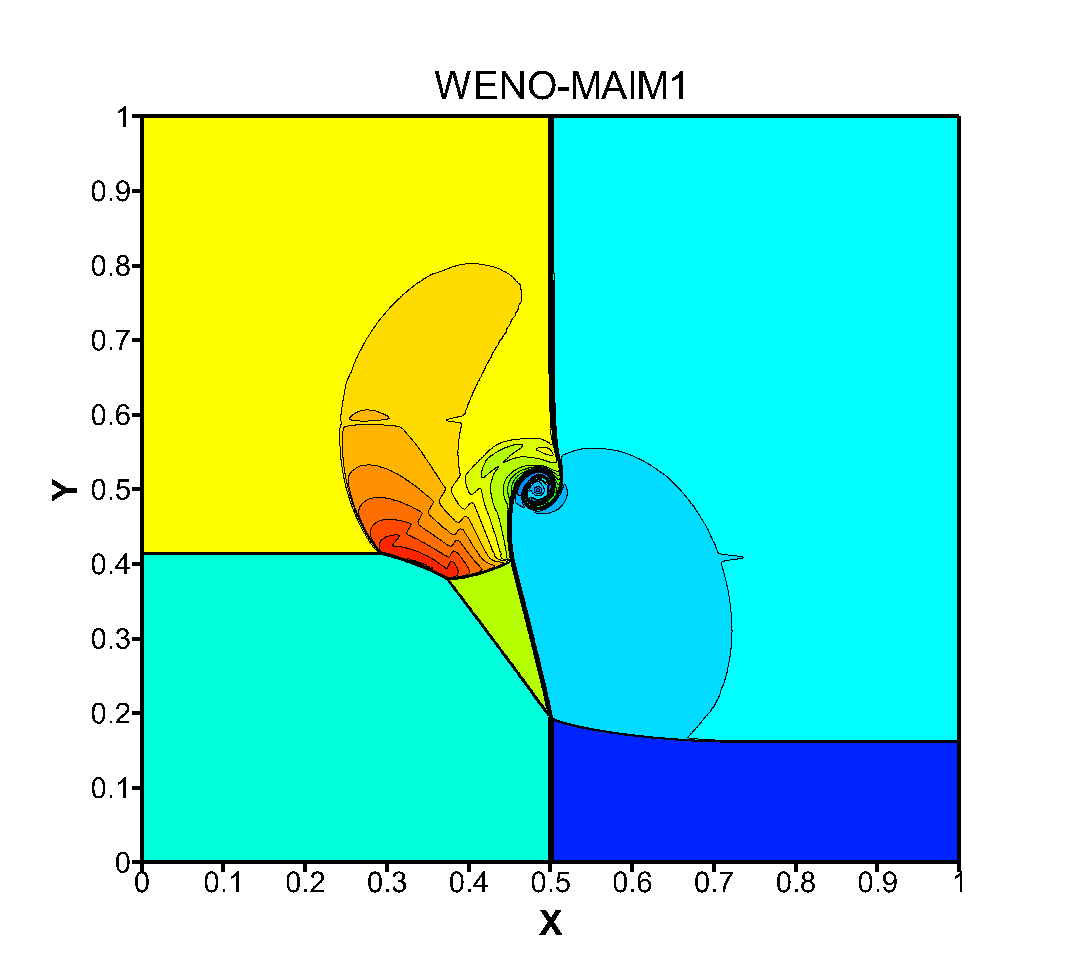}
  \includegraphics[height=0.29\textwidth]
  {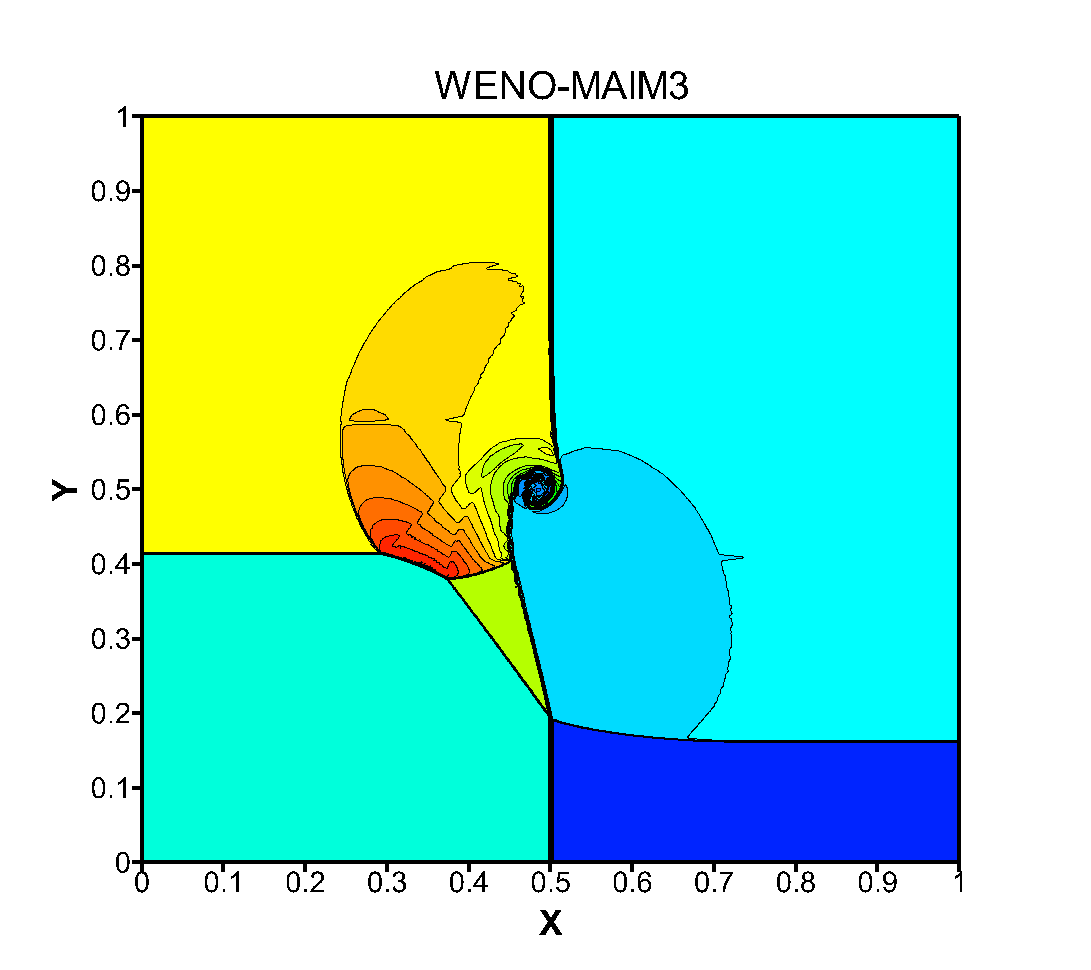}
  \includegraphics[height=0.29\textwidth]
  {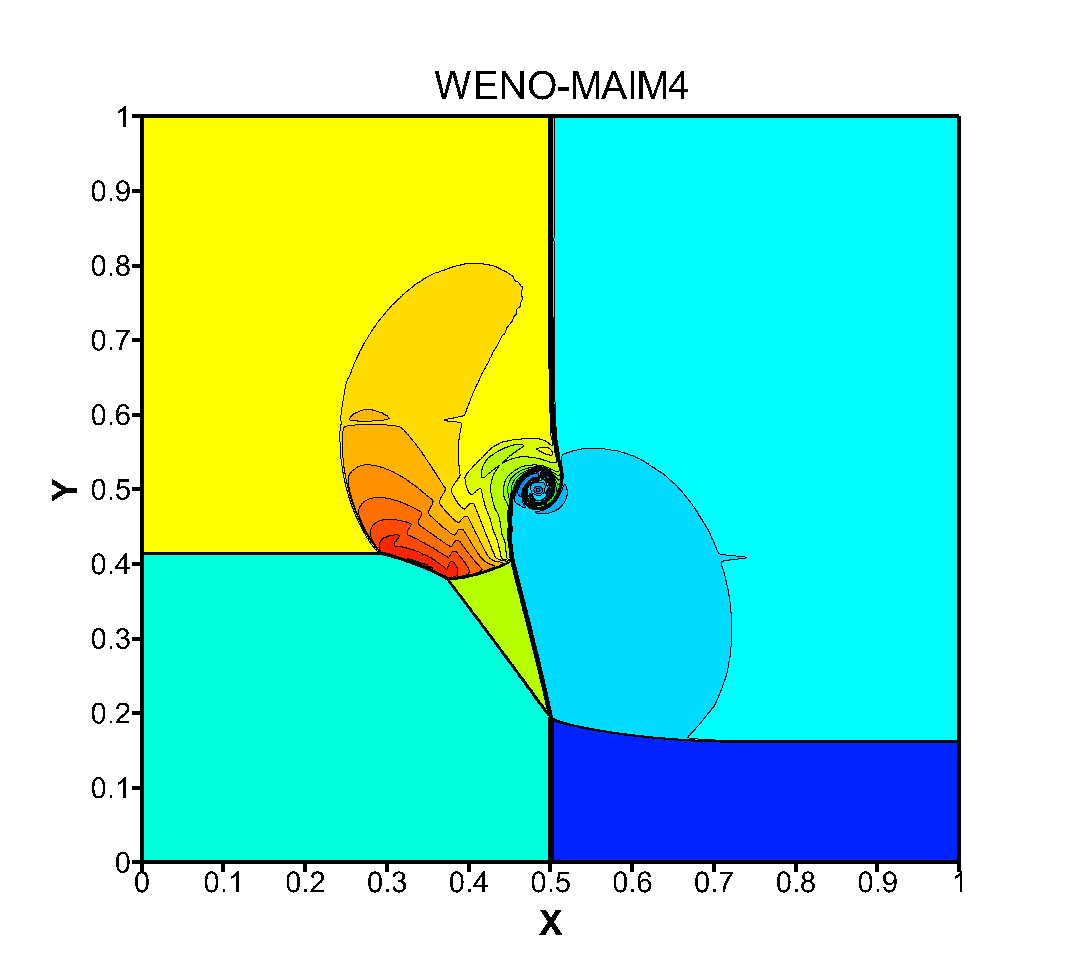}    
\caption{Density plots for the 2D Riemann problem using $30$ contour
lines with range from $0.5$ to $2.4$, computed using the WENO-JS, 
WENO-M, WENO-MAIM$i(i=1,2,3,4)$ schemes at output time $t = 0.3$ 
with a uniform mesh size of $1200\times1200$.}
\label{fig:ex:Riemann2D}
\end{figure}

\begin{figure}[ht]
\centering
  \includegraphics[height=0.28\textwidth]
  {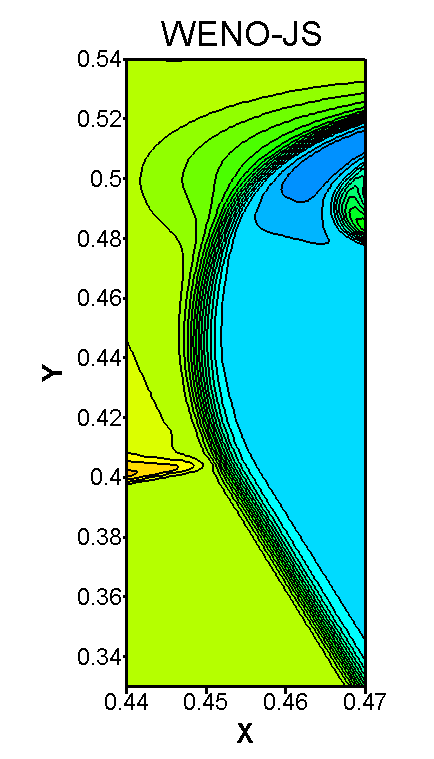}
  \includegraphics[height=0.28\textwidth]
  {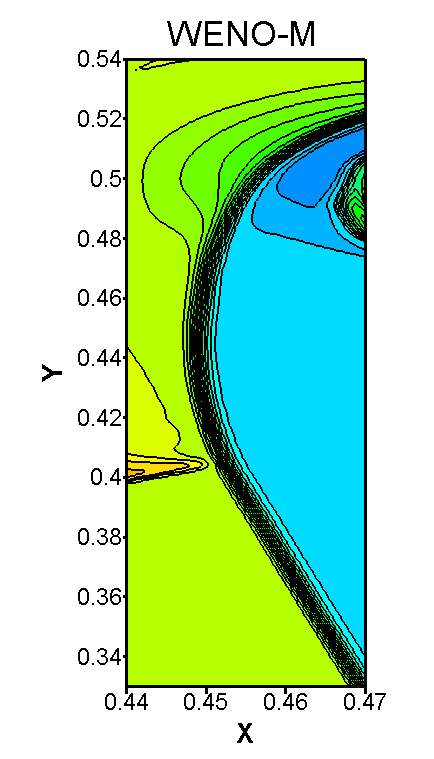}
  \includegraphics[height=0.28\textwidth]
  {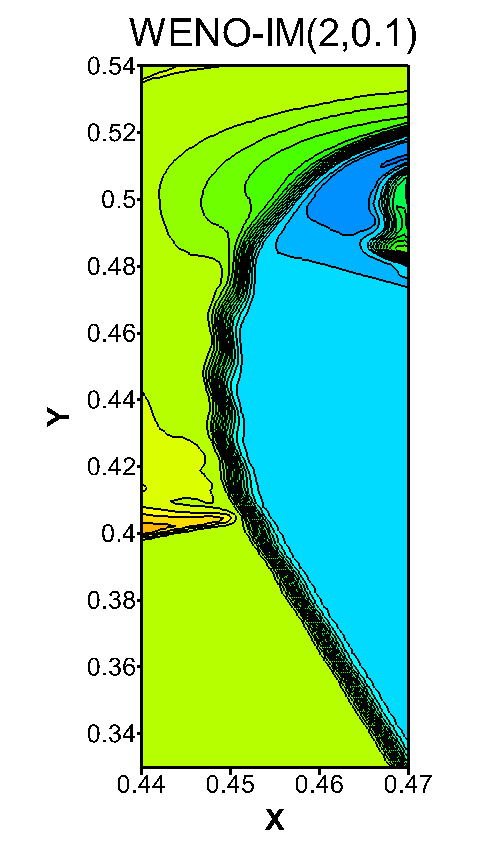}
  \includegraphics[height=0.28\textwidth]
  {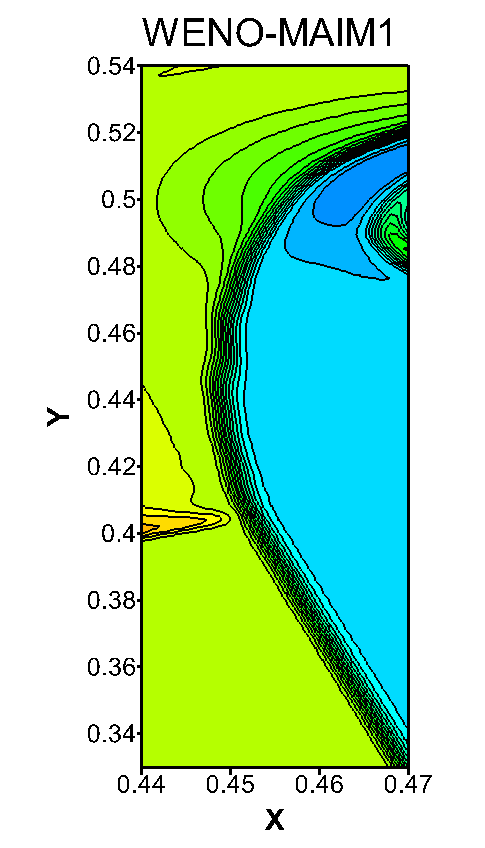}
  \includegraphics[height=0.28\textwidth]
  {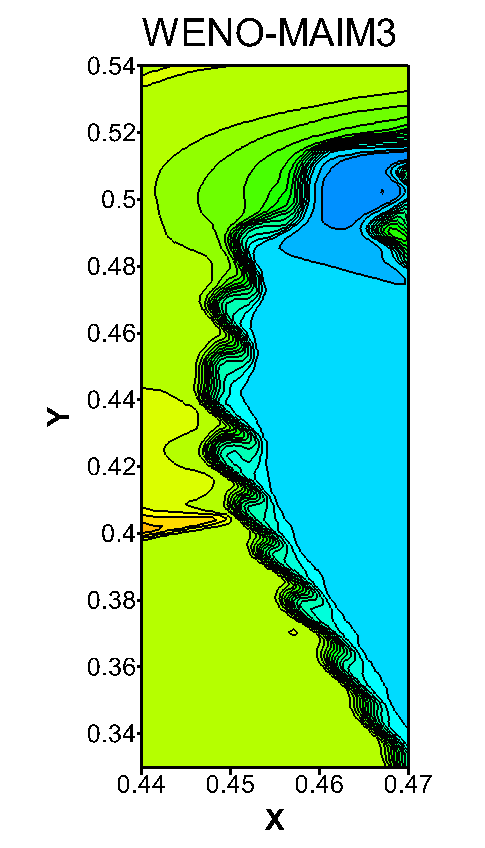}
  \includegraphics[height=0.28\textwidth]
  {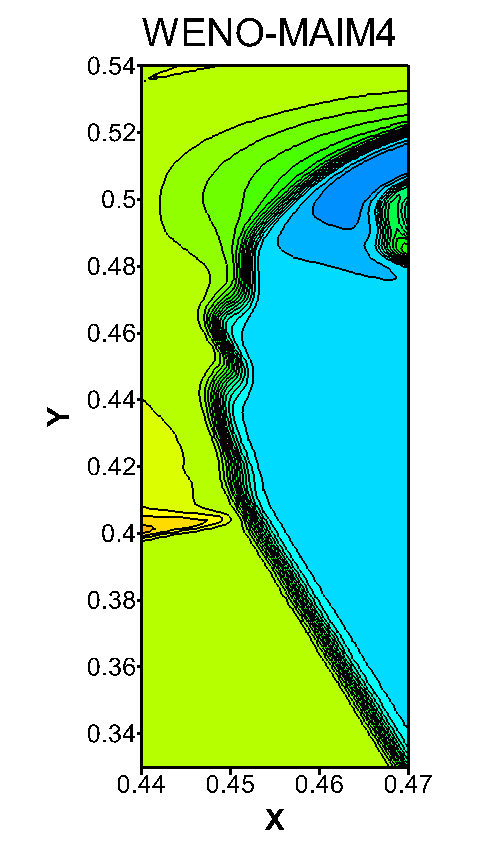}    
\caption{The zoomed-in density plots for the 2D Riemann problem.}
\label{fig:ex:Riemann2D-ZoomedIn}
\end{figure}

The numerical results of density obtained using the WENO-JS, WENO-M, 
WENO-IM(2, 0.1), WENO-MAIM1, WENO-MAIM3 and WENO-MAIM4 schemes have 
been shown in Fig. \ref{fig:ex:Riemann2D}. All these considered WENO 
schemes can captured the main structure of the solution. 
Furthermore, as mentioned in \cite{WENO-PPM5,Riemann2D-04}, this 
example is commonly focused on the description of the instability of 
the slip line (marked by the pink dashed box in the first figure of 
Fig. \ref{fig:ex:Riemann2D}). Thus, we have displayed the close-up 
view of this instability in Fig. \ref{fig:ex:Riemann2D-ZoomedIn}. We 
can see that the WENO-JS, WENO-M and WENO-MAIM1 schemes failed to 
resolve the instability of the slip line under current spatial 
resolution. The WENO-IM(2, 0.1) scheme is able to resolve this 
instability, but it is not very noticeable. However, the structure 
of the instability in the solution of the WENO-MAIM3 scheme is 
evidently larger in size and has more waves when compared with that 
of the WENO-IM(2, 0.1) scheme, and this demonstrates the advantage 
of the WENO-MAIM3 scheme that has the lower dissipation and hence 
has a higher resolution in capturing details of the complicated flow 
structures. The solution of the WENO-MAIM4 scheme appears to resolve 
the instability with the waves larger in size but less in number 
than that of the WENO-IM(2, 0.1) scheme.

\begin{example}
\bf{(Double Mach reflection, DMR)}
\rm{We solve the two dimensional double Mach reflection problem 
\cite{interactingBlastWaves-Woodward-Colella} using the considered 
WENO schemes. It is a commonly used test where a vertical shock 
wave moves horizontally into a wedge that is inclined by some angle. 
The computational domain is $[0,4]\times[0,1]$ and the initial 
condition is given by}
\label{ex:DMR}
\end{example}
\begin{equation*}
\big( \rho, u, v, p \big)(x, y, 0) = \left\{
\begin{aligned}
\begin{array}{ll}
(8.0, 8.25\cos \dfrac{\pi}{6},- 8.25\sin \dfrac{\pi}{6}, 116.5)
,         & x < x_{0} + \dfrac{y}{\sqrt{3}}, \\
(1.4, 0.0, 0.0, 1.0),    & x \geq x_{0} + \dfrac{y}{\sqrt{3}}, \\
\end{array}
\end{aligned}
\right.
\label{eq:initial_Euer2D:DMR}
\end{equation*}
where $x_{0} = \frac{1}{6}$. At $x = 0, 4$, the inflow boundary 
condition with the post-shock values as stated above and the outflow 
boundary condition are used respectively. At $y = 0$, the reflective 
boundary condition is applied to the interval $[x_{0}, 4]$, while at 
$(0, x_{0})$, the post-shock values are imposed. The boundary 
condition on the upper boundary $y = 1$ is implemented as follows
\begin{equation*}
\begin{aligned}
\begin{array}{l}
\big( \rho, u, v, p \big)(x, 1, t) = \left\{
\begin{array}{ll}
(8.0, 8.25\cos \dfrac{\pi}{6}, -8.25\sin \dfrac{\pi}{6}, 116.5), & x \in [0, s(t)), \\
(1.4, 0.0, 0.0, 1.0), & x \in [s(t), 4].
\end{array}
\right.
\end{array}
\end{aligned}
\label{eq:boundary:Euler2D:DMR:UB}
\end{equation*}
where $s(t) = x_{0} + \frac{1 + 20t}{\sqrt{3}}$ and it is the 
position of shock wave at time $t$ on the upper boundary. We 
discretize the computational domain using a uniform mesh size of 
$2000\times500$ and the output time is set to be $t=0.2$.

\begin{figure}[ht]
\centering
  \includegraphics[height=0.28\textwidth]
  {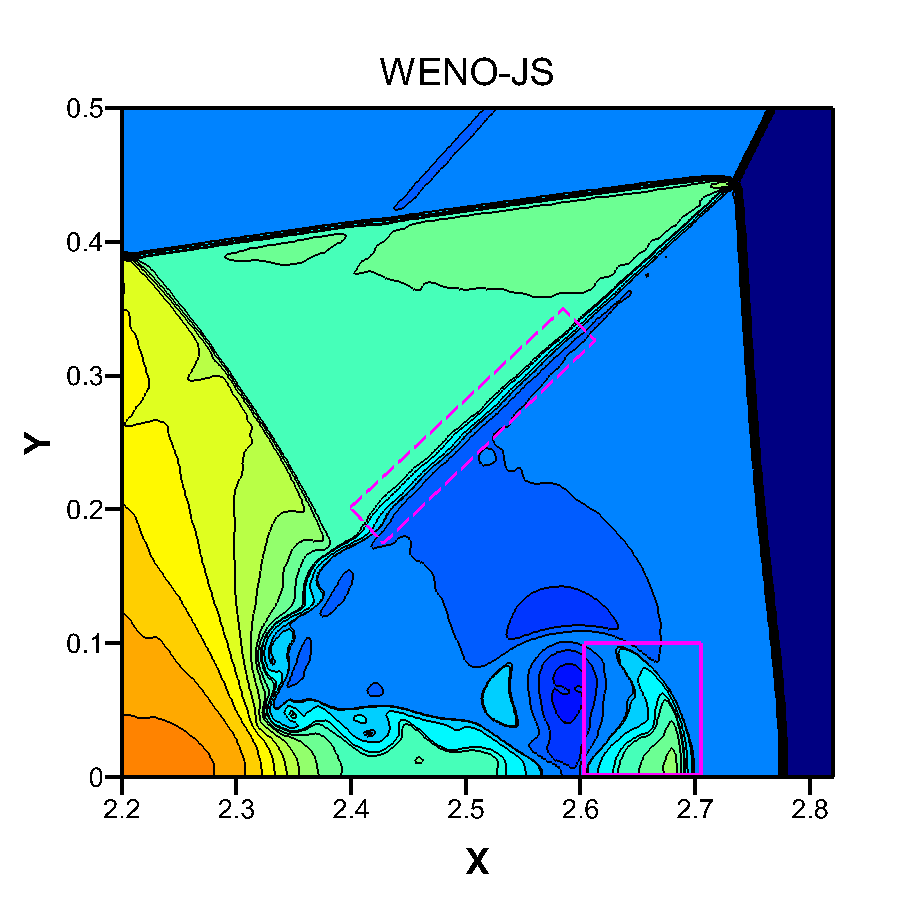}
  \includegraphics[height=0.28\textwidth]
  {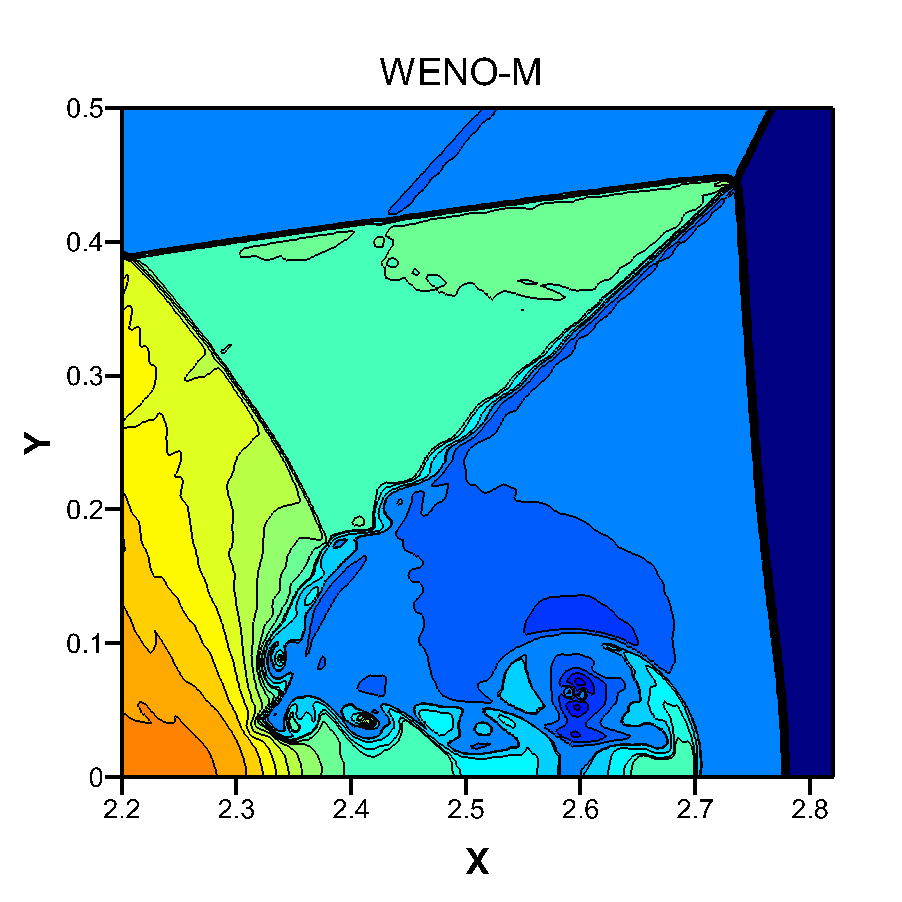}
  \includegraphics[height=0.28\textwidth]
  {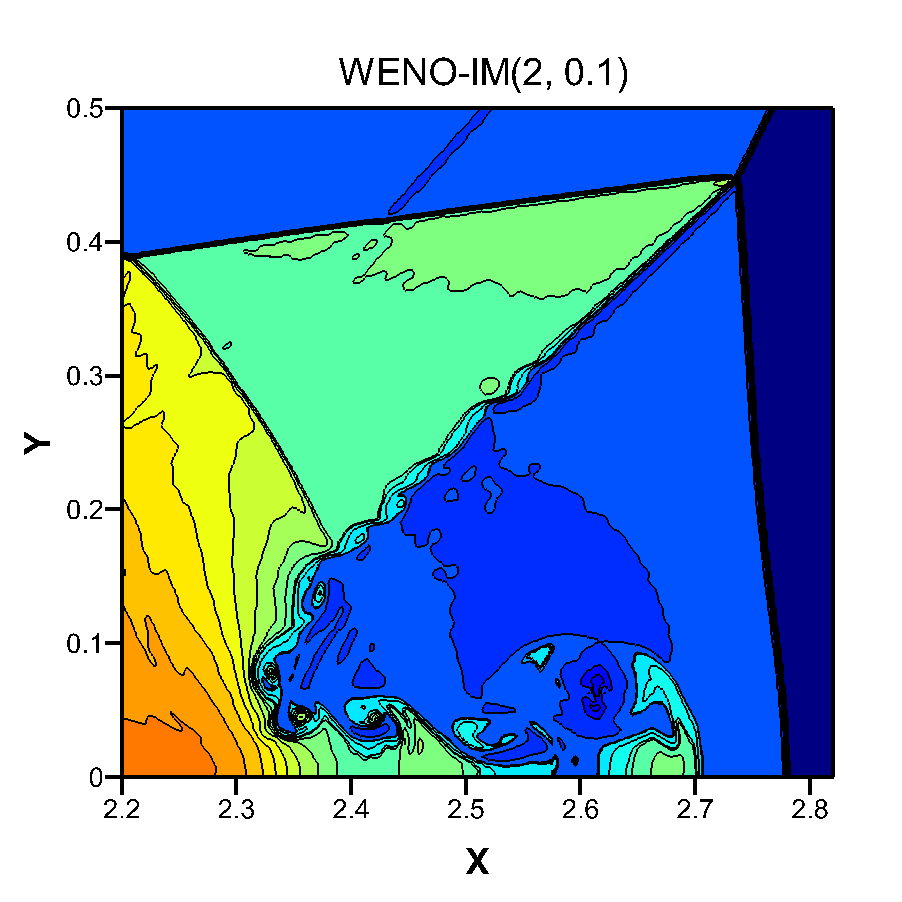}\\
  \includegraphics[height=0.28\textwidth]
  {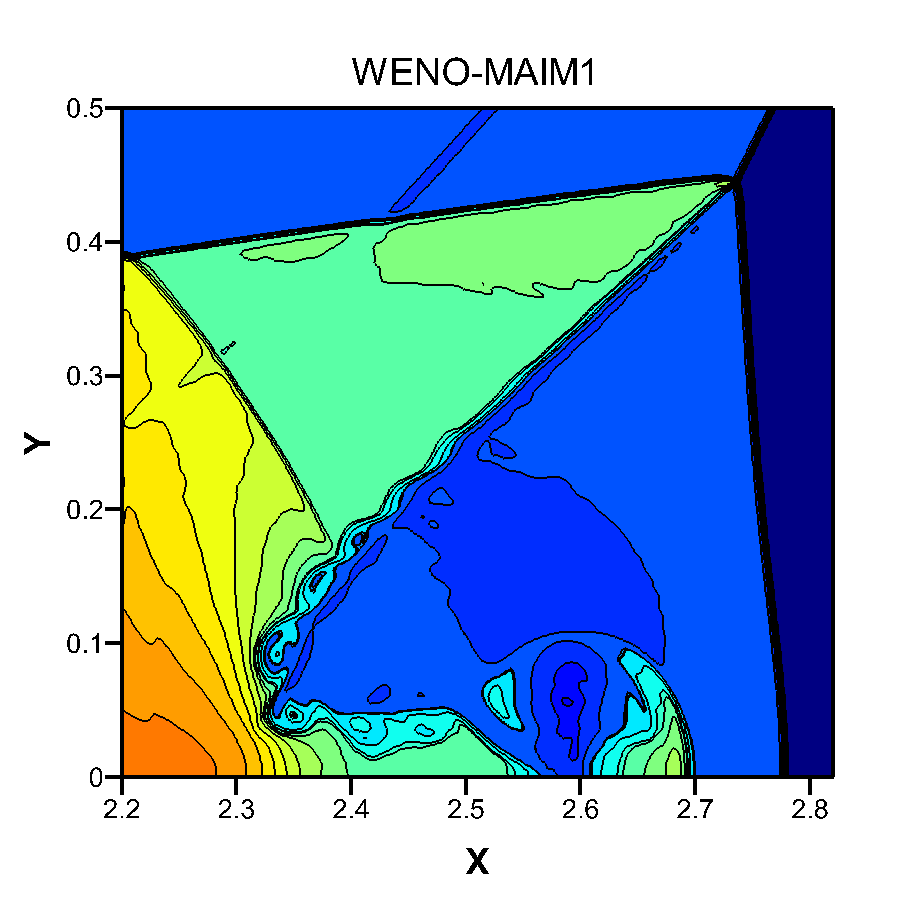}
  \includegraphics[height=0.28\textwidth]
  {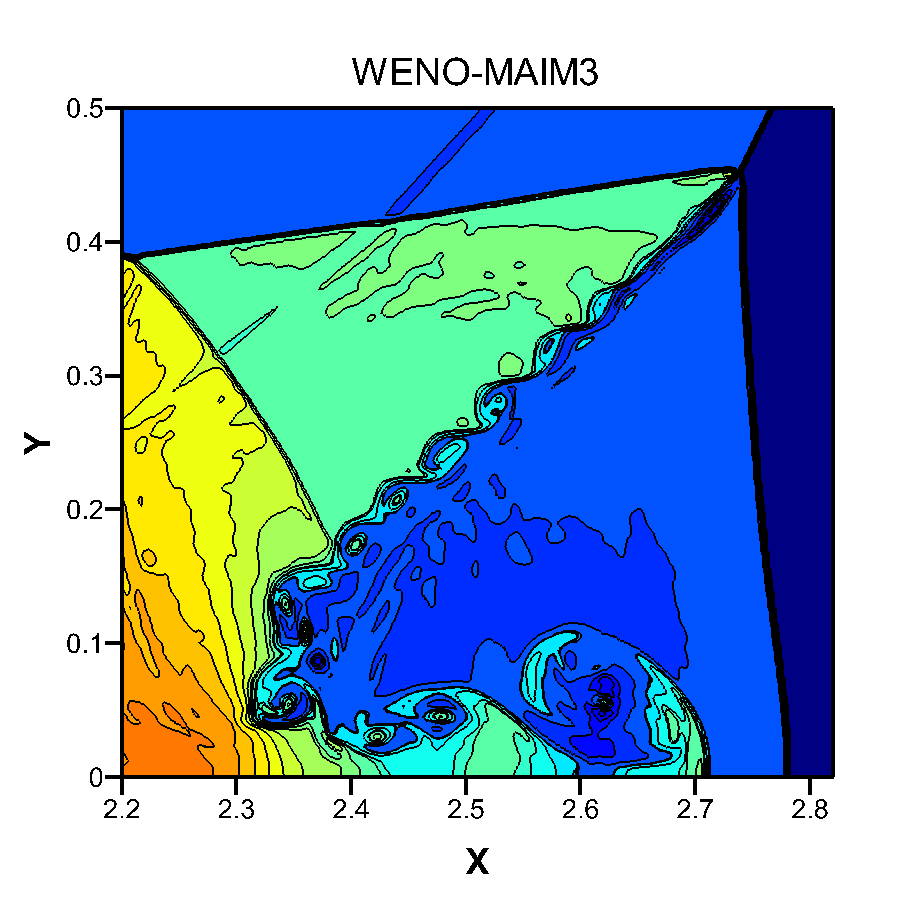}
  \includegraphics[height=0.28\textwidth]
  {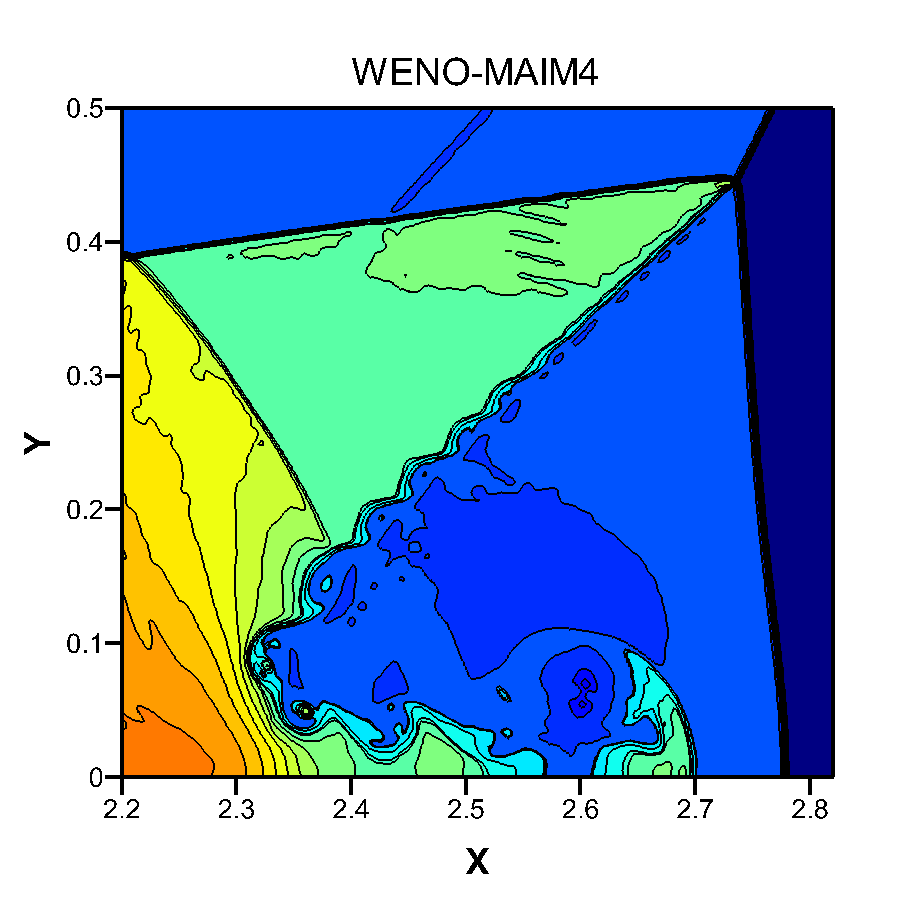}
\caption{Close-up views of performance of the WENO-JS, WENO-M, 
WENO-IM(2, 0.1) and WENO-MAIM$i(i=1,3,4)$ schemes for the DMR 
problem at output time $t = 0.2$ with a uniform mesh size of $2000
\times500$.}
\label{fig:ex:DMR:ZoomedIn}
\end{figure}

As the instability of the slip line and the companion structures 
after the primary reflection shock (marked by the pink dashed box 
and the pink solid box in the first figure of Fig. 
\ref{fig:ex:DMR:ZoomedIn} respectively) are the usual two concerns 
in this problem, we give the zoomed-in contours around the double 
Mach reflection region in Fig. \ref{fig:ex:DMR:ZoomedIn}. All the 
considered schemes can present the main structure of the double Mach 
reflection in general, and they can also successfully capture the 
companion structures behind the lower half of the right moving 
reflection shock. However, one can distinguish the dissipation of 
the various schemes by the number and size of the small vortices 
generated along the slip lines. It is clear that the WENO-MAIM3 
scheme captures most in number and biggest in size of the small 
vortices. We can also see that all the mapped WENO schemes capture 
more in number and bigger in size of the small vortices than the 
WENO-JS scheme. The WENO-MAIM4 and WENO-IM(2, 0.1) schemes appear to 
give slightly better flow features near the slip lines than the 
WENO-MAIM1 and WENO-M schemes.

\begin{example}
\bf{(Forward facing step problem, FFS)} 
\rm{We solve the forward facing step problem first presented by 
Woodward and Colella \cite{interactingBlastWaves-Woodward-Colella}. 
Recently, some important details, like the physical instability and 
roll-up of the vortex sheet emanating from the Mach stem, have been 
successfully captured by various high order schemes 
\cite{ForwardFacingStep-Cockburn_Shu,ADER-WENO-1,ADER-WENO-2, 
WENO-eta,1DEuler_exact}. We will show that our new method is also 
able to successfully capture these important details.

The setup of this problem is as follows: a step with a height of 
$0.2$ length units locates $0.6$ length units from the left-hand end 
of a wind tunnel that is $1$ length unit wide and $3$ length units 
long. The initial condition is specified by}
\label{ex:FFS}
\end{example}
\begin{equation*}
\big( \rho, u, v, p \big)(x, y, 0) = (1.4, 3.0, 0.0, 1.0), \quad
(x,y) \in \Omega.
\label{eq:initial_Euer2D:FFS}
\end{equation*}
where $\Omega = [0,0.6]\times[0,1]\cup[0.6,3]\times[0.2,1]$ is the 
computational domain. The reflective boundary condition is used 
along the walls of the wind tunnel and the step, and the inflow and 
outflow conditions are used at the entrance and the exit of the wind 
tunnel respectively. We discretize the computational domain using a 
uniform mesh size of $900\times300$ and evolve the initial data 
until time $t=4$. 

\begin{figure}[ht]
\centering
  \includegraphics[height=0.193\textwidth]
  {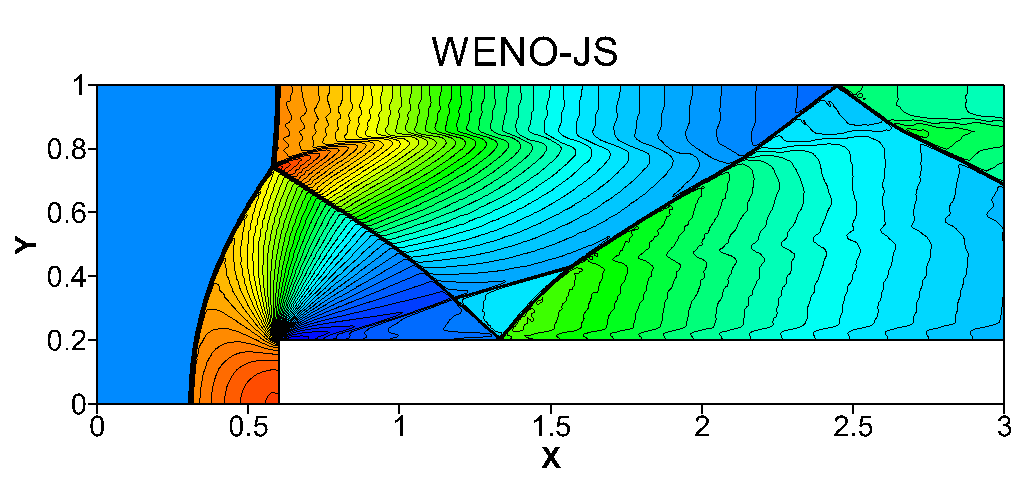}
  \includegraphics[height=0.193\textwidth]
  {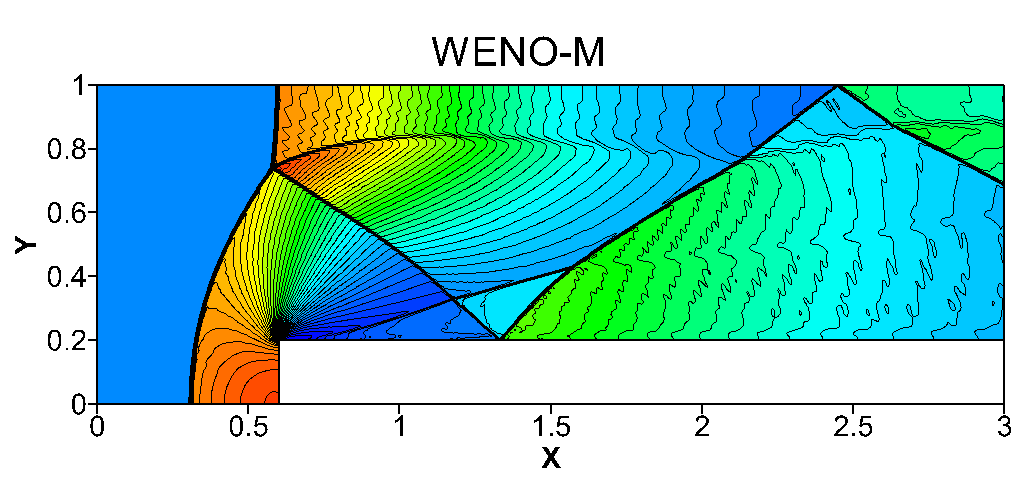}\\
  \includegraphics[height=0.193\textwidth]
  {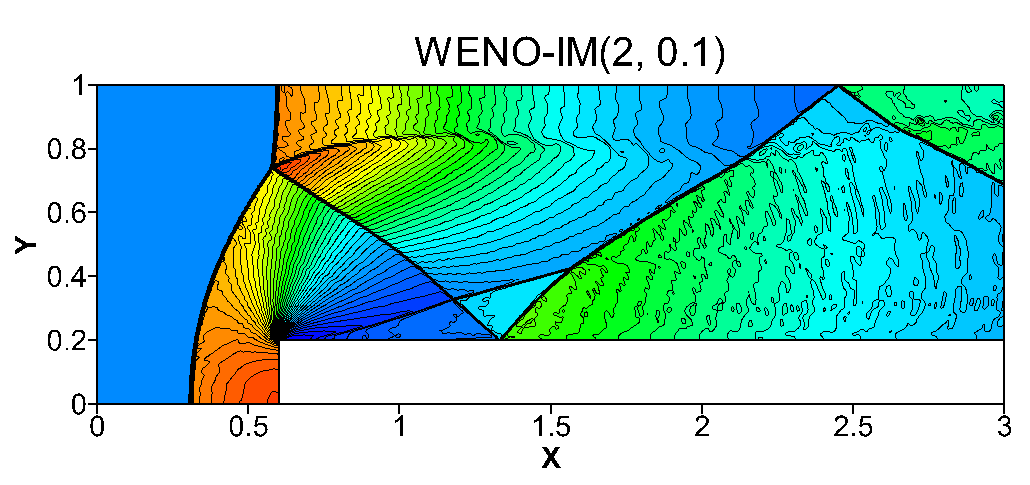}
  \includegraphics[height=0.193\textwidth]
  {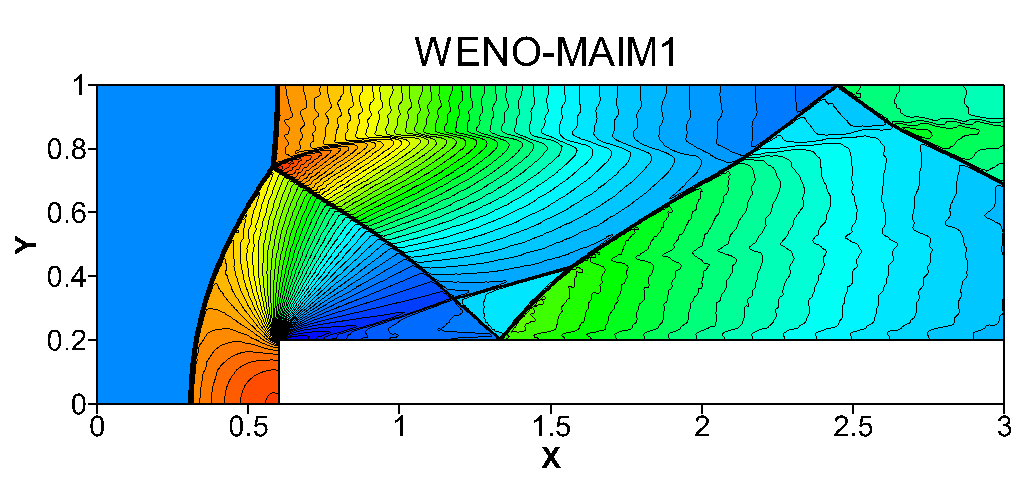}\\
  \includegraphics[height=0.193\textwidth]
  {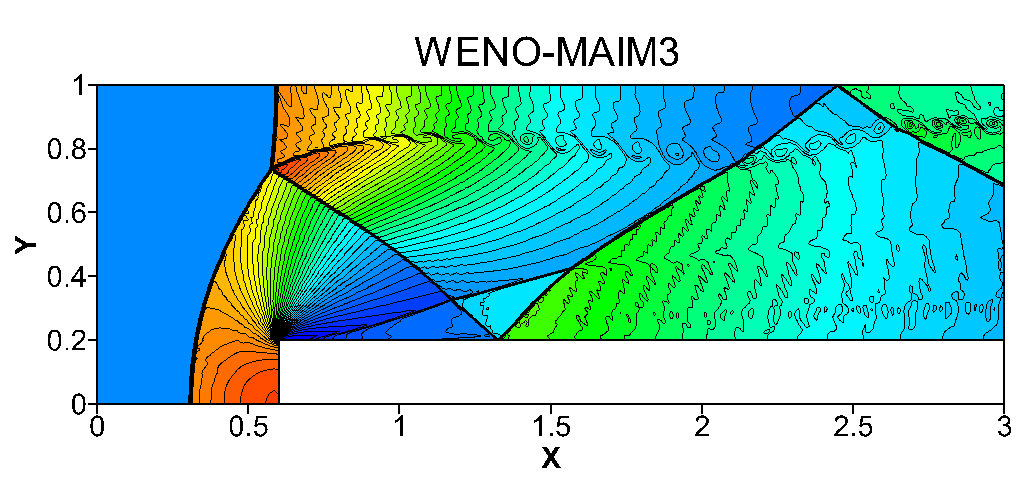}
  \includegraphics[height=0.193\textwidth]
  {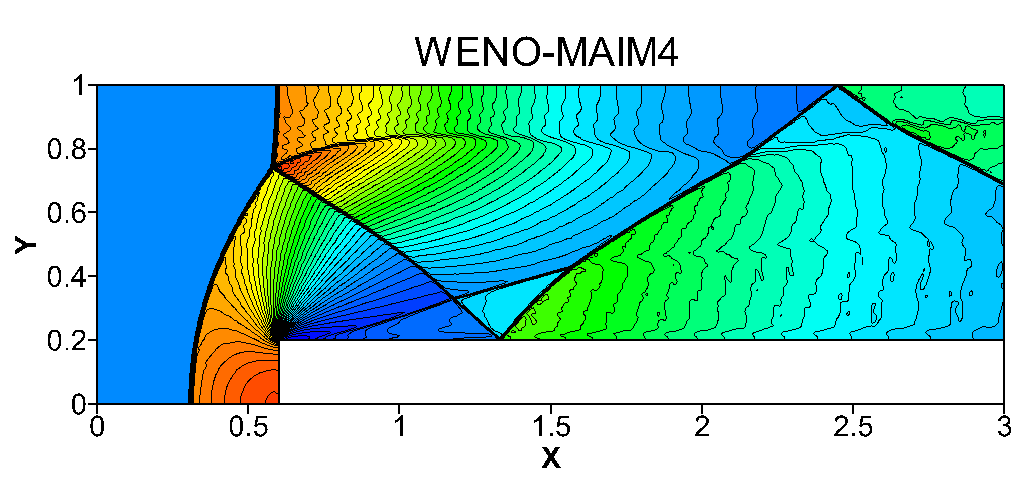}  
\caption{Density plots for the forward facing step problem using 50
density contour lines with range from $0.32$ to $6.5$, computed
using the WENO-JS, WENO-M, WENO-IM(2, 0.1) and WENO-MAIM$i(i=1,3,4)$ 
schemes at output time $t = 4.0$ with a uniform mesh size of $900
\times300$.}
\label{fig:ex:FFS:900x300}
\end{figure}

In Fig. \ref{fig:ex:FFS:900x300}, we have shown the density 
contours computed by considered WENO schemes. We can see that the 
considered schemes capture all the shocks properly with sharp 
profiles. Moreover, the roll-up of the vortex sheet in the solutions 
of the WENO-MAIM3 scheme is most visible, and this demonstrates the 
ability of the WENO-MAIM3 scheme to provide a better resolution in 
solving problems with complicated flow structures. We can see that 
the WENO-IM(2, 0.1) scheme also shows evidence of the vortex sheet’s 
roll-up although it is less evident than that of the WENO-MAIM3 
scheme. Meanwhile, the vortex sheet’s roll-up is not observed in the 
solutions of the WENO-JS, WENO-M, WENO-MAIM1 and WENO-MAIM4 schemes. 

%% file: article_conclusion.tex
\section{Conclusions}
\label{secConclusions} 
We have devised a modified adaptive improved mapped WENO method, 
that is, WENO-MAIMi, by employing a new family of mapping functions. 
The new mapping functions are based on some adaptive control 
functions and a smoothing approximation of the signum function. This 
helps to introduce the adaptive nature and provide a wider selection 
of the parameters in the schemes. In several numerical experiments, 
one-dimensional scalar, as well as one- and two-dimensional system 
cases, is considered to show the performance of our new method. Four 
specified WENO-MAIM$i$ schemes with fine-tuned parameters are used, 
and all these proposed schemes can achieve optimal convergence 
orders even near critical points in smooth regions. When a short 
time simulation is desired, all the four proposed WENO-MAIM$i$ 
schemes are good choices in most hyperbolic conservation 
simulations. In summary,  we recommend to use the WENO-MAIM$3$ 
scheme with parameters as $k=10$, $A=1\mathrm{e-}6$ for pratical 
computation, because it exhibits improved comprehensive performance 
in most hyperbolic conservation simulations. It shows clear 
advantages on approximating near discontinuities, especially for 
one-dimensional linear advection cases with long output time and 
two-dimensional cases with complicated solution structures.


%% file: article_appendix.tex
\section*{Appendix A}
\label{Appendix:recommended_alpha_s}
The objective of this appendix is to present the results of the zero 
points of $q_{s}'(\omega)$, the monotone intervals and maximum 
values of $q_{s}(\omega)$ in $\omega \in (0, 1)$, and the 
recommended values of $\alpha_{s}$ in Property \ref
{property:helperFunctions01} and Property \ref
{property:helperFunctions02}.

\begin{table}[!th]
\begin{myFontSize}
\centering
\caption{The results of the zero points of $q_{s}'(\omega)$, the 
monotone intervals and maximum values of $q_{s}(\omega)$ in $\omega 
\in (0, 1)$, and the recommended values of $\alpha_{s}$ in Property 
\ref{property:helperFunctions01} and Property \ref
{property:helperFunctions02}.}
\label{table:recommended_alpha_s}
\begin{tabular*}{\hsize}
{@{}@{\extracolsep{\fill}}lllllllll@{}}
\hline
\multicolumn{1}{l}{\multirow{2}{*}{$r$}} & \multicolumn{1}{l}{
\multirow{2}{*}{$s$}} & \multicolumn{1}{l}{\multirow{2}{*}{$d_{s}$}} 
& \multicolumn{1}{l}{\multirow{2}{*}{$\omega_{\mathrm{crit}}$}} 
& \multicolumn{3}{c}{$q_{s}'(\omega)$} 
&\multicolumn{1}{l}{\multirow{2}{*}{$q_{\mathrm{max}}(\omega)$}} & 
\multicolumn{1}{l}{\multirow{2}{*}{\makecell[l]{recommended \\ value 
of $\alpha_{s}$}}}
\\
\cline{5-7} 
\multicolumn{1}{l}{} & \multicolumn{1}{l}{} & \multicolumn{1}{l}{} & 
\multicolumn{1}{l}{} 
& \multicolumn{1}{l}{$\omega \in (0, \omega_{\mathrm{crit}})$}  &  {$
\omega = \omega_{\mathrm{crit}}$}  & {$\omega \in (\omega_{\mathrm{
crit}}, 1)$}  
&\multicolumn{1}{l}{} &\multicolumn{1}{l}{}
\\   \hline
$1$    & $0$   & $1$           & -             & $> 0$ &  $> 0$ &  $>
 0$  &  $< 0$          & $0.0001$   \\
$2$    & $0$   & $1/3$         & $0.380873415$ & $> 0$ &  $= 0$ &  $<
 0$  &  $0.092362764$  & $0.0924$ \\
\quad  & $1$   & $2/3$         & $0.619126585$ & $> 0$ &  $= 0$ &  $<
 0$  &  $0.092362764$  & $0.0924$ \\
$3$    & $0$   & $1/10$        & $0.156515839$ & $> 0$ &  $= 0$ &  $<
 0$  &  $0.574903653$  & $0.5750$ \\
\quad  & $1$   & $6/10$        & $0.570057856$ & $> 0$ &  $= 0$ &  $<
 0$  &  $0.033129692$  & $0.0332$ \\
\quad  & $2$   & $3/10$        & $0.354983613$ & $> 0$ &  $= 0$ &  $<
 0$  &  $0.133395837$  & $0.1334$ \\
$4$    & $0$   & $1/35$        & $0.050334141$ & $> 0$ &  $= 0$ &  $<
 0$  &  $0.954367890$  & $0.9544$ \\
\quad  & $1$   & $12/35$       & $0.388078022$ & $> 0$ &  $= 0$ &  $<
 0$  &  $0.082052483$  & $0.0821$ \\
\quad  & $2$   & $18/35$       & $0.509902230$ & $> 0$ &  $= 0$ &  $<
 0$  &  $0.000674949$  & $0.0007$ \\
\quad  & $3$   & $4/35$        & $0.174716514$ & $> 0$ &  $= 0$ &  $<
 0$  &  $0.527277883$  & $0.5273$ \\
$5$    & $0$   & $1/126$       & $0.014533733$ & $> 0$ &  $= 0$ &  $<
 0$  &  $1.289877215$  & $1.2899$ \\
\quad  & $1$   & $10/63$       & $0.226106436$ & $> 0$ &  $= 0$ &  $<
 0$  &  $0.401303129$  & $0.4014$ \\
\quad  & $2$   & $10/21$       & $0.483489850$ & $> 0$ &  $= 0$ &  $<
 0$  &  $0.001874973$  & $0.0019$ \\
\quad  & $3$   & $20/63$       & $0.368684520$ & $> 0$ &  $= 0$ &  $<
 0$  &  $0.110936397$  & $0.1110$ \\
\quad  & $4$   & $5/126$       & $0.068631349$ & $> 0$ &  $= 0$ &  $<
 0$  &  $0.862927630$  & $0.8630$ \\
$6$    & $0$   & $1/462$       & $0.004036768$ & $> 0$ &  $= 0$ &  $<
 0$  &  $1.617654227$  & $1.6177$ \\
\quad  & $1$   & $5/77$        & $0.107726908$ & $> 0$ &  $= 0$ &  $<
 0$  &  $0.716954249$  & $0.7170$ \\
\quad  & $2$   & $25/77$       & $0.374254393$ & $> 0$ &  $= 0$ &  $<
 0$  &  $0.102277562$  & $0.1023$ \\
\quad  & $3$   & $100/231$     & $0.453273191$ & $> 0$ &  $= 0$ &  $<
 0$  &  $0.014901369$  & $0.0150$ \\
\quad  & $4$   & $25/154$      & $0.229973747$ & $> 0$ &  $= 0$ &  $<
 0$  &  $0.392223085$  & $0.3923$ \\
\quad  & $5$   & $1/77$        & $0.023524596$ & $> 0$ &  $= 0$ &  $<
 0$  &  $1.163386235$  & $1.1634$ \\
$7$    & $0$   & $1/1716$      & $0.001099119$ & $> 0$ &  $= 0$ &  $<
 0$  &  $1.945668808$  & $1.9457$ \\
\quad  & $1$   & $7/286$       & $0.043419023$ & $> 0$ &  $= 0$ &  $<
 0$  &  $0.996333645$  & $0.9964$ \\
\quad  & $2$   & $105/572$     & $0.251943992$ & $> 0$ &  $= 0$ &  $<
 0$  &  $0.341617878$  & $0.3417$ \\
\quad  & $3$   & $175/429$     & $0.435602337$ & $> 0$ &  $= 0$ &  $<
 0$  &  $0.028078584$  & $0.0281$ \\
\quad  & $4$   & $175/572$     & $0.359684813$ & $> 0$ &  $= 0$ &  $<
 0$  &  $0.125509465$  & $0.1256$ \\
\quad  & $5$   & $21/286$      & $0.120110409$ & $> 0$ &  $= 0$ &  $<
 0$  &  $0.678155990$  & $0.6782$ \\
\quad  & $6$   & $7/1716$      & $0.007550349$ & $> 0$ &  $= 0$ &  $<
 0$  &  $1.458456912$  & $1.4585$ \\
$8$    & $0$   & $1/6435$      & $0.000295439$ & $> 0$ &  $= 0$ &  $<
 0$  &  $2.275453190$  & $2.2755$ \\
\quad  & $1$   & $56/6435$     & $0.015907656$ & $> 0$ &  $= 0$ &  $<
 0$  &  $1.266360740$  & $1.2664$ \\
\quad  & $2$   & $196/2145$    & $0.145076866$ & $> 0$ &  $= 0$ &  $<
 0$  &  $0.605994603$  & $0.6060$ \\
\quad  & $3$   & $392/1287$    & $0.358612853$ & $> 0$ &  $= 0$ &  $<
 0$  &  $0.127291552$  & $0.1273$ \\
\quad  & $4$   & $490/1287$    & $0.416050440$ & $> 0$ &  $= 0$ &  $<
 0$  &  $0.047165647$  & $0.0472$ \\
\quad  & $5$   & $392/2145$    & $0.251123145$ & $> 0$ &  $= 0$ &  $<
 0$  &  $0.343478122$  & $0.3435$ \\
\quad  & $6$   & $196/6435$    & $0.053489044$ & $> 0$ &  $= 0$ &  $<
 0$  &  $0.936837253$  & $0.9369$ \\
\quad  & $7$   & $8/6435$      & $0.002330913$ & $> 0$ &  $= 0$ &  $<
 0$  &  $1.756424091$  & $1.7565$ \\
$9$    & $0$   & $1/24310$     & $0.000078677$ & $> 0$ &  $= 0$ &  $<
 0$  &  $2.607051525$  & $2.6071$ \\
\quad  & $1$   & $36/12155$    & $0.005504277$ & $> 0$ &  $= 0$ &  $<
 0$  &  $2.607051528$  & $2.6071$ \\
\quad  & $2$   & $504/12155$   & $0.071503078$ & $> 0$ &  $= 0$ &  $<
 0$  &  $0.850403988$  & $0.8505$ \\
\quad  & $3$   & $2352/12155$  & $0.261799125$ & $> 0$ &  $= 0$ &  $<
 0$  &  $0.319476358$  & $0.3195$ \\
\quad  & $4$   & $882/2431$    & $0.402941270$ & $> 0$ &  $= 0$ &  $<
 0$  &  $0.062457658$  & $0.0625$ \\
\quad  & $5$   & $3528/12155$  & $0.347184466$ & $> 0$ &  $= 0$ &  $<
 0$  &  $0.146872875$  & $0.1469$ \\
\quad  & $6$   & $1176/12155$  & $0.152246813$ & $> 0$ &  $= 0$ &  $<
 0$  &  $0.586388730$  & $0.5864$ \\
\quad  & $7$   & $144/12155$   & $0.021508966$ & $> 0$ &  $= 0$ &  $<
 0$  &  $1.187147894$  & $1.1872$ \\
\quad  & $8$   & $9/24310$     & $0.000700372$ & $> 0$ &  $= 0$ &  $<
 0$  &  $2.058883535$  & $2.0589$ \\
\hline
\end{tabular*}
\end{myFontSize}
\end{table}

\section*{Appendix B}
\label{Appendix_Proofsof_convergenceRateofWENO-M}
In this appendix, we present the proof of Lemma
\ref{Lemma_convergenceRateofWENO-M}. \\ \\
\textbf{Proof of Lemma \ref{Lemma_convergenceRateofWENO-M}.} \\

First, by employing the mapping function Eq.(\ref
{mappingFunctionWENO-M}), the mapped weights with $n$ times mapping 
are given by
\begin{equation*}
\omega_{s}^{\mathrm{M}(n)} = \dfrac{\alpha _{s}^{\mathrm{M}(n)}}{\sum
_{l = 0}^{r - 1} \alpha _{l}^{\mathrm{M}(n)}}, \alpha_{s}^{\mathrm{M}
(n)} = \big( g^{\mathrm{M}} \big)_{s}(\omega^{\mathrm{M}(n - 1)}_{s})
, \quad s = 0, 1, \cdots, r - 1,
\end{equation*}
where $\omega^{\mathrm{M}(n - 1)}_{s}$ are the mapped weights with $n
- 1$ times mapping. Clearly, the mapped weights with $0$ times 
mapping are determined by
\begin{equation}
\omega^{\mathrm{M}(0)}_{s} = \omega_{s}^{\mathrm{JS}}, 
\label{Eq2_ProofCorollary_convergenceRateofWENO-M}
\end{equation}
which are computed by Eq.(\ref{nonlinearWeightsWENO-JS}). Similarly,
as $\Big( g^{\mathrm{M}} \Big)_{s}(d_{s}) = d_{s}$ and $
\Big( g^{\mathrm{M}} \Big)'_{s}(d_{s}) = \Big( g^{\mathrm{M}} \Big)''
_{s}(d_{s}) = 0, \big( g^{\mathrm{M}} \big)'''_{s}(d_{s}) \neq 0$ in 
the WENO-M scheme, evaluation at $\omega_{s}^{\mathrm{M}(n - 1)}$ of 
the Taylor series approximations of $\Big( g^{\mathrm{M}} \Big)_{s}(
\omega)$ about $d_{s}$ yields
\begin{equation*}
\begin{aligned}
\alpha_{s}^{\mathrm{M}(n)} &= \big( g^{\mathrm{M}} \big)_{s}(d_{s}) +
 \big( g^{\mathrm{M}} \big)'_{s}(d_{s})\Big( \omega_{s}^{\mathrm{M}(n
 - 1)} - d_{s} \Big) + \dfrac{\big( g^{\mathrm{M}} \big)''_{s}(d_{s})
 }{2!}\Big( \omega_{s}^{\mathrm{M}(n - 1)} - d_{s} \Big)^{2} + \dfrac
 {\big( g^{\mathrm{M}} \big)'''_{s}(d_{s})}{3!}\Big( \omega_{s}^{
 \mathrm{M}(n - 1)} - d_{s} \Big)^{3} + \cdots \\
&= d_{s} + \dfrac{\big( g^{\mathrm{M}} \big)'''_{s}(d_{s})}{3!}\Big( 
\omega_{s}^{\mathrm{M}(n - 1)} - d_{s} \Big)^{3} + \cdots \\
\end{aligned}
\end{equation*}
By employing the fact that $\sum_{s = 0}^{r - 1}d_{s} = 1$, we obtain
\begin{equation*}
\dfrac{1}{\sum_{l = 0}^{r - 1}\alpha _{l}^{\mathrm{M}(n)}} = 1 + O
\Big( \big( \omega_{s}^{\mathrm{M}(n - 1)} - d_{s} \big)^{3} \Big).
\end{equation*}
and then,
\begin{equation}
\omega_{s}^{\mathrm{M}(n)} - d_{s} = O\Big( \big( \omega_{s}^{\mathrm
{M}(n - 1)} - d_{s} \big)^{3} \Big).
\label{Eq5_ProofCorollary_convergenceRateofWENO-M}
\end{equation}
As Eq.(\ref{Eq5_ProofCorollary_convergenceRateofWENO-M}) is a 
recurrence formula, it is easy to obtain
\begin{equation}
\omega_{s}^{\mathrm{M}(n)} - d_{s} = O\Big( \big( \omega_{s}^{\mathrm
{M}(0)} - d_{s} \big)^{3^{n}} \Big).
\label{Eq6_ProofCorollary_convergenceRateofWENO-M}
\end{equation}
From Eq.(\ref{Eq1_lemma_ncpWENO-JS})(\ref
{Eq2_ProofCorollary_convergenceRateofWENO-M})(\ref
{Eq6_ProofCorollary_convergenceRateofWENO-M}), we have
\begin{equation}
\omega_{s}^{\mathrm{M}(n)} - d_{s} = O\Big( \big( \Delta x \big)^{3^{
n}\times (r - 1 - n_{\mathrm{cp}})} \Big).
\label{Eq7_ProofCorollary_convergenceRateofWENO-M}
\end{equation}
Therefore, according to Lemma \ref{lemma_SufficientCondition} and 
Eq.(\ref{Eq7_ProofCorollary_convergenceRateofWENO-M}), the $(2r - 1)$
th-order WENO-M schemes can achieve the optimal order of accuracy, 
namely, $r_{\mathrm{c}} = 2r - 1$, with the requirement
\begin{equation*}
3^{n}\times(r - 1 - n_{\mathrm{cp}}) \geq r,
\end{equation*}
then $n_{\mathrm{cp}}$ is limited to $n_{\mathrm{cp}} \leq \dfrac{3^{
n}-1}{3^{n}}r - 1$, which we rewrite as
\begin{equation*}
n_{\mathrm{cp}} = 0, 1, \cdots, \Bigg\lfloor \dfrac{3^{n}-1}{3^{n}}r 
- 1 \Bigg\rfloor.
\end{equation*}
Now, we consider the convergence order for $n_{\mathrm{cp}} > \dfrac{
3^{n}-1}{3^{n}}r - 1$. As the convergence order is $O(\Delta x^{r - 1
}) \cdot O(\omega_{s}^{\mathrm{M}(n)} - d_{s})$, which can be found 
in the statement of page 549 in \cite{WENO-M} or from the proof of 
Lemma 1 in \cite{WENO-IM}, we obtain
\begin{equation*}
r_{\mathrm{c}} = (3^{n} + 1)(r - 1) - 3^{n} \times n_{\mathrm{cp}}, 
\quad \quad \mathrm{if} \quad n_{\mathrm{cp}} = \Bigg\lfloor \dfrac{3
^{n} - 1}{3^{n}}r - 1 \Bigg\rfloor + 1, \cdots, r - 1.
\end{equation*}
$\hfill\square$ \\ \\


%% file: article.bbl
\begin{thebibliography}{36}
\expandafter\ifx\csname natexlab\endcsname\relax\def\natexlab#1{#1}\fi
\providecommand{\bibinfo}[2]{#2}
\ifx\xfnm\relax \def\xfnm[#1]{\unskip,\space#1}\fi
\bibitem[{Arandiga et~al.(2011)Arandiga, Baeza, Belda and
  Mulet}]{WENO-GlobalAccuracy}
\bibinfo{author}{F.~Arandiga}, \bibinfo{author}{A.~Baeza},
  \bibinfo{author}{A.M. Belda}, \bibinfo{author}{P.~Mulet},
  \bibinfo{title}{Analysis of {WENO} schemes for full and global accuracy},
  \bibinfo{shortjournal}{SIAM J. Numer. Anal.} \bibinfo{volume}{49}
  (\bibinfo{year}{2011}) \bibinfo{pages}{893--915}.
\bibitem[{Balsara et~al.(2009)Balsara, Rumpf, Dumbser and Munz}]{ADER-WENO-1}
\bibinfo{author}{D.~Balsara}, \bibinfo{author}{D.~Rumpf},
  \bibinfo{author}{M.~Dumbser}, \bibinfo{author}{C.D. Munz},
  \bibinfo{title}{Efficient, high accuracy {ADER-WENO} schemes for
  hydrodynamics and divergence-free magnetohydrodynamics},
  \bibinfo{shortjournal}{J. Comput. Phys.} \bibinfo{volume}{228}
  (\bibinfo{year}{2009}) \bibinfo{pages}{2480--2516}.
\bibitem[{Balsara et~al.(2013)Balsara, Rumpf, Dumbser and Munz}]{ADER-WENO-2}
\bibinfo{author}{D.~Balsara}, \bibinfo{author}{D.~Rumpf},
  \bibinfo{author}{M.~Dumbser}, \bibinfo{author}{C.D. Munz},
  \bibinfo{title}{Efficient implementation of {ADER} schemes for {E}uler and
  magnetohydrodynamical flows on structured meshes - {S}peed comparisons with
  {R}unge-{K}utta methods}, \bibinfo{shortjournal}{J. Comput. Phys.}
  \bibinfo{volume}{235} (\bibinfo{year}{2013}) \bibinfo{pages}{934--969}.
\bibitem[{Borges et~al.(2008)Borges, Carmona, Costa and S.}]{WENO-Z}
\bibinfo{author}{R.~Borges}, \bibinfo{author}{M.~Carmona},
  \bibinfo{author}{B.~Costa}, \bibinfo{author}{D.W. S.}, \bibinfo{title}{An
  improved weighted essentially non-oscillatory scheme for hyperbolic
  conservation laws}, \bibinfo{shortjournal}{J. Comput. Phys.}
  \bibinfo{volume}{227} (\bibinfo{year}{2008}) \bibinfo{pages}{3101--3211}.
\bibitem[{Castro et~al.(2011)Castro, Costa and S.}]{WENO-Z01}
\bibinfo{author}{M.~Castro}, \bibinfo{author}{B.~Costa}, \bibinfo{author}{D.W.
  S.}, \bibinfo{title}{High order weighted essentially non-oscillatory {WENO-Z}
  schemes for hyperbolic conservation laws}, \bibinfo{shortjournal}{J. Comput.
  Phys.} \bibinfo{volume}{230} (\bibinfo{year}{2011})
  \bibinfo{pages}{1766--1792}.
\bibitem[{Cockburn and Shu(1998)}]{ForwardFacingStep-Cockburn_Shu}
\bibinfo{author}{B.~Cockburn}, \bibinfo{author}{C.W. Shu}, \bibinfo{title}{The
  {R}unge-{K}utta discontinuous galerkin method for {C}onservation {L}aws {V}},
  \bibinfo{shortjournal}{J. Comput. Phys.} \bibinfo{volume}{141}
  (\bibinfo{year}{1998}) \bibinfo{pages}{199--224}.
\bibitem[{Don and Borges(2013)}]{WENO-Accuracy}
\bibinfo{author}{W.S. Don}, \bibinfo{author}{R.~Borges},
  \bibinfo{title}{Accuracy of the weighted essentially non-oscillatory
  conservative finite difference schemes}, \bibinfo{shortjournal}{J. Comput.
  Phys.} \bibinfo{volume}{250} (\bibinfo{year}{2013})
  \bibinfo{pages}{347--372}.
\bibitem[{Fan et~al.(2014)Fan, Shen and Tian}]{WENO-eta}
\bibinfo{author}{P.~Fan}, \bibinfo{author}{Y.Q. Shen}, \bibinfo{author}{B.L.
  Tian}, \bibinfo{title}{A new smoothness indicator for improving the weighted
  essentially non-oscillatory scheme}, \bibinfo{shortjournal}{J. Comput. Phys.}
  \bibinfo{volume}{269} (\bibinfo{year}{2014}) \bibinfo{pages}{329--354}.
\bibitem[{Feng et~al.(2012)Feng, Hu and Wang}]{WENO-PM}
\bibinfo{author}{H.~Feng}, \bibinfo{author}{F.~Hu}, \bibinfo{author}{R.~Wang},
  \bibinfo{title}{A new mapped weighted essentially non-oscillatory scheme},
  \bibinfo{shortjournal}{J. Sci. Comput.} \bibinfo{volume}{51}
  (\bibinfo{year}{2012}) \bibinfo{pages}{449--473}.
\bibitem[{Feng et~al.(2014)Feng, Huang and Wang}]{WENO-IM}
\bibinfo{author}{H.~Feng}, \bibinfo{author}{C.~Huang},
  \bibinfo{author}{R.~Wang}, \bibinfo{title}{An improved mapped weighted
  essentially non-oscillatory scheme}, \bibinfo{shortjournal}{Appl. Math.
  Comput.} \bibinfo{volume}{232} (\bibinfo{year}{2014})
  \bibinfo{pages}{453--468}.
\bibitem[{Gerolymos et~al.(2009)Gerolymos, Senechal and
  Vallet}]{veryHighOrderWENO}
\bibinfo{author}{G.A. Gerolymos}, \bibinfo{author}{D.~Senechal},
  \bibinfo{author}{I.~Vallet}, \bibinfo{title}{Very-high-order {WENO} schemes},
  \bibinfo{shortjournal}{J. Comput. Phys.} \bibinfo{volume}{228}
  (\bibinfo{year}{2009}) \bibinfo{pages}{8481--8524}.
\bibitem[{Gottlied and Shu(1998)}]{SSPRK1998}
\bibinfo{author}{S.~Gottlied}, \bibinfo{author}{C.W. Shu},
  \bibinfo{title}{Totalvariation diminishing runge-kutta schemes},
  \bibinfo{shortjournal}{Math. Comput.} \bibinfo{volume}{67}
  (\bibinfo{year}{1998}) \bibinfo{pages}{73--85}.
\bibitem[{Gottlied et~al.(2001)Gottlied, Shu and Tadmor}]{SSPRK2001}
\bibinfo{author}{S.~Gottlied}, \bibinfo{author}{C.W. Shu},
  \bibinfo{author}{E.~Tadmor}, \bibinfo{title}{Strong stability-preserving
  high-order time discretization methods}, \bibinfo{shortjournal}{SIAM Rev.}
  \bibinfo{volume}{43} (\bibinfo{year}{2001}) \bibinfo{pages}{89--112}.
\bibitem[{Harten(1987)}]{ENO1987JCP83}
\bibinfo{author}{A.~Harten}, \bibinfo{title}{{ENO} schemes with subcell
  resolution}, \bibinfo{shortjournal}{J. Comput. Phys.} \bibinfo{volume}{83}
  (\bibinfo{year}{1987}) \bibinfo{pages}{148--184}.
\bibitem[{Harten et~al.(1987)Harten, Engquist, Osher and
  Chakravarthy}]{ENO1987JCP71}
\bibinfo{author}{A.~Harten}, \bibinfo{author}{B.~Engquist},
  \bibinfo{author}{S.~Osher}, \bibinfo{author}{S.~Chakravarthy},
  \bibinfo{title}{Uniformly high order essentially non-oscillatory schemes
  {III}}, \bibinfo{shortjournal}{J. Comput. Phys.} \bibinfo{volume}{71}
  (\bibinfo{year}{1987}) \bibinfo{pages}{231--303}.
\bibitem[{Harten et~al.(1986)Harten, Osher, Engquist and
  Chakravarthy}]{ENO1986}
\bibinfo{author}{A.~Harten}, \bibinfo{author}{B.~Osher},
  \bibinfo{author}{S.~Engquist}, \bibinfo{author}{S.~Chakravarthy},
  \bibinfo{title}{Some results on uniformly high order accurate essentially
  non-oscillatory schemes}, \bibinfo{shortjournal}{Appl. Numer. Math.}
  \bibinfo{volume}{2} (\bibinfo{year}{1986}) \bibinfo{pages}{347--377}.
\bibitem[{Harten and Osher(1987)}]{ENO1987V24}
\bibinfo{author}{A.~Harten}, \bibinfo{author}{S.~Osher},
  \bibinfo{title}{Uniformly high order essentially non-oscillatory schemes
  {I}}, \bibinfo{shortjournal}{SIAM J. Numer. Anal.} \bibinfo{volume}{24}
  (\bibinfo{year}{1987}) \bibinfo{pages}{279--309}.
\bibitem[{Henrick et~al.(2005)Henrick, Aslam and Powers}]{WENO-M}
\bibinfo{author}{A.K. Henrick}, \bibinfo{author}{T.D. Aslam},
  \bibinfo{author}{J.M. Powers}, \bibinfo{title}{Mapped weighted essentially
  non-oscillatory schemes: Achieving optimal order near critical points},
  \bibinfo{shortjournal}{J. Comput. Phys.} \bibinfo{volume}{207}
  (\bibinfo{year}{2005}) \bibinfo{pages}{542--567}.
\bibitem[{Hu et~al.(2016)Hu, Wang and X.}]{WENO-Z02}
\bibinfo{author}{F.~Hu}, \bibinfo{author}{R.~Wang}, \bibinfo{author}{C.~X.},
  \bibinfo{title}{A modified fifth-order {WENO-Z} method hyperbolic
  conservation laws}, \bibinfo{shortjournal}{J. Comput. Appl. Math.}
  \bibinfo{volume}{303} (\bibinfo{year}{2016}) \bibinfo{pages}{56--68}.
\bibitem[{Jiang and Shu(1996)}]{WENO-JS}
\bibinfo{author}{G.S. Jiang}, \bibinfo{author}{C.W. Shu},
  \bibinfo{title}{Efficient implementation of weighted {ENO} schemes},
  \bibinfo{shortjournal}{J. Comput. Phys.} \bibinfo{volume}{126}
  (\bibinfo{year}{1996}) \bibinfo{pages}{202--228}.
\bibitem[{LeVeque(2002)}]{LeVeque-FVM}
\bibinfo{author}{R.J. LeVeque}, \bibinfo{title}{Finite {Volume} {Methods} for
  {Hyperbolic} {Problems}}, \bibinfo{publisher}{Cambridge University Press},
  \bibinfo{year}{2002}.
\bibitem[{Liu et~al.(2015)Liu, Liu and Zhang}]{WENO-PPM5}
\bibinfo{author}{Q.~Liu}, \bibinfo{author}{P.~Liu}, \bibinfo{author}{H.~Zhang},
  \bibinfo{title}{Piecewise {Polynomial} {Mapping} {Method} and {Corresponding}
  {WENO} {Scheme} with {Improved} {Resolution}}, \bibinfo{shortjournal}{Commun.
  Comput. Phys.} \bibinfo{volume}{18} (\bibinfo{year}{2015})
  \bibinfo{pages}{1417--1444}.
\bibitem[{Liu et~al.(1994)Liu, Osher and Chan}]{WENO-LiuXD}
\bibinfo{author}{X.D. Liu}, \bibinfo{author}{S.~Osher},
  \bibinfo{author}{T.~Chan}, \bibinfo{title}{Weighted essentially
  non-oscillatory schemes}, \bibinfo{shortjournal}{J. Comput. Phys.}
  \bibinfo{volume}{115} (\bibinfo{year}{1994}) \bibinfo{pages}{200--212}.
\bibitem[{Pirozzoli(2010)}]{Riemann2D-04}
\bibinfo{author}{S.~Pirozzoli}, \bibinfo{title}{Numerical methods for
  high-speed flows}, \bibinfo{shortjournal}{Annu. Rev. Fluid Mech.}
  \bibinfo{volume}{43} (\bibinfo{year}{2010}) \bibinfo{pages}{163--194}.
\bibitem[{Schulz-Rinne(1993)}]{Riemann2D-02}
\bibinfo{author}{C.W. Schulz-Rinne}, \bibinfo{title}{Classification of the
  {R}iemann problem for two-dimensional gas dynamics},
  \bibinfo{shortjournal}{SIAM J. Math. Anal.} \bibinfo{volume}{24}
  (\bibinfo{year}{1993}) \bibinfo{pages}{76--88}.
\bibitem[{Schulz-Rinne et~al.(1993)Schulz-Rinne, Collins and
  Glaz}]{Riemann-2D-01}
\bibinfo{author}{C.W. Schulz-Rinne}, \bibinfo{author}{J.P. Collins},
  \bibinfo{author}{H.M. Glaz}, \bibinfo{title}{Numerical solution of the
  {R}iemann problem for two-dimensional gas dynamics},
  \bibinfo{shortjournal}{SIAM J. Sci. Comput.} \bibinfo{volume}{14}
  (\bibinfo{year}{1993}) \bibinfo{pages}{1394--1414}.
\bibitem[{Shu(1998)}]{WENOoverview}
\bibinfo{author}{C.W. Shu}, \bibinfo{title}{Essentially non-oscillatory and
  weighted essentially non-oscillatory schemes for hyperbolic conservation
  laws}, in: \bibinfo{booktitle}{{Advanced Numerical Approximation of Nonlinear
  Hyperbolic Equations. Lecture Notes in Mathematics}}, volume
  \bibinfo{volume}{1697}, \bibinfo{publisher}{Springer},
  \bibinfo{address}{Berlin}, \bibinfo{year}{1998}, pp.
  \bibinfo{pages}{325--432}.
\bibitem[{Shu and Osher(1988)}]{ENO-Shu1988}
\bibinfo{author}{C.W. Shu}, \bibinfo{author}{S.~Osher},
  \bibinfo{title}{Efficient implementation of essentially non-oscillatory
  shock-capturing schemes}, \bibinfo{shortjournal}{J. Comput. Phys.}
  \bibinfo{volume}{77} (\bibinfo{year}{1988}) \bibinfo{pages}{439--471}.
\bibitem[{Shu and Osher(1989)}]{ENO-Shu1989}
\bibinfo{author}{C.W. Shu}, \bibinfo{author}{S.~Osher},
  \bibinfo{title}{Efficient implementation of essentially non-oscillatory
  shock-capturing schemes {II}}, \bibinfo{shortjournal}{J. Comput. Phys.}
  \bibinfo{volume}{83} (\bibinfo{year}{1989}) \bibinfo{pages}{32--78}.
\bibitem[{Titarev and Toro(2004)}]{Titarev-Toro-1}
\bibinfo{author}{V.~Titarev}, \bibinfo{author}{E.~Toro},
  \bibinfo{title}{Finite-volume {WENO} schemes for three-dimensional
  conservation laws}, \bibinfo{shortjournal}{J. Comput. Phys.}
  \bibinfo{volume}{201} (\bibinfo{year}{2004}) \bibinfo{pages}{238--260}.
\bibitem[{Titarev and Toro(2005)}]{Titarev-Toro-3}
\bibinfo{author}{V.~Titarev}, \bibinfo{author}{E.~Toro}, \bibinfo{title}{{WENO}
  schemes based on upwind and centred {TVD} fluxes},
  \bibinfo{shortjournal}{Comput. Fluids} \bibinfo{volume}{34}
  (\bibinfo{year}{2005}) \bibinfo{pages}{705--720}.
\bibitem[{Toro and Titarev(2005)}]{Titarev-Toro-2}
\bibinfo{author}{E.~Toro}, \bibinfo{author}{V.~Titarev}, \bibinfo{title}{{TVD}
  {Fluxes} for the {High}-{Order} {ADER} {Schemes}}, \bibinfo{shortjournal}{J.
  Sci. Comput.} \bibinfo{volume}{24} (\bibinfo{year}{2005})
  \bibinfo{pages}{285--309}.
\bibitem[{Wang et~al.(2016)Wang, Feng and Huang}]{WENO-RM260}
\bibinfo{author}{R.~Wang}, \bibinfo{author}{H.~Feng},
  \bibinfo{author}{C.~Huang}, \bibinfo{title}{A {New} {Mapped} {Weighted}
  {Essentially} {Non-oscillatory} {Method} {Using} {Rational} {Function}},
  \bibinfo{shortjournal}{J. Sci. Comput.} \bibinfo{volume}{67}
  (\bibinfo{year}{2016}) \bibinfo{pages}{540--580}.
\bibitem[{Woodward and Colella(1984)}]{interactingBlastWaves-Woodward-Colella}
\bibinfo{author}{P.~Woodward}, \bibinfo{author}{P.~Colella},
  \bibinfo{title}{The numerical simulation of two-dimensional fluid flow with
  strong shocks}, \bibinfo{shortjournal}{J. Comput. Phys.} \bibinfo{volume}{54}
  (\bibinfo{year}{1984}) \bibinfo{pages}{115--173}.
\bibitem[{Zhang et~al.(2011)Zhang, Zhang and Shu}]{WENO-accuracy_fv-WENO}
\bibinfo{author}{R.~Zhang}, \bibinfo{author}{M.~Zhang}, \bibinfo{author}{C.W.
  Shu}, \bibinfo{title}{On the {Order} of {Accuracy} and {Numerical}
  {Performance} of {Two} {Classes} of {Finite} {Volume} {WENO} {Schemes}},
  \bibinfo{shortjournal}{Commun. Comput. Phys.} \bibinfo{volume}{9}
  (\bibinfo{year}{2011}) \bibinfo{pages}{807--827}.
\bibitem[{Zhu and Qiu(2016)}]{1DEuler_exact}
\bibinfo{author}{J.~Zhu}, \bibinfo{author}{J.~Qiu}, \bibinfo{title}{A new fifth
  order finite difference weno scheme for solving hyperbolic conservation
  laws}, \bibinfo{shortjournal}{J. Comput. Phys.} \bibinfo{volume}{318}
  (\bibinfo{year}{2016}) \bibinfo{pages}{110--121}.

\end{thebibliography}
